\documentclass[a4paper,12pt,leqno]{amsart}
\usepackage{latexsym}
\usepackage[all]{xy}

\usepackage{amssymb} 
\usepackage{amsmath} 
\usepackage{color}
\usepackage{comment}
\usepackage{mathtools}

\definecolor{gray}{gray}{0.7}
\definecolor{Gray}{gray}{0.3}

\usepackage{amscd}

\numberwithin{equation}{section}

\theoremstyle{break}
 \newtheorem{theorem}{Theorem}[section]
 \newtheorem{proposition}[theorem]{Proposition}
 \newtheorem{corollary}[theorem]{Corollary}
 
 \newtheorem{conjecture}[theorem]{Conjecture}

 \theoremstyle{definition}
 
 \newtheorem{remark}[theorem]{Remark}
 \newtheorem{example}[theorem]{Example}

\allowdisplaybreaks[3]

\def\C{\mathbb C}
\def\Q{\mathbb Q}

\def\R{\mathbb R}
\def\A{\mathcal{A}}
\def\CR{\mathcal{R}}

\DeclareMathOperator{\Der}{Der}
\DeclareMathOperator{\Poin}{Poin}
\DeclareMathOperator{\Hess}{Hess}

\newcommand{\Flags}{Flag}

\def\Sn{S_n}
\def\Tn{T}
\def\q{q}

\begin{document}
  
\title[A survey of Hessenberg varieties]{A survey of recent developments on \\Hessenberg varieties}
\author {Hiraku Abe}
\address{Faculty of Liberal Arts and Sciences, Osaka Prefecture University, 1-1 Gakuen-cho, Naka-ku, Sakai, Osaka 599-8531, Japan}
\email{hirakuabe@globe.ocn.ne.jp}

\author {Tatsuya Horiguchi}
\address{Department of Pure and Applied Mathematics,
Graduate School of Information Science and Technology,
Osaka University, 
1-5, Yamadaoka, Suita, Osaka, 565-0871, Japan}
\email{tatsuya.horiguchi0103@gmail.com}


\begin{abstract}
This article surveys recent developments on Hessenberg varieties, emphasizing some of the rich connections of their cohomology and combinatorics. In particular, we will see how hyperplane arrangements, representations of symmetric groups, and Stanley's chromatic symmetric functions are related to the cohomology rings of Hessenberg varieties. We also include several other topics on Hessenberg varieties to cover recent developments.
\end{abstract}

\maketitle

\setcounter{tocdepth}{1}

\section{Introduction}
\label{sec:introduction}
Hessenberg varieties are subvarieties of the full flag variety which was introduced by F. De Mari, C. Procesi, and M. A. Shayman (\cite{ma-sh,ma-pr-sh}) around 1990. They provide a relatively new research subject, and similarly to Schubert varieties it has been found that geometry, combinatorics, and representation theory interact nicely on Hessenberg varieties. Let $X$ be a complex $n\times n$ matrix considered as a linear map $X:\C^n\rightarrow \C^n$ and $h:\{1,2,\ldots,n\}\rightarrow \{1,2,\ldots,n\}$ a Hessenberg function, i.e.\ a non-decreasing function satisfying $h(j)\geq j$ for $1\leq j\leq n$. The \textbf{Hessenberg variety} (in type $A_{n-1}$) associated to $X$ and $h$ is defined as follows:
\begin{align*}
\Hess(X,h) := \{ V_{\bullet}\in Fl(\C^n) \mid XV_j\subseteq V_{h(j)} \text{ for all } 1\leq j\leq n\},
\end{align*}
where $Fl(\C^n)$ is the flag variety of $\C^n$ consisting of sequences $V_{\bullet}=(V_1\subset V_2\subset \cdots \subset V_n=\C^n)$ of linear subspaces of $\C^n$ with $\dim V_i=i$ for $1\leq i\leq n$. 
Particular examples include the full flag variety itself, Springer fibres, the Peterson variety, and the permutohedral variety.

Over the past 20 years, plentiful developments of Hessenberg varieties has been made. For example, it has been discovered that hyperplane arrangements and representations of symmetric groups appear when we deal with the cohomology rings of Hessenberg varieties. Also, these representations are determined by Stanley's chromatic symmetric functions of certain graphs, and is related with the Stanley-Stembridge conjecture in graph theory. 

This article is a survey of recent developments on Hessenberg varieties, and it is intended to stimulate future research. We keep our explanations concise and include concrete examples so as to make the important ideas accessible, especially for young mathematicians (e.g. graduate students, postdoctoral fellows, and so on). We also include several other topics on Hessenberg varieties to cover recent developments.
For simplicity, we explain most of the results in type $A$, but we make comments for results which hold in arbitrary Lie type. 

\bigskip
\noindent \textbf{Acknowledgements.}  
We offer our immense gratitude to those who helped in preparation of the talks and this article including Takuro Abe, Peter Crooks, Mikiya Masuda, Haozhi Zeng, and all of the organizers of \textit{International Festival in Schubert Calculus} Jianxun Hu, Changzheng Li, and Leonardo C. Mihalcea. We also want to thank the students in Sun Yat-sen University who helped organize the conference, the audience who attended the talks at the conference, and the readers of this article. 
The first author is supported in part by JSPS Grant-in-Aid for Early-Career Scientists: 18K13413. The second author is supported in part by JSPS Grant-in-Aid for JSPS Fellows: 17J04330.

\bigskip

\section{Background and notations}
\label{sec:background}

In this section, we recall some background, and establish some notations for the rest of the document. 

\subsection{Definitions and basic properties}

Let $n$ be a positive integer, and we use the notation $[n]:=\{1,2,\ldots,n\}$ throughout this document.
A function $h: [n] \to [n]$ is a \textbf{Hessenberg function} if it satisfies the following two conditions:
\begin{enumerate}
\item[(i)] $h(1) \leq h(2) \leq \cdots \leq h(n)$,
\item[(ii)] $h(j) \geq j$ for all $j \in [n]$. 
\end{enumerate}
Note that $h(n)=n$ by definition.
We frequently write a Hessenberg function by listing its values in a sequence, i.e. $h=(h(1),h(2),\ldots,h(n))$.
We may identify a Hessenberg function $h$ with a configuration of (shaded) boxes on a square grid of size $n \times n$ which consists of boxes in the $i$-th row and the $j$-th column satisfying $i \leq h(j)$ for $i,j\in[n]$, as we illustrate in the following example.

\begin{example}\label{example:HessenbergFunction}
Let $n=5$. The Hessenberg function $h=(3,3,4,5,5)$ corresponds to the configuration of the shaded boxes drawn in Figure $\ref{pic:stair-shape}$. 
\begin{figure}[h]
\begin{center}
\begin{picture}(75,75)
\put(0,63){\colorbox{gray}}
\put(0,67){\colorbox{gray}}
\put(0,72){\colorbox{gray}}
\put(4,63){\colorbox{gray}}
\put(4,67){\colorbox{gray}}
\put(4,72){\colorbox{gray}}
\put(9,63){\colorbox{gray}}
\put(9,67){\colorbox{gray}}
\put(9,72){\colorbox{gray}}

\put(15,63){\colorbox{gray}}
\put(15,67){\colorbox{gray}}
\put(15,72){\colorbox{gray}}
\put(19,63){\colorbox{gray}}
\put(19,67){\colorbox{gray}}
\put(19,72){\colorbox{gray}}
\put(24,63){\colorbox{gray}}
\put(24,67){\colorbox{gray}}
\put(24,72){\colorbox{gray}}

\put(30,63){\colorbox{gray}}
\put(30,67){\colorbox{gray}}
\put(30,72){\colorbox{gray}}
\put(34,63){\colorbox{gray}}
\put(34,67){\colorbox{gray}}
\put(34,72){\colorbox{gray}}
\put(39,63){\colorbox{gray}}
\put(39,67){\colorbox{gray}}
\put(39,72){\colorbox{gray}}

\put(45,63){\colorbox{gray}}
\put(45,67){\colorbox{gray}}
\put(45,72){\colorbox{gray}}
\put(49,63){\colorbox{gray}}
\put(49,67){\colorbox{gray}}
\put(49,72){\colorbox{gray}}
\put(54,63){\colorbox{gray}}
\put(54,67){\colorbox{gray}}
\put(54,72){\colorbox{gray}}

\put(60,63){\colorbox{gray}}
\put(60,67){\colorbox{gray}}
\put(60,72){\colorbox{gray}}
\put(64,63){\colorbox{gray}}
\put(64,67){\colorbox{gray}}
\put(64,72){\colorbox{gray}}
\put(69,63){\colorbox{gray}}
\put(69,67){\colorbox{gray}}
\put(69,72){\colorbox{gray}}

\put(0,48){\colorbox{gray}}
\put(0,52){\colorbox{gray}}
\put(0,57){\colorbox{gray}}
\put(4,48){\colorbox{gray}}
\put(4,52){\colorbox{gray}}
\put(4,57){\colorbox{gray}}
\put(9,48){\colorbox{gray}}
\put(9,52){\colorbox{gray}}
\put(9,57){\colorbox{gray}}

\put(15,48){\colorbox{gray}}
\put(15,52){\colorbox{gray}}
\put(15,57){\colorbox{gray}}
\put(19,48){\colorbox{gray}}
\put(19,52){\colorbox{gray}}
\put(19,57){\colorbox{gray}}
\put(24,48){\colorbox{gray}}
\put(24,52){\colorbox{gray}}
\put(24,57){\colorbox{gray}}

\put(30,48){\colorbox{gray}}
\put(30,52){\colorbox{gray}}
\put(30,57){\colorbox{gray}}
\put(34,48){\colorbox{gray}}
\put(34,52){\colorbox{gray}}
\put(34,57){\colorbox{gray}}
\put(39,48){\colorbox{gray}}
\put(39,52){\colorbox{gray}}
\put(39,57){\colorbox{gray}}

\put(45,48){\colorbox{gray}}
\put(45,52){\colorbox{gray}}
\put(45,57){\colorbox{gray}}
\put(49,48){\colorbox{gray}}
\put(49,52){\colorbox{gray}}
\put(49,57){\colorbox{gray}}
\put(54,48){\colorbox{gray}}
\put(54,52){\colorbox{gray}}
\put(54,57){\colorbox{gray}}

\put(60,48){\colorbox{gray}}
\put(60,52){\colorbox{gray}}
\put(60,57){\colorbox{gray}}
\put(64,48){\colorbox{gray}}
\put(64,52){\colorbox{gray}}
\put(64,57){\colorbox{gray}}
\put(69,48){\colorbox{gray}}
\put(69,52){\colorbox{gray}}
\put(69,57){\colorbox{gray}}

\put(0,33){\colorbox{gray}}
\put(0,37){\colorbox{gray}}
\put(0,42){\colorbox{gray}}
\put(4,33){\colorbox{gray}}
\put(4,37){\colorbox{gray}}
\put(4,42){\colorbox{gray}}
\put(9,33){\colorbox{gray}}
\put(9,37){\colorbox{gray}}
\put(9,42){\colorbox{gray}}

\put(15,33){\colorbox{gray}}
\put(15,37){\colorbox{gray}}
\put(15,42){\colorbox{gray}}
\put(19,33){\colorbox{gray}}
\put(19,37){\colorbox{gray}}
\put(19,42){\colorbox{gray}}
\put(24,33){\colorbox{gray}}
\put(24,37){\colorbox{gray}}
\put(24,42){\colorbox{gray}}

\put(30,33){\colorbox{gray}}
\put(30,37){\colorbox{gray}}
\put(30,42){\colorbox{gray}}
\put(34,33){\colorbox{gray}}
\put(34,37){\colorbox{gray}}
\put(34,42){\colorbox{gray}}
\put(39,33){\colorbox{gray}}
\put(39,37){\colorbox{gray}}
\put(39,42){\colorbox{gray}}

\put(45,33){\colorbox{gray}}
\put(45,37){\colorbox{gray}}
\put(45,42){\colorbox{gray}}
\put(49,33){\colorbox{gray}}
\put(49,37){\colorbox{gray}}
\put(49,42){\colorbox{gray}}
\put(54,33){\colorbox{gray}}
\put(54,37){\colorbox{gray}}
\put(54,42){\colorbox{gray}}

\put(60,33){\colorbox{gray}}
\put(60,37){\colorbox{gray}}
\put(60,42){\colorbox{gray}}
\put(64,33){\colorbox{gray}}
\put(64,37){\colorbox{gray}}
\put(64,42){\colorbox{gray}}
\put(69,33){\colorbox{gray}}
\put(69,37){\colorbox{gray}}
\put(69,42){\colorbox{gray}}

%
%
\put(30,18){\colorbox{gray}}
\put(30,22){\colorbox{gray}}
\put(30,27){\colorbox{gray}}
\put(34,18){\colorbox{gray}}
\put(34,22){\colorbox{gray}}
\put(34,27){\colorbox{gray}}
\put(39,18){\colorbox{gray}}
\put(39,22){\colorbox{gray}}
\put(39,27){\colorbox{gray}}

\put(45,18){\colorbox{gray}}
\put(45,22){\colorbox{gray}}
\put(45,27){\colorbox{gray}}
\put(49,18){\colorbox{gray}}
\put(49,22){\colorbox{gray}}
\put(49,27){\colorbox{gray}}
\put(54,18){\colorbox{gray}}
\put(54,22){\colorbox{gray}}
\put(54,27){\colorbox{gray}}

\put(60,18){\colorbox{gray}}
\put(60,22){\colorbox{gray}}
\put(60,27){\colorbox{gray}}
\put(64,18){\colorbox{gray}}
\put(64,22){\colorbox{gray}}
\put(64,27){\colorbox{gray}}
\put(69,18){\colorbox{gray}}
\put(69,22){\colorbox{gray}}
\put(69,27){\colorbox{gray}}

\put(45,3){\colorbox{gray}}
\put(45,7){\colorbox{gray}}
\put(45,12){\colorbox{gray}}
\put(49,3){\colorbox{gray}}
\put(49,7){\colorbox{gray}}
\put(49,12){\colorbox{gray}}
\put(54,3){\colorbox{gray}}
\put(54,7){\colorbox{gray}}
\put(54,12){\colorbox{gray}}

\put(60,3){\colorbox{gray}}
\put(60,7){\colorbox{gray}}
\put(60,12){\colorbox{gray}}
\put(64,3){\colorbox{gray}}
\put(64,7){\colorbox{gray}}
\put(64,12){\colorbox{gray}}
\put(69,3){\colorbox{gray}}
\put(69,7){\colorbox{gray}}
\put(69,12){\colorbox{gray}}

\put(0,0){\framebox(15,15)}
\put(15,0){\framebox(15,15)}
\put(30,0){\framebox(15,15)}
\put(45,0){\framebox(15,15)}
\put(60,0){\framebox(15,15)}
\put(0,15){\framebox(15,15)}
\put(15,15){\framebox(15,15)}
\put(30,15){\framebox(15,15)}
\put(45,15){\framebox(15,15)}
\put(60,15){\framebox(15,15)}
\put(0,30){\framebox(15,15)}
\put(15,30){\framebox(15,15)}
\put(30,30){\framebox(15,15)}
\put(45,30){\framebox(15,15)}
\put(60,30){\framebox(15,15)}
\put(0,45){\framebox(15,15)}
\put(15,45){\framebox(15,15)}
\put(30,45){\framebox(15,15)}
\put(45,45){\framebox(15,15)}
\put(60,45){\framebox(15,15)}
\put(0,60){\framebox(15,15)}
\put(15,60){\framebox(15,15)}
\put(30,60){\framebox(15,15)}
\put(45,60){\framebox(15,15)}
\put(60,60){\framebox(15,15)}
\end{picture}
\end{center}
\caption{The configuration corresponding to $h=(3,3,4,5,5)$.}
\label{pic:stair-shape}
\end{figure}
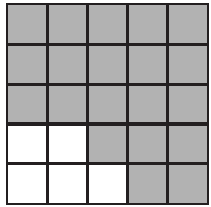
\end{example}

In particular, this identification implies that the set of Hessenberg functions and the set of Dyck paths are in one-to-one correspondence. 
That is, the number of Hessenberg functions is the Catalan number:
$$
\#\{h:[n] \to [n]: \ \mbox{Hessenberg functions} \}=\frac{1}{n+1}\binom{2n}{n}.
$$
There is a natural partial order $\subset$ on Hessenberg functions defined as follows. For any two Hessenberg functions $h$ and $h'$, we define $h \subset h'$ by
\begin{equation*}
h \subset h' \iff h(j) \leq h'(j) \ \ \ \forall j\in[n].
\end{equation*}
We use the symbol $\subset$ for this order since it corresponds to the inclusion of the configurations of boxes under the above visualization of Hessenberg functions. 

The (full) flag variety $Fl(\C^n)$ of $\C^n$ is the collection of nested linear subspaces $V_\bullet= (V_1 \subset V_2 \subset \cdots \subset V_n = \C^n)$ with $\dim V_i=i$ for $i\in[n]$. For an $n\times n$ matrix $X$ considered as a linear map $X:\C^n \to \C^n$ and a Hessenberg function $h: [n] \to [n]$, the \textbf{Hessenberg variety}\footnote{For the origin of the name \textit{Hessenberg varieties}, see \cite{ma-sh}.} (in type $A_{n-1}$) associated with $X$ and $h$ is defined as
\begin{equation} \label{eq:DefinitionHessenbergVariety}
\Hess(X,h)=\{V_\bullet \in Fl(\C^n) \mid XV_j \subset V_{h(j)} \ \mbox{for all} \ j\in[n] \}
\end{equation}
(\cite{ma-sh, ma-pr-sh}).
If $X$ is the zero matrix or $h=(n,n,\ldots,n)$, 
then it is clear that $\Hess(X,h)=Fl(\C^n)$ is the flag variety itself from the definition \eqref{eq:DefinitionHessenbergVariety}. If $X$ is nilpotent and $h=id=(1,2,\ldots,n)$, then $\Hess(X,h)$ is called a \textbf{Springer fiber} which plays an important role in the geometric representation theory of the symmetric group $S_n$ (\cite{Spr1}, \cite{Spr2}). 

\begin{remark}\label{rem: def of Hess}
The definition \eqref{eq:DefinitionHessenbergVariety} can be rephrased in terms of the adjoint representation of $\mbox{GL}(n,\C)$. See \cite{ma-pr-sh} for the definition in arbitrary Lie type.
Also, Goresky-Kottwitz-MacPherson (\cite[Section 2]{G-Ko-Ma06}) considered more general Hessenberg varieties which are defined for arbitrary representations of reductive algebraic groups (cf. Chen-Vilonen-Xue \cite[Section 2]{ch-vi-xu15-2}).
\end{remark}

As a general picture of Hessenberg varieties, we remark the following two properties. 
Suppose that one considers Hessenberg varieties for a fixed matrix $X$. Then, Hessenberg varieties preserves the inclusions:
\begin{equation} \label{eq:inclusion-preserving}
h \subset h' \ \Rightarrow \ \Hess(X,h) \subset \Hess(X,h').
\end{equation}
Also, if $g \in\mbox{GL}(n,\C)$, then we have an isomorphism 
\begin{equation*} 
\Hess(X,h) \cong \Hess(gXg^{-1},h)
\end{equation*}
by sending $V_\bullet$ to $gV_\bullet$. 
This implies that we may assume that $X$ is in a Jordan canonical form.

J. Tymoczko lays the foundation for the study of Hessenberg varieties as follows: 

\begin{theorem}[\cite{ty}] \label{theorem:Tymoczko}
Every Hessenberg variety $\Hess(X,h)$ is paved by affines\footnote{A paving by affines means a ``(complex) cellular decomposition'' in algebraic geometry. See \cite[Definition~2.1]{ty} for the details.}.
In particular, the integral cohomology group of $\Hess(X,h)$ is torsion-free, and the odd-degree cohomology groups vanish. 
\end{theorem}

\begin{remark}
This generalizes the work of Spaltenstein \cite{Spa} for the Springer fibers and De Mari-Procesi-Shayman \cite{ma-pr-sh} for the regular semisimple Hessenberg varieties (which we will define below). 
For generalizations to arbitrary Lie type, see Precup \cite{precup13a} (cf.\ Tymoczko \cite{Ty2}, De Mari-Procesi-Shayman \cite{ma-pr-sh}, 
and De Concini-Lusztig-Procesi \cite{De Concini-Lusztig-Procesi}). 
\end{remark}

By Theorem~\ref{theorem:Tymoczko} we may denote the Poincar\'e polynomial of $\Hess(X,h)$ by
\begin{equation} \label{eq:PoincarePolynomial}
\Poin(\Hess(X,h),\q):=\sum_{i=0}^{m} \dim H^{2i}(\Hess(X,h);\Q) \, \q^i
\end{equation}
where $m:=\dim_\C \Hess(X,h)$ and the variable $q$ stands for the grading with $\text{deg}(q)=2$.\\

In this survey, we will focus on so-called regular nilpotent Hessenberg varieties and regular semisimple Hessenberg varieties which are particularly well-studied. As we will see in Section \ref{sect: cohomology} and Section \ref{sect: combinatorics}, their cohomology rings has an interesting relation between each other, and these are related with other research areas such as hyperplane arrangements, representation theory, and graph theory.

\smallskip

\subsection{Regular nilpotent Hessenberg varieties}
\label{subsec:Regular nilpotent Hessenberg varieties}

Let $N$ a regular nilpotent matrix of size $n\times n$, i.e.\ a nilpotent matrix with a single Jordan block. In Jordan canonical form, it is given by
\begin{equation*} 
N = 
\begin{pmatrix}
0 & 1 & &  & \\
 & 0 &  1  & \\
 &    & \ddots & \ddots & \\
 &    &   & 0 & 1 \\
 &    &  & & 0 \\ 
\end{pmatrix}.
\end{equation*}
For a Hessenberg function $h:[n]\rightarrow [n]$,  $\Hess(N,h)$ is called a \textbf{regular nilpotent Hessenberg variety}. 
We have the following two main examples for this class of Hessenberg varieties. 
For $h=(n,n,\ldots,n)$, $\Hess(N,h)$ is the flag variety $Fl(\C^n)$ itself.
For $h=(2,3,4,\ldots,n,n)$, i.e.\ $h(j)=j+1$ for $1\leq j <n$, 
$\Hess(N,h)$ is called the \textbf{Peterson variety} which is related to the quantum cohomology of partial flag varieties (cf. \cite{Ko}, \cite{R}, \cite{Balibanu}).
Here, $h=(n,n,\ldots,n)$ is the maximum Hessenberg function, and we may think that $h=(2,3,4,\ldots,n,n)$ gives the minimum one in the following sense. 

If we have $h(j)=j$ for some $j<n$, then the Hessenberg function $h$ can be decomposed into two Hessenberg functions $h^{(1)}$ and $h^{(2)}$ of smaller sizes defined as follows: 
\begin{align*}
h^{(1)}:=&(h(1),h(2),\ldots,h(j)), \\ 
h^{(2)}:=&(h(j+1)-j, h(j+2)-j, \ldots, h(n)-j). 
\end{align*}
\vspace{-15pt}
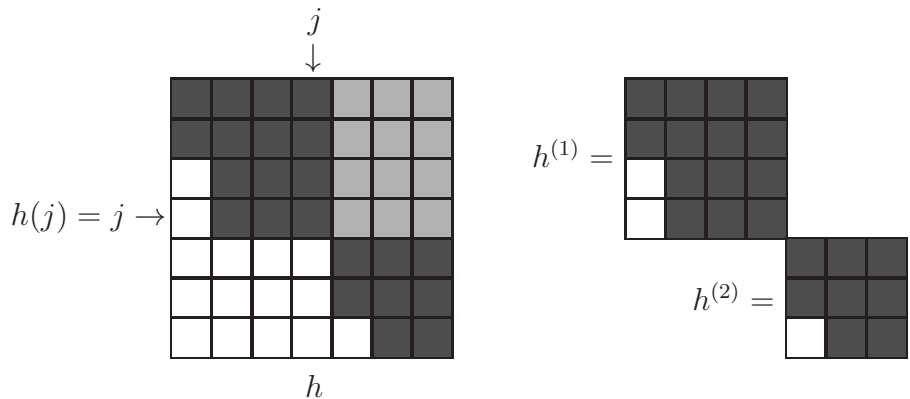
\begin{figure}[h]
\begin{center}
\begin{picture}(200,130)
\put(-10,93){\colorbox{gray}}
\put(-10,97){\colorbox{gray}}
\put(-10,102){\colorbox{gray}}
\put(-6,93){\colorbox{gray}}
\put(-6,97){\colorbox{gray}}
\put(-6,102){\colorbox{gray}}
\put(-1,93){\colorbox{gray}}
\put(-1,97){\colorbox{gray}}
\put(-1,102){\colorbox{gray}}

\put(5,93){\colorbox{gray}}
\put(5,97){\colorbox{gray}}
\put(5,102){\colorbox{gray}}
\put(9,93){\colorbox{gray}}
\put(9,97){\colorbox{gray}}
\put(9,102){\colorbox{gray}}
\put(14,93){\colorbox{gray}}
\put(14,97){\colorbox{gray}}
\put(14,102){\colorbox{gray}}

\put(20,93){\colorbox{gray}}
\put(20,97){\colorbox{gray}}
\put(20,102){\colorbox{gray}}
\put(24,93){\colorbox{gray}}
\put(24,97){\colorbox{gray}}
\put(24,102){\colorbox{gray}}
\put(29,93){\colorbox{gray}}
\put(29,97){\colorbox{gray}}
\put(29,102){\colorbox{gray}}

\put(35,93){\colorbox{gray}}
\put(35,97){\colorbox{gray}}
\put(35,102){\colorbox{gray}}
\put(39,93){\colorbox{gray}}
\put(39,97){\colorbox{gray}}
\put(39,102){\colorbox{gray}}
\put(44,93){\colorbox{gray}}
\put(44,97){\colorbox{gray}}
\put(44,102){\colorbox{gray}}

\put(50,93){\colorbox{gray}}
\put(50,97){\colorbox{gray}}
\put(50,102){\colorbox{gray}}
\put(54,93){\colorbox{gray}}
\put(54,97){\colorbox{gray}}
\put(54,102){\colorbox{gray}}
\put(59,93){\colorbox{gray}}
\put(59,97){\colorbox{gray}}
\put(59,102){\colorbox{gray}}

\put(65,93){\colorbox{gray}}
\put(65,97){\colorbox{gray}}
\put(65,102){\colorbox{gray}}
\put(69,93){\colorbox{gray}}
\put(69,97){\colorbox{gray}}
\put(69,102){\colorbox{gray}}
\put(74,93){\colorbox{gray}}
\put(74,97){\colorbox{gray}}
\put(74,102){\colorbox{gray}}

\put(80,93){\colorbox{gray}}
\put(80,97){\colorbox{gray}}
\put(80,102){\colorbox{gray}}
\put(84,93){\colorbox{gray}}
\put(84,97){\colorbox{gray}}
\put(84,102){\colorbox{gray}}
\put(89,93){\colorbox{gray}}
\put(89,97){\colorbox{gray}}
\put(89,102){\colorbox{gray}}

\put(-10,78){\colorbox{gray}}
\put(-10,82){\colorbox{gray}}
\put(-10,87){\colorbox{gray}}
\put(-6,78){\colorbox{gray}}
\put(-6,82){\colorbox{gray}}
\put(-6,87){\colorbox{gray}}
\put(-1,78){\colorbox{gray}}
\put(-1,82){\colorbox{gray}}
\put(-1,87){\colorbox{gray}}

\put(5,78){\colorbox{gray}}
\put(5,82){\colorbox{gray}}
\put(5,87){\colorbox{gray}}
\put(9,78){\colorbox{gray}}
\put(9,82){\colorbox{gray}}
\put(9,87){\colorbox{gray}}
\put(14,78){\colorbox{gray}}
\put(14,82){\colorbox{gray}}
\put(14,87){\colorbox{gray}}

\put(20,78){\colorbox{gray}}
\put(20,82){\colorbox{gray}}
\put(20,87){\colorbox{gray}}
\put(24,78){\colorbox{gray}}
\put(24,82){\colorbox{gray}}
\put(24,87){\colorbox{gray}}
\put(29,78){\colorbox{gray}}
\put(29,82){\colorbox{gray}}
\put(29,87){\colorbox{gray}}

\put(35,78){\colorbox{gray}}
\put(35,82){\colorbox{gray}}
\put(35,87){\colorbox{gray}}
\put(39,78){\colorbox{gray}}
\put(39,82){\colorbox{gray}}
\put(39,87){\colorbox{gray}}
\put(44,78){\colorbox{gray}}
\put(44,82){\colorbox{gray}}
\put(44,87){\colorbox{gray}}

\put(50,78){\colorbox{gray}}
\put(50,82){\colorbox{gray}}
\put(50,87){\colorbox{gray}}
\put(54,78){\colorbox{gray}}
\put(54,82){\colorbox{gray}}
\put(54,87){\colorbox{gray}}
\put(59,78){\colorbox{gray}}
\put(59,82){\colorbox{gray}}
\put(59,87){\colorbox{gray}}

\put(65,78){\colorbox{gray}}
\put(65,82){\colorbox{gray}}
\put(65,87){\colorbox{gray}}
\put(69,78){\colorbox{gray}}
\put(69,82){\colorbox{gray}}
\put(69,87){\colorbox{gray}}
\put(74,78){\colorbox{gray}}
\put(74,82){\colorbox{gray}}
\put(74,87){\colorbox{gray}}

\put(80,78){\colorbox{gray}}
\put(80,82){\colorbox{gray}}
\put(80,87){\colorbox{gray}}
\put(84,78){\colorbox{gray}}
\put(84,82){\colorbox{gray}}
\put(84,87){\colorbox{gray}}
\put(89,78){\colorbox{gray}}
\put(89,82){\colorbox{gray}}
\put(89,87){\colorbox{gray}}

%
\put(5,63){\colorbox{gray}}
\put(5,67){\colorbox{gray}}
\put(5,72){\colorbox{gray}}
\put(9,63){\colorbox{gray}}
\put(9,67){\colorbox{gray}}
\put(9,72){\colorbox{gray}}
\put(14,63){\colorbox{gray}}
\put(14,67){\colorbox{gray}}
\put(14,72){\colorbox{gray}}

\put(20,63){\colorbox{gray}}
\put(20,67){\colorbox{gray}}
\put(20,72){\colorbox{gray}}
\put(24,63){\colorbox{gray}}
\put(24,67){\colorbox{gray}}
\put(24,72){\colorbox{gray}}
\put(29,63){\colorbox{gray}}
\put(29,67){\colorbox{gray}}
\put(29,72){\colorbox{gray}}

\put(35,63){\colorbox{gray}}
\put(35,67){\colorbox{gray}}
\put(35,72){\colorbox{gray}}
\put(39,63){\colorbox{gray}}
\put(39,67){\colorbox{gray}}
\put(39,72){\colorbox{gray}}
\put(44,63){\colorbox{gray}}
\put(44,67){\colorbox{gray}}
\put(44,72){\colorbox{gray}}

\put(50,63){\colorbox{gray}}
\put(50,67){\colorbox{gray}}
\put(50,72){\colorbox{gray}}
\put(54,63){\colorbox{gray}}
\put(54,67){\colorbox{gray}}
\put(54,72){\colorbox{gray}}
\put(59,63){\colorbox{gray}}
\put(59,67){\colorbox{gray}}
\put(59,72){\colorbox{gray}}

\put(65,63){\colorbox{gray}}
\put(65,67){\colorbox{gray}}
\put(65,72){\colorbox{gray}}
\put(69,63){\colorbox{gray}}
\put(69,67){\colorbox{gray}}
\put(69,72){\colorbox{gray}}
\put(74,63){\colorbox{gray}}
\put(74,67){\colorbox{gray}}
\put(74,72){\colorbox{gray}}

\put(80,63){\colorbox{gray}}
\put(80,67){\colorbox{gray}}
\put(80,72){\colorbox{gray}}
\put(84,63){\colorbox{gray}}
\put(84,67){\colorbox{gray}}
\put(84,72){\colorbox{gray}}
\put(89,63){\colorbox{gray}}
\put(89,67){\colorbox{gray}}
\put(89,72){\colorbox{gray}}

\put(5,48){\colorbox{gray}}
\put(5,52){\colorbox{gray}}
\put(5,57){\colorbox{gray}}
\put(9,48){\colorbox{gray}}
\put(9,52){\colorbox{gray}}
\put(9,57){\colorbox{gray}}
\put(14,48){\colorbox{gray}}
\put(14,52){\colorbox{gray}}
\put(14,57){\colorbox{gray}}

\put(20,48){\colorbox{gray}}
\put(20,52){\colorbox{gray}}
\put(20,57){\colorbox{gray}}
\put(24,48){\colorbox{gray}}
\put(24,52){\colorbox{gray}}
\put(24,57){\colorbox{gray}}
\put(29,48){\colorbox{gray}}
\put(29,52){\colorbox{gray}}
\put(29,57){\colorbox{gray}}

\put(35,48){\colorbox{gray}}
\put(35,52){\colorbox{gray}}
\put(35,57){\colorbox{gray}}
\put(39,48){\colorbox{gray}}
\put(39,52){\colorbox{gray}}
\put(39,57){\colorbox{gray}}
\put(44,48){\colorbox{gray}}
\put(44,52){\colorbox{gray}}
\put(44,57){\colorbox{gray}}

\put(50,48){\colorbox{gray}}
\put(50,52){\colorbox{gray}}
\put(50,57){\colorbox{gray}}
\put(54,48){\colorbox{gray}}
\put(54,52){\colorbox{gray}}
\put(54,57){\colorbox{gray}}
\put(59,48){\colorbox{gray}}
\put(59,52){\colorbox{gray}}
\put(59,57){\colorbox{gray}}

\put(65,48){\colorbox{gray}}
\put(65,52){\colorbox{gray}}
\put(65,57){\colorbox{gray}}
\put(69,48){\colorbox{gray}}
\put(69,52){\colorbox{gray}}
\put(69,57){\colorbox{gray}}
\put(74,48){\colorbox{gray}}
\put(74,52){\colorbox{gray}}
\put(74,57){\colorbox{gray}}

\put(80,48){\colorbox{gray}}
\put(80,52){\colorbox{gray}}
\put(80,57){\colorbox{gray}}
\put(84,48){\colorbox{gray}}
\put(84,52){\colorbox{gray}}
\put(84,57){\colorbox{gray}}
\put(89,48){\colorbox{gray}}
\put(89,52){\colorbox{gray}}
\put(89,57){\colorbox{gray}}

\put(50,33){\colorbox{gray}}
\put(50,37){\colorbox{gray}}
\put(50,42){\colorbox{gray}}
\put(54,33){\colorbox{gray}}
\put(54,37){\colorbox{gray}}
\put(54,42){\colorbox{gray}}
\put(59,33){\colorbox{gray}}
\put(59,37){\colorbox{gray}}
\put(59,42){\colorbox{gray}}

\put(65,33){\colorbox{gray}}
\put(65,37){\colorbox{gray}}
\put(65,42){\colorbox{gray}}
\put(69,33){\colorbox{gray}}
\put(69,37){\colorbox{gray}}
\put(69,42){\colorbox{gray}}
\put(74,33){\colorbox{gray}}
\put(74,37){\colorbox{gray}}
\put(74,42){\colorbox{gray}}

\put(80,33){\colorbox{gray}}
\put(80,37){\colorbox{gray}}
\put(80,42){\colorbox{gray}}
\put(84,33){\colorbox{gray}}
\put(84,37){\colorbox{gray}}
\put(84,42){\colorbox{gray}}
\put(89,33){\colorbox{gray}}
\put(89,37){\colorbox{gray}}
\put(89,42){\colorbox{gray}}

\put(50,18){\colorbox{gray}}
\put(50,22){\colorbox{gray}}
\put(50,27){\colorbox{gray}}
\put(54,18){\colorbox{gray}}
\put(54,22){\colorbox{gray}}
\put(54,27){\colorbox{gray}}
\put(59,18){\colorbox{gray}}
\put(59,22){\colorbox{gray}}
\put(59,27){\colorbox{gray}}

\put(65,18){\colorbox{gray}}
\put(65,22){\colorbox{gray}}
\put(65,27){\colorbox{gray}}
\put(69,18){\colorbox{gray}}
\put(69,22){\colorbox{gray}}
\put(69,27){\colorbox{gray}}
\put(74,18){\colorbox{gray}}
\put(74,22){\colorbox{gray}}
\put(74,27){\colorbox{gray}}

\put(80,18){\colorbox{gray}}
\put(80,22){\colorbox{gray}}
\put(80,27){\colorbox{gray}}
\put(84,18){\colorbox{gray}}
\put(84,22){\colorbox{gray}}
\put(84,27){\colorbox{gray}}
\put(89,18){\colorbox{gray}}
\put(89,22){\colorbox{gray}}
\put(89,27){\colorbox{gray}}

\put(65,3){\colorbox{gray}}
\put(65,7){\colorbox{gray}}
\put(65,12){\colorbox{gray}}
\put(69,3){\colorbox{gray}}
\put(69,7){\colorbox{gray}}
\put(69,12){\colorbox{gray}}
\put(74,3){\colorbox{gray}}
\put(74,7){\colorbox{gray}}
\put(74,12){\colorbox{gray}}

\put(80,3){\colorbox{gray}}
\put(80,7){\colorbox{gray}}
\put(80,12){\colorbox{gray}}
\put(84,3){\colorbox{gray}}
\put(84,7){\colorbox{gray}}
\put(84,12){\colorbox{gray}}
\put(89,3){\colorbox{gray}}
\put(89,7){\colorbox{gray}}
\put(89,12){\colorbox{gray}}

\put(-10,0){\framebox(15,15)}
\put(5,0){\framebox(15,15)}
\put(20,0){\framebox(15,15)}
\put(35,0){\framebox(15,15)}
\put(50,0){\framebox(15,15)}
\put(65,0){\framebox(15,15)}
\put(80,0){\framebox(15,15)}

\put(-10,15){\framebox(15,15)}
\put(5,15){\framebox(15,15)}
\put(20,15){\framebox(15,15)}
\put(35,15){\framebox(15,15)}
\put(50,15){\framebox(15,15)}
\put(65,15){\framebox(15,15)}
\put(80,15){\framebox(15,15)}

\put(-10,30){\framebox(15,15)}
\put(5,30){\framebox(15,15)}
\put(20,30){\framebox(15,15)}
\put(35,30){\framebox(15,15)}
\put(50,30){\framebox(15,15)}
\put(65,30){\framebox(15,15)}
\put(80,30){\framebox(15,15)}

\put(-10,45){\framebox(15,15)}
\put(5,45){\framebox(15,15)}
\put(20,45){\framebox(15,15)}
\put(35,45){\framebox(15,15)}
\put(50,45){\framebox(15,15)}
\put(65,45){\framebox(15,15)}
\put(80,45){\framebox(15,15)}

\put(-10,60){\framebox(15,15)}
\put(5,60){\framebox(15,15)}
\put(20,60){\framebox(15,15)}
\put(35,60){\framebox(15,15)}
\put(50,60){\framebox(15,15)}
\put(65,60){\framebox(15,15)}
\put(80,60){\framebox(15,15)}

\put(-10,75){\framebox(15,15)}
\put(5,75){\framebox(15,15)}
\put(20,75){\framebox(15,15)}
\put(35,75){\framebox(15,15)}
\put(50,75){\framebox(15,15)}
\put(65,75){\framebox(15,15)}
\put(80,75){\framebox(15,15)}

\put(-10,90){\framebox(15,15)}
\put(5,90){\framebox(15,15)}
\put(20,90){\framebox(15,15)}
\put(35,90){\framebox(15,15)}
\put(50,90){\framebox(15,15)}
\put(65,90){\framebox(15,15)}
\put(80,90){\framebox(15,15)}

\put(40,-15){$h$}
\put(-70,50){$h(j)=j \rightarrow$}
\put(40,110){$\downarrow$}
\put(41,124){$j$}

\put(-10,93){\colorbox{Gray}}
\put(-10,97){\colorbox{Gray}}
\put(-10,102){\colorbox{Gray}}
\put(-6,93){\colorbox{Gray}}
\put(-6,97){\colorbox{Gray}}
\put(-6,102){\colorbox{Gray}}
\put(-1,93){\colorbox{Gray}}
\put(-1,97){\colorbox{Gray}}
\put(-1,102){\colorbox{Gray}}

\put(5,93){\colorbox{Gray}}
\put(5,97){\colorbox{Gray}}
\put(5,102){\colorbox{Gray}}
\put(9,93){\colorbox{Gray}}
\put(9,97){\colorbox{Gray}}
\put(9,102){\colorbox{Gray}}
\put(14,93){\colorbox{Gray}}
\put(14,97){\colorbox{Gray}}
\put(14,102){\colorbox{Gray}}

\put(20,93){\colorbox{Gray}}
\put(20,97){\colorbox{Gray}}
\put(20,102){\colorbox{Gray}}
\put(24,93){\colorbox{Gray}}
\put(24,97){\colorbox{Gray}}
\put(24,102){\colorbox{Gray}}
\put(29,93){\colorbox{Gray}}
\put(29,97){\colorbox{Gray}}
\put(29,102){\colorbox{Gray}}

\put(35,93){\colorbox{Gray}}
\put(35,97){\colorbox{Gray}}
\put(35,102){\colorbox{Gray}}
\put(39,93){\colorbox{Gray}}
\put(39,97){\colorbox{Gray}}
\put(39,102){\colorbox{Gray}}
\put(44,93){\colorbox{Gray}}
\put(44,97){\colorbox{Gray}}
\put(44,102){\colorbox{Gray}}

\put(50,93){\colorbox{gray}}
\put(50,97){\colorbox{gray}}
\put(50,102){\colorbox{gray}}
\put(54,93){\colorbox{gray}}
\put(54,97){\colorbox{gray}}
\put(54,102){\colorbox{gray}}
\put(59,93){\colorbox{gray}}
\put(59,97){\colorbox{gray}}
\put(59,102){\colorbox{gray}}

\put(65,93){\colorbox{gray}}
\put(65,97){\colorbox{gray}}
\put(65,102){\colorbox{gray}}
\put(69,93){\colorbox{gray}}
\put(69,97){\colorbox{gray}}
\put(69,102){\colorbox{gray}}
\put(74,93){\colorbox{gray}}
\put(74,97){\colorbox{gray}}
\put(74,102){\colorbox{gray}}

\put(80,93){\colorbox{gray}}
\put(80,97){\colorbox{gray}}
\put(80,102){\colorbox{gray}}
\put(84,93){\colorbox{gray}}
\put(84,97){\colorbox{gray}}
\put(84,102){\colorbox{gray}}
\put(89,93){\colorbox{gray}}
\put(89,97){\colorbox{gray}}
\put(89,102){\colorbox{gray}}

\put(-10,78){\colorbox{Gray}}
\put(-10,82){\colorbox{Gray}}
\put(-10,87){\colorbox{Gray}}
\put(-6,78){\colorbox{Gray}}
\put(-6,82){\colorbox{Gray}}
\put(-6,87){\colorbox{Gray}}
\put(-1,78){\colorbox{Gray}}
\put(-1,82){\colorbox{Gray}}
\put(-1,87){\colorbox{Gray}}

\put(5,78){\colorbox{Gray}}
\put(5,82){\colorbox{Gray}}
\put(5,87){\colorbox{Gray}}
\put(9,78){\colorbox{Gray}}
\put(9,82){\colorbox{Gray}}
\put(9,87){\colorbox{Gray}}
\put(14,78){\colorbox{Gray}}
\put(14,82){\colorbox{Gray}}
\put(14,87){\colorbox{Gray}}

\put(20,78){\colorbox{Gray}}
\put(20,82){\colorbox{Gray}}
\put(20,87){\colorbox{Gray}}
\put(24,78){\colorbox{Gray}}
\put(24,82){\colorbox{Gray}}
\put(24,87){\colorbox{Gray}}
\put(29,78){\colorbox{Gray}}
\put(29,82){\colorbox{Gray}}
\put(29,87){\colorbox{Gray}}

\put(35,78){\colorbox{Gray}}
\put(35,82){\colorbox{Gray}}
\put(35,87){\colorbox{Gray}}
\put(39,78){\colorbox{Gray}}
\put(39,82){\colorbox{Gray}}
\put(39,87){\colorbox{Gray}}
\put(44,78){\colorbox{Gray}}
\put(44,82){\colorbox{Gray}}
\put(44,87){\colorbox{Gray}}

\put(50,78){\colorbox{gray}}
\put(50,82){\colorbox{gray}}
\put(50,87){\colorbox{gray}}
\put(54,78){\colorbox{gray}}
\put(54,82){\colorbox{gray}}
\put(54,87){\colorbox{gray}}
\put(59,78){\colorbox{gray}}
\put(59,82){\colorbox{gray}}
\put(59,87){\colorbox{gray}}

\put(65,78){\colorbox{gray}}
\put(65,82){\colorbox{gray}}
\put(65,87){\colorbox{gray}}
\put(69,78){\colorbox{gray}}
\put(69,82){\colorbox{gray}}
\put(69,87){\colorbox{gray}}
\put(74,78){\colorbox{gray}}
\put(74,82){\colorbox{gray}}
\put(74,87){\colorbox{gray}}

\put(80,78){\colorbox{gray}}
\put(80,82){\colorbox{gray}}
\put(80,87){\colorbox{gray}}
\put(84,78){\colorbox{gray}}
\put(84,82){\colorbox{gray}}
\put(84,87){\colorbox{gray}}
\put(89,78){\colorbox{gray}}
\put(89,82){\colorbox{gray}}
\put(89,87){\colorbox{gray}}

\put(5,63){\colorbox{Gray}}
\put(5,67){\colorbox{Gray}}
\put(5,72){\colorbox{Gray}}
\put(9,63){\colorbox{Gray}}
\put(9,67){\colorbox{Gray}}
\put(9,72){\colorbox{Gray}}
\put(14,63){\colorbox{Gray}}
\put(14,67){\colorbox{Gray}}
\put(14,72){\colorbox{Gray}}

\put(20,63){\colorbox{Gray}}
\put(20,67){\colorbox{Gray}}
\put(20,72){\colorbox{Gray}}
\put(24,63){\colorbox{Gray}}
\put(24,67){\colorbox{Gray}}
\put(24,72){\colorbox{Gray}}
\put(29,63){\colorbox{Gray}}
\put(29,67){\colorbox{Gray}}
\put(29,72){\colorbox{Gray}}

\put(35,63){\colorbox{Gray}}
\put(35,67){\colorbox{Gray}}
\put(35,72){\colorbox{Gray}}
\put(39,63){\colorbox{Gray}}
\put(39,67){\colorbox{Gray}}
\put(39,72){\colorbox{Gray}}
\put(44,63){\colorbox{Gray}}
\put(44,67){\colorbox{Gray}}
\put(44,72){\colorbox{Gray}}

\put(50,63){\colorbox{gray}}
\put(50,67){\colorbox{gray}}
\put(50,72){\colorbox{gray}}
\put(54,63){\colorbox{gray}}
\put(54,67){\colorbox{gray}}
\put(54,72){\colorbox{gray}}
\put(59,63){\colorbox{gray}}
\put(59,67){\colorbox{gray}}
\put(59,72){\colorbox{gray}}

\put(65,63){\colorbox{gray}}
\put(65,67){\colorbox{gray}}
\put(65,72){\colorbox{gray}}
\put(69,63){\colorbox{gray}}
\put(69,67){\colorbox{gray}}
\put(69,72){\colorbox{gray}}
\put(74,63){\colorbox{gray}}
\put(74,67){\colorbox{gray}}
\put(74,72){\colorbox{gray}}

\put(80,63){\colorbox{gray}}
\put(80,67){\colorbox{gray}}
\put(80,72){\colorbox{gray}}
\put(84,63){\colorbox{gray}}
\put(84,67){\colorbox{gray}}
\put(84,72){\colorbox{gray}}
\put(89,63){\colorbox{gray}}
\put(89,67){\colorbox{gray}}
\put(89,72){\colorbox{gray}}

\put(5,48){\colorbox{Gray}}
\put(5,52){\colorbox{Gray}}
\put(5,57){\colorbox{Gray}}
\put(9,48){\colorbox{Gray}}
\put(9,52){\colorbox{Gray}}
\put(9,57){\colorbox{Gray}}
\put(14,48){\colorbox{Gray}}
\put(14,52){\colorbox{Gray}}
\put(14,57){\colorbox{Gray}}

\put(20,48){\colorbox{Gray}}
\put(20,52){\colorbox{Gray}}
\put(20,57){\colorbox{Gray}}
\put(24,48){\colorbox{Gray}}
\put(24,52){\colorbox{Gray}}
\put(24,57){\colorbox{Gray}}
\put(29,48){\colorbox{Gray}}
\put(29,52){\colorbox{Gray}}
\put(29,57){\colorbox{Gray}}

\put(35,48){\colorbox{Gray}}
\put(35,52){\colorbox{Gray}}
\put(35,57){\colorbox{Gray}}
\put(39,48){\colorbox{Gray}}
\put(39,52){\colorbox{Gray}}
\put(39,57){\colorbox{Gray}}
\put(44,48){\colorbox{Gray}}
\put(44,52){\colorbox{Gray}}
\put(44,57){\colorbox{Gray}}


\put(50,48){\colorbox{gray}}
\put(50,52){\colorbox{gray}}
\put(50,57){\colorbox{gray}}
\put(54,48){\colorbox{gray}}
\put(54,52){\colorbox{gray}}
\put(54,57){\colorbox{gray}}
\put(59,48){\colorbox{gray}}
\put(59,52){\colorbox{gray}}
\put(59,57){\colorbox{gray}}

\put(65,48){\colorbox{gray}}
\put(65,52){\colorbox{gray}}
\put(65,57){\colorbox{gray}}
\put(69,48){\colorbox{gray}}
\put(69,52){\colorbox{gray}}
\put(69,57){\colorbox{gray}}
\put(74,48){\colorbox{gray}}
\put(74,52){\colorbox{gray}}
\put(74,57){\colorbox{gray}}

\put(80,48){\colorbox{gray}}
\put(80,52){\colorbox{gray}}
\put(80,57){\colorbox{gray}}
\put(84,48){\colorbox{gray}}
\put(84,52){\colorbox{gray}}
\put(84,57){\colorbox{gray}}
\put(89,48){\colorbox{gray}}
\put(89,52){\colorbox{gray}}
\put(89,57){\colorbox{gray}}

\put(50,33){\colorbox{Gray}}
\put(50,37){\colorbox{Gray}}
\put(50,42){\colorbox{Gray}}
\put(54,33){\colorbox{Gray}}
\put(54,37){\colorbox{Gray}}
\put(54,42){\colorbox{Gray}}
\put(59,33){\colorbox{Gray}}
\put(59,37){\colorbox{Gray}}
\put(59,42){\colorbox{Gray}}

\put(65,33){\colorbox{Gray}}
\put(65,37){\colorbox{Gray}}
\put(65,42){\colorbox{Gray}}
\put(69,33){\colorbox{Gray}}
\put(69,37){\colorbox{Gray}}
\put(69,42){\colorbox{Gray}}
\put(74,33){\colorbox{Gray}}
\put(74,37){\colorbox{Gray}}
\put(74,42){\colorbox{Gray}}

\put(80,33){\colorbox{Gray}}
\put(80,37){\colorbox{Gray}}
\put(80,42){\colorbox{Gray}}
\put(84,33){\colorbox{Gray}}
\put(84,37){\colorbox{Gray}}
\put(84,42){\colorbox{Gray}}
\put(89,33){\colorbox{Gray}}
\put(89,37){\colorbox{Gray}}
\put(89,42){\colorbox{Gray}}

\put(50,18){\colorbox{Gray}}
\put(50,22){\colorbox{Gray}}
\put(50,27){\colorbox{Gray}}
\put(54,18){\colorbox{Gray}}
\put(54,22){\colorbox{Gray}}
\put(54,27){\colorbox{Gray}}
\put(59,18){\colorbox{Gray}}
\put(59,22){\colorbox{Gray}}
\put(59,27){\colorbox{Gray}}

\put(65,18){\colorbox{Gray}}
\put(65,22){\colorbox{Gray}}
\put(65,27){\colorbox{Gray}}
\put(69,18){\colorbox{Gray}}
\put(69,22){\colorbox{Gray}}
\put(69,27){\colorbox{Gray}}
\put(74,18){\colorbox{Gray}}
\put(74,22){\colorbox{Gray}}
\put(74,27){\colorbox{Gray}}

\put(80,18){\colorbox{Gray}}
\put(80,22){\colorbox{Gray}}
\put(80,27){\colorbox{Gray}}
\put(84,18){\colorbox{Gray}}
\put(84,22){\colorbox{Gray}}
\put(84,27){\colorbox{Gray}}
\put(89,18){\colorbox{Gray}}
\put(89,22){\colorbox{Gray}}
\put(89,27){\colorbox{Gray}}

\put(65,3){\colorbox{Gray}}
\put(65,7){\colorbox{Gray}}
\put(65,12){\colorbox{Gray}}
\put(69,3){\colorbox{Gray}}
\put(69,7){\colorbox{Gray}}
\put(69,12){\colorbox{Gray}}
\put(74,3){\colorbox{Gray}}
\put(74,7){\colorbox{Gray}}
\put(74,12){\colorbox{Gray}}

\put(80,3){\colorbox{Gray}}
\put(80,7){\colorbox{Gray}}
\put(80,12){\colorbox{Gray}}
\put(84,3){\colorbox{Gray}}
\put(84,7){\colorbox{Gray}}
\put(84,12){\colorbox{Gray}}
\put(89,3){\colorbox{Gray}}
\put(89,7){\colorbox{Gray}}
\put(89,12){\colorbox{Gray}}

\put(-10,0){\framebox(15,15)}
\put(5,0){\framebox(15,15)}
\put(20,0){\framebox(15,15)}
\put(35,0){\framebox(15,15)}
\put(50,0){\framebox(15,15)}
\put(65,0){\framebox(15,15)}
\put(80,0){\framebox(15,15)}

\put(-10,15){\framebox(15,15)}
\put(5,15){\framebox(15,15)}
\put(20,15){\framebox(15,15)}
\put(35,15){\framebox(15,15)}
\put(50,15){\framebox(15,15)}
\put(65,15){\framebox(15,15)}
\put(80,15){\framebox(15,15)}

\put(-10,30){\framebox(15,15)}
\put(5,30){\framebox(15,15)}
\put(20,30){\framebox(15,15)}
\put(35,30){\framebox(15,15)}
\put(50,30){\framebox(15,15)}
\put(65,30){\framebox(15,15)}
\put(80,30){\framebox(15,15)}

\put(-10,45){\framebox(15,15)}
\put(5,45){\framebox(15,15)}
\put(20,45){\framebox(15,15)}
\put(35,45){\framebox(15,15)}
\put(50,45){\framebox(15,15)}
\put(65,45){\framebox(15,15)}
\put(80,45){\framebox(15,15)}

\put(-10,60){\framebox(15,15)}
\put(5,60){\framebox(15,15)}
\put(20,60){\framebox(15,15)}
\put(35,60){\framebox(15,15)}
\put(50,60){\framebox(15,15)}
\put(65,60){\framebox(15,15)}
\put(80,60){\framebox(15,15)}

\put(-10,75){\framebox(15,15)}
\put(5,75){\framebox(15,15)}
\put(20,75){\framebox(15,15)}
\put(35,75){\framebox(15,15)}
\put(50,75){\framebox(15,15)}
\put(65,75){\framebox(15,15)}
\put(80,75){\framebox(15,15)}

\put(-10,90){\framebox(15,15)}
\put(5,90){\framebox(15,15)}
\put(20,90){\framebox(15,15)}
\put(35,90){\framebox(15,15)}
\put(50,90){\framebox(15,15)}
\put(65,90){\framebox(15,15)}
\put(80,90){\framebox(15,15)}

\put(160,93){\colorbox{Gray}}
\put(160,97){\colorbox{Gray}}
\put(160,102){\colorbox{Gray}}
\put(164,93){\colorbox{Gray}}
\put(164,97){\colorbox{Gray}}
\put(164,102){\colorbox{Gray}}
\put(169,93){\colorbox{Gray}}
\put(169,97){\colorbox{Gray}}
\put(169,102){\colorbox{Gray}}

\put(175,93){\colorbox{Gray}}
\put(175,97){\colorbox{Gray}}
\put(175,102){\colorbox{Gray}}
\put(179,93){\colorbox{Gray}}
\put(179,97){\colorbox{Gray}}
\put(179,102){\colorbox{Gray}}
\put(184,93){\colorbox{Gray}}
\put(184,97){\colorbox{Gray}}
\put(184,102){\colorbox{Gray}}

\put(190,93){\colorbox{Gray}}
\put(190,97){\colorbox{Gray}}
\put(190,102){\colorbox{Gray}}
\put(194,93){\colorbox{Gray}}
\put(194,97){\colorbox{Gray}}
\put(194,102){\colorbox{Gray}}
\put(199,93){\colorbox{Gray}}
\put(199,97){\colorbox{Gray}}
\put(199,102){\colorbox{Gray}}

\put(205,93){\colorbox{Gray}}
\put(205,97){\colorbox{Gray}}
\put(205,102){\colorbox{Gray}}
\put(209,93){\colorbox{Gray}}
\put(209,97){\colorbox{Gray}}
\put(209,102){\colorbox{Gray}}
\put(214,93){\colorbox{Gray}}
\put(214,97){\colorbox{Gray}}
\put(214,102){\colorbox{Gray}}

\put(160,78){\colorbox{Gray}}
\put(160,82){\colorbox{Gray}}
\put(160,87){\colorbox{Gray}}
\put(164,78){\colorbox{Gray}}
\put(164,82){\colorbox{Gray}}
\put(164,87){\colorbox{Gray}}
\put(169,78){\colorbox{Gray}}
\put(169,82){\colorbox{Gray}}
\put(169,87){\colorbox{Gray}}

\put(175,78){\colorbox{Gray}}
\put(175,82){\colorbox{Gray}}
\put(175,87){\colorbox{Gray}}
\put(179,78){\colorbox{Gray}}
\put(179,82){\colorbox{Gray}}
\put(179,87){\colorbox{Gray}}
\put(184,78){\colorbox{Gray}}
\put(184,82){\colorbox{Gray}}
\put(184,87){\colorbox{Gray}}

\put(190,78){\colorbox{Gray}}
\put(190,82){\colorbox{Gray}}
\put(190,87){\colorbox{Gray}}
\put(194,78){\colorbox{Gray}}
\put(194,82){\colorbox{Gray}}
\put(194,87){\colorbox{Gray}}
\put(199,78){\colorbox{Gray}}
\put(199,82){\colorbox{Gray}}
\put(199,87){\colorbox{Gray}}

\put(205,78){\colorbox{Gray}}
\put(205,82){\colorbox{Gray}}
\put(205,87){\colorbox{Gray}}
\put(209,78){\colorbox{Gray}}
\put(209,82){\colorbox{Gray}}
\put(209,87){\colorbox{Gray}}
\put(214,78){\colorbox{Gray}}
\put(214,82){\colorbox{Gray}}
\put(214,87){\colorbox{Gray}}

\put(175,63){\colorbox{Gray}}
\put(175,67){\colorbox{Gray}}
\put(175,72){\colorbox{Gray}}
\put(179,63){\colorbox{Gray}}
\put(179,67){\colorbox{Gray}}
\put(179,72){\colorbox{Gray}}
\put(184,63){\colorbox{Gray}}
\put(184,67){\colorbox{Gray}}
\put(184,72){\colorbox{Gray}}

\put(190,63){\colorbox{Gray}}
\put(190,67){\colorbox{Gray}}
\put(190,72){\colorbox{Gray}}
\put(194,63){\colorbox{Gray}}
\put(194,67){\colorbox{Gray}}
\put(194,72){\colorbox{Gray}}
\put(199,63){\colorbox{Gray}}
\put(199,67){\colorbox{Gray}}
\put(199,72){\colorbox{Gray}}

\put(205,63){\colorbox{Gray}}
\put(205,67){\colorbox{Gray}}
\put(205,72){\colorbox{Gray}}
\put(209,63){\colorbox{Gray}}
\put(209,67){\colorbox{Gray}}
\put(209,72){\colorbox{Gray}}
\put(214,63){\colorbox{Gray}}
\put(214,67){\colorbox{Gray}}
\put(214,72){\colorbox{Gray}}

\put(175,48){\colorbox{Gray}}
\put(175,52){\colorbox{Gray}}
\put(175,57){\colorbox{Gray}}
\put(179,48){\colorbox{Gray}}
\put(179,52){\colorbox{Gray}}
\put(179,57){\colorbox{Gray}}
\put(184,48){\colorbox{Gray}}
\put(184,52){\colorbox{Gray}}
\put(184,57){\colorbox{Gray}}

\put(190,48){\colorbox{Gray}}
\put(190,52){\colorbox{Gray}}
\put(190,57){\colorbox{Gray}}
\put(194,48){\colorbox{Gray}}
\put(194,52){\colorbox{Gray}}
\put(194,57){\colorbox{Gray}}
\put(199,48){\colorbox{Gray}}
\put(199,52){\colorbox{Gray}}
\put(199,57){\colorbox{Gray}}

\put(205,48){\colorbox{Gray}}
\put(205,52){\colorbox{Gray}}
\put(205,57){\colorbox{Gray}}
\put(209,48){\colorbox{Gray}}
\put(209,52){\colorbox{Gray}}
\put(209,57){\colorbox{Gray}}
\put(214,48){\colorbox{Gray}}
\put(214,52){\colorbox{Gray}}
\put(214,57){\colorbox{Gray}}


\put(220,33){\colorbox{Gray}}
\put(220,37){\colorbox{Gray}}
\put(220,42){\colorbox{Gray}}
\put(224,33){\colorbox{Gray}}
\put(224,37){\colorbox{Gray}}
\put(224,42){\colorbox{Gray}}
\put(229,33){\colorbox{Gray}}
\put(229,37){\colorbox{Gray}}
\put(229,42){\colorbox{Gray}}

\put(235,33){\colorbox{Gray}}
\put(235,37){\colorbox{Gray}}
\put(235,42){\colorbox{Gray}}
\put(239,33){\colorbox{Gray}}
\put(239,37){\colorbox{Gray}}
\put(239,42){\colorbox{Gray}}
\put(244,33){\colorbox{Gray}}
\put(244,37){\colorbox{Gray}}
\put(244,42){\colorbox{Gray}}

\put(250,33){\colorbox{Gray}}
\put(250,37){\colorbox{Gray}}
\put(250,42){\colorbox{Gray}}
\put(254,33){\colorbox{Gray}}
\put(254,37){\colorbox{Gray}}
\put(254,42){\colorbox{Gray}}
\put(259,33){\colorbox{Gray}}
\put(259,37){\colorbox{Gray}}
\put(259,42){\colorbox{Gray}}

\put(220,18){\colorbox{Gray}}
\put(220,22){\colorbox{Gray}}
\put(220,27){\colorbox{Gray}}
\put(224,18){\colorbox{Gray}}
\put(224,22){\colorbox{Gray}}
\put(224,27){\colorbox{Gray}}
\put(229,18){\colorbox{Gray}}
\put(229,22){\colorbox{Gray}}
\put(229,27){\colorbox{Gray}}

\put(235,18){\colorbox{Gray}}
\put(235,22){\colorbox{Gray}}
\put(235,27){\colorbox{Gray}}
\put(239,18){\colorbox{Gray}}
\put(239,22){\colorbox{Gray}}
\put(239,27){\colorbox{Gray}}
\put(244,18){\colorbox{Gray}}
\put(244,22){\colorbox{Gray}}
\put(244,27){\colorbox{Gray}}

\put(250,18){\colorbox{Gray}}
\put(250,22){\colorbox{Gray}}
\put(250,27){\colorbox{Gray}}
\put(254,18){\colorbox{Gray}}
\put(254,22){\colorbox{Gray}}
\put(254,27){\colorbox{Gray}}
\put(259,18){\colorbox{Gray}}
\put(259,22){\colorbox{Gray}}
\put(259,27){\colorbox{Gray}}

\put(235,3){\colorbox{Gray}}
\put(235,7){\colorbox{Gray}}
\put(235,12){\colorbox{Gray}}
\put(239,3){\colorbox{Gray}}
\put(239,7){\colorbox{Gray}}
\put(239,12){\colorbox{Gray}}
\put(244,3){\colorbox{Gray}}
\put(244,7){\colorbox{Gray}}
\put(244,12){\colorbox{Gray}}

\put(250,3){\colorbox{Gray}}
\put(250,7){\colorbox{Gray}}
\put(250,12){\colorbox{Gray}}
\put(254,3){\colorbox{Gray}}
\put(254,7){\colorbox{Gray}}
\put(254,12){\colorbox{Gray}}
\put(259,3){\colorbox{Gray}}
\put(259,7){\colorbox{Gray}}
\put(259,12){\colorbox{Gray}}

\put(220,0){\framebox(15,15)}
\put(235,0){\framebox(15,15)}
\put(250,0){\framebox(15,15)}

\put(220,15){\framebox(15,15)}
\put(235,15){\framebox(15,15)}
\put(250,15){\framebox(15,15)}

\put(220,30){\framebox(15,15)}
\put(235,30){\framebox(15,15)}
\put(250,30){\framebox(15,15)}

\put(160,45){\framebox(15,15)}
\put(175,45){\framebox(15,15)}
\put(190,45){\framebox(15,15)}
\put(205,45){\framebox(15,15)}

\put(160,60){\framebox(15,15)}
\put(175,60){\framebox(15,15)}
\put(190,60){\framebox(15,15)}
\put(205,60){\framebox(15,15)}

\put(160,75){\framebox(15,15)}
\put(175,75){\framebox(15,15)}
\put(190,75){\framebox(15,15)}
\put(205,75){\framebox(15,15)}

\put(160,90){\framebox(15,15)}
\put(175,90){\framebox(15,15)}
\put(190,90){\framebox(15,15)}
\put(205,90){\framebox(15,15)}

\put(125,72){$h^{(1)}=$}
\put(185,18){$h^{(2)}=$}

\end{picture}
\end{center}
\vspace{5pt}
\caption{Decomposition of $h$ into $h^{(1)}$ and $h^{(2)}$.}
\label{pic:Decomposition of $h$ by two smaller Hessenberg functions $h^{(1)}$ and $h^{(2)}$}
\end{figure}
\vspace{-5pt}

\noindent
Then, $\Hess(N,h)$ is decomposed as the product of regular nilpotent Hessenberg varieties associated with $h^{(1)}:[j]\rightarrow[j]$ and $h^{(2)}:[n-j]\rightarrow[n-j]$ of smaller sizes (\cite[Theorem~4.5]{Dre1}). In this sense, the condition $h(j) \geq j+1$ for any $j<n$ is essential, and for such Hessenberg functions, $h=(2,3,4,\ldots,n,n)$ is the minimum one.
Hence, from $\eqref{eq:inclusion-preserving}$, we may regard regular nilpotent Hessenberg varieties $\Hess(N,h)$ as a (discrete) family of subvarieties of $Fl(\C^n)$ connecting the Peterson variety and the flag variety itself.

In the following theorem, we summarize some basic properties of $\Hess(N,h)$. For this purpose, let us prepare some notations.
Given a Hessenberg function $h$, we denote by $S_n^{h}$ the subset of the $n$-th symmetric group $S_n$ defined as 
\begin{equation} \label{eq:RegularNilpotentAndSchubertCell}
S_n^{h}:=\{w \in S_n \mid w^{-1}(w(j)-1) \leq h(j) \ \ {\rm for \ all \ } j \in [n] \},
\end{equation}
where we take by convention $w^{-1}(w(j)-1)=0$ whenever $w(j)-1=0$. The condition for $w \in S_n^h$ in \eqref{eq:RegularNilpotentAndSchubertCell} is exactly the condition that the permutation flag\footnote{For $w\in S_n$, the associated permutation flag $V_{\bullet}$ is given by $V_i:={ \rm span }_{\C}\{e_{w(1)},\ldots,e_{w(i)}\}$ where $\{e_{1},\ldots,e_{n}\}$ is the standard basis of $\C^n$.} 
associated with $w$ belongs to $\Hess(N,h)$, and this condition is equivalent to the condition that the Schubert cell $X_w^{\circ}$ intersects with $\Hess(N,h)$ (\cite{Ty2}). That is, we have
\begin{equation*} 
w \in S_n^h \iff X_w^{\circ} \cap \Hess(N,h) \neq \emptyset.
\end{equation*} 
Here, the dimension of the Schubert cell $X_w^{\circ}$ is equal to $\ell(w)$, which is the length of $w$. 
It is known from \cite{Ty2} that the dimension of the intersection $X_w^{\circ} \cap \Hess(N,h)$ is equal to 
\begin{equation} \label{eq:HessenbergLength}
\ell_h(w):=\#\{(j,i) \mid 1\leq j < i \leq n, \ w(j)>w(i), \ i\leq h(j) \},
\end{equation}
and that the intersections $X_w^{\circ} \cap \Hess(N,h)$ form a paving by affines of $\Hess(N,h)$. Combining the works of D. Anderson, E. Insko, B. Kostant, E. Sommers, J. Tymoczko, and A. Yong, we now give some basic properties of $\Hess(N,h)$.

\begin{theorem}[\cite{an-ty}, \cite{Ko,InskoYong}, \cite{Ty2}, \cite{so-ty}] \label{theorem:PropertyRegularNilpotent}
Let $\Hess(N,h)$ be a regular nilpotent Hessenberg variety. 
Then the following hold. 
\begin{enumerate}
\item $\Hess(N,h)$ is irreducible, and it is singular in general.  
\item The $($complex$)$ dimension of $\Hess(N,h)$ is equal to $\sum _{j=1}^n (h(j)-j)$.
\item The Poincar\'{e} polynomial \eqref{eq:PoincarePolynomial} for $X=N$ has the following two types of expressions: 
\begin{align}
\Poin(\Hess(N,h), \q)=&\sum_{w\in S_n^{h}} \q^{\ell_h(w)} \label{eq:PoinHess(N,h)Sum} \\
=&\displaystyle\prod _{j=1}^n (1+\q+\q^2+\cdots+\q^{h(j)-j}) \label{eq:PoinHess(N,h)Prod} 
\end{align}
\end{enumerate}
\end{theorem}

For generalizations of Theorem \ref{theorem:PropertyRegularNilpotent} to arbitrary Lie type, see \cite{precup18, Ty2,precup13a,so-ty,Roe,AHMMS}.

Note that the (complex) dimension of $\Hess(N,h)$ given in Theorem \ref{theorem:PropertyRegularNilpotent} (2) is equal to the number of boxes in which lie strictly below the diagonal under the identification of a Hessenberg function $h$ and a configuration of boxes. 

\begin{example}
Let $h=(3,3,4,5,5)$. Then we have that $\dim_\C \Hess(N,h)=5$, which is the number of boxes which lie strictly below the diagonal. 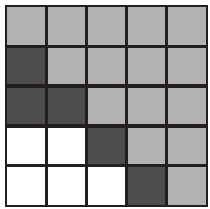
\begin{figure}[h]
\begin{center}
\begin{picture}(75,75)
\put(0,63){\colorbox{gray}}
\put(0,67){\colorbox{gray}}
\put(0,72){\colorbox{gray}}
\put(4,63){\colorbox{gray}}
\put(4,67){\colorbox{gray}}
\put(4,72){\colorbox{gray}}
\put(9,63){\colorbox{gray}}
\put(9,67){\colorbox{gray}}
\put(9,72){\colorbox{gray}}

\put(15,63){\colorbox{gray}}
\put(15,67){\colorbox{gray}}
\put(15,72){\colorbox{gray}}
\put(19,63){\colorbox{gray}}
\put(19,67){\colorbox{gray}}
\put(19,72){\colorbox{gray}}
\put(24,63){\colorbox{gray}}
\put(24,67){\colorbox{gray}}
\put(24,72){\colorbox{gray}}

\put(30,63){\colorbox{gray}}
\put(30,67){\colorbox{gray}}
\put(30,72){\colorbox{gray}}
\put(34,63){\colorbox{gray}}
\put(34,67){\colorbox{gray}}
\put(34,72){\colorbox{gray}}
\put(39,63){\colorbox{gray}}
\put(39,67){\colorbox{gray}}
\put(39,72){\colorbox{gray}}

\put(45,63){\colorbox{gray}}
\put(45,67){\colorbox{gray}}
\put(45,72){\colorbox{gray}}
\put(49,63){\colorbox{gray}}
\put(49,67){\colorbox{gray}}
\put(49,72){\colorbox{gray}}
\put(54,63){\colorbox{gray}}
\put(54,67){\colorbox{gray}}
\put(54,72){\colorbox{gray}}

\put(60,63){\colorbox{gray}}
\put(60,67){\colorbox{gray}}
\put(60,72){\colorbox{gray}}
\put(64,63){\colorbox{gray}}
\put(64,67){\colorbox{gray}}
\put(64,72){\colorbox{gray}}
\put(69,63){\colorbox{gray}}
\put(69,67){\colorbox{gray}}
\put(69,72){\colorbox{gray}}

\put(0,48){\colorbox{Gray}}
\put(0,52){\colorbox{Gray}}
\put(0,57){\colorbox{Gray}}
\put(4,48){\colorbox{Gray}}
\put(4,52){\colorbox{Gray}}
\put(4,57){\colorbox{Gray}}
\put(9,48){\colorbox{Gray}}
\put(9,52){\colorbox{Gray}}
\put(9,57){\colorbox{Gray}}

\put(15,48){\colorbox{gray}}
\put(15,52){\colorbox{gray}}
\put(15,57){\colorbox{gray}}
\put(19,48){\colorbox{gray}}
\put(19,52){\colorbox{gray}}
\put(19,57){\colorbox{gray}}
\put(24,48){\colorbox{gray}}
\put(24,52){\colorbox{gray}}
\put(24,57){\colorbox{gray}}

\put(30,48){\colorbox{gray}}
\put(30,52){\colorbox{gray}}
\put(30,57){\colorbox{gray}}
\put(34,48){\colorbox{gray}}
\put(34,52){\colorbox{gray}}
\put(34,57){\colorbox{gray}}
\put(39,48){\colorbox{gray}}
\put(39,52){\colorbox{gray}}
\put(39,57){\colorbox{gray}}

\put(45,48){\colorbox{gray}}
\put(45,52){\colorbox{gray}}
\put(45,57){\colorbox{gray}}
\put(49,48){\colorbox{gray}}
\put(49,52){\colorbox{gray}}
\put(49,57){\colorbox{gray}}
\put(54,48){\colorbox{gray}}
\put(54,52){\colorbox{gray}}
\put(54,57){\colorbox{gray}}

\put(60,48){\colorbox{gray}}
\put(60,52){\colorbox{gray}}
\put(60,57){\colorbox{gray}}
\put(64,48){\colorbox{gray}}
\put(64,52){\colorbox{gray}}
\put(64,57){\colorbox{gray}}
\put(69,48){\colorbox{gray}}
\put(69,52){\colorbox{gray}}
\put(69,57){\colorbox{gray}}

\put(0,33){\colorbox{Gray}}
\put(0,37){\colorbox{Gray}}
\put(0,42){\colorbox{Gray}}
\put(4,33){\colorbox{Gray}}
\put(4,37){\colorbox{Gray}}
\put(4,42){\colorbox{Gray}}
\put(9,33){\colorbox{Gray}}
\put(9,37){\colorbox{Gray}}
\put(9,42){\colorbox{Gray}}

\put(15,33){\colorbox{Gray}}
\put(15,37){\colorbox{Gray}}
\put(15,42){\colorbox{Gray}}
\put(19,33){\colorbox{Gray}}
\put(19,37){\colorbox{Gray}}
\put(19,42){\colorbox{Gray}}
\put(24,33){\colorbox{Gray}}
\put(24,37){\colorbox{Gray}}
\put(24,42){\colorbox{Gray}}

\put(30,33){\colorbox{gray}}
\put(30,37){\colorbox{gray}}
\put(30,42){\colorbox{gray}}
\put(34,33){\colorbox{gray}}
\put(34,37){\colorbox{gray}}
\put(34,42){\colorbox{gray}}
\put(39,33){\colorbox{gray}}
\put(39,37){\colorbox{gray}}
\put(39,42){\colorbox{gray}}

\put(45,33){\colorbox{gray}}
\put(45,37){\colorbox{gray}}
\put(45,42){\colorbox{gray}}
\put(49,33){\colorbox{gray}}
\put(49,37){\colorbox{gray}}
\put(49,42){\colorbox{gray}}
\put(54,33){\colorbox{gray}}
\put(54,37){\colorbox{gray}}
\put(54,42){\colorbox{gray}}

\put(60,33){\colorbox{gray}}
\put(60,37){\colorbox{gray}}
\put(60,42){\colorbox{gray}}
\put(64,33){\colorbox{gray}}
\put(64,37){\colorbox{gray}}
\put(64,42){\colorbox{gray}}
\put(69,33){\colorbox{gray}}
\put(69,37){\colorbox{gray}}
\put(69,42){\colorbox{gray}}

\put(30,18){\colorbox{Gray}}
\put(30,22){\colorbox{Gray}}
\put(30,27){\colorbox{Gray}}
\put(34,18){\colorbox{Gray}}
\put(34,22){\colorbox{Gray}}
\put(34,27){\colorbox{Gray}}
\put(39,18){\colorbox{Gray}}
\put(39,22){\colorbox{Gray}}
\put(39,27){\colorbox{Gray}}

\put(45,18){\colorbox{gray}}
\put(45,22){\colorbox{gray}}
\put(45,27){\colorbox{gray}}
\put(49,18){\colorbox{gray}}
\put(49,22){\colorbox{gray}}
\put(49,27){\colorbox{gray}}
\put(54,18){\colorbox{gray}}
\put(54,22){\colorbox{gray}}
\put(54,27){\colorbox{gray}}

\put(60,18){\colorbox{gray}}
\put(60,22){\colorbox{gray}}
\put(60,27){\colorbox{gray}}
\put(64,18){\colorbox{gray}}
\put(64,22){\colorbox{gray}}
\put(64,27){\colorbox{gray}}
\put(69,18){\colorbox{gray}}
\put(69,22){\colorbox{gray}}
\put(69,27){\colorbox{gray}}

\put(45,3){\colorbox{Gray}}
\put(45,7){\colorbox{Gray}}
\put(45,12){\colorbox{Gray}}
\put(49,3){\colorbox{Gray}}
\put(49,7){\colorbox{Gray}}
\put(49,12){\colorbox{Gray}}
\put(54,3){\colorbox{Gray}}
\put(54,7){\colorbox{Gray}}
\put(54,12){\colorbox{Gray}}

\put(60,3){\colorbox{gray}}
\put(60,7){\colorbox{gray}}
\put(60,12){\colorbox{gray}}
\put(64,3){\colorbox{gray}}
\put(64,7){\colorbox{gray}}
\put(64,12){\colorbox{gray}}
\put(69,3){\colorbox{gray}}
\put(69,7){\colorbox{gray}}
\put(69,12){\colorbox{gray}}

\put(0,0){\framebox(15,15)}
\put(15,0){\framebox(15,15)}
\put(30,0){\framebox(15,15)}
\put(45,0){\framebox(15,15)}
\put(60,0){\framebox(15,15)}
\put(0,15){\framebox(15,15)}
\put(15,15){\framebox(15,15)}
\put(30,15){\framebox(15,15)}
\put(45,15){\framebox(15,15)}
\put(60,15){\framebox(15,15)}
\put(0,30){\framebox(15,15)}
\put(15,30){\framebox(15,15)}
\put(30,30){\framebox(15,15)}
\put(45,30){\framebox(15,15)}
\put(60,30){\framebox(15,15)}
\put(0,45){\framebox(15,15)}
\put(15,45){\framebox(15,15)}
\put(30,45){\framebox(15,15)}
\put(45,45){\framebox(15,15)}
\put(60,45){\framebox(15,15)}
\put(0,60){\framebox(15,15)}
\put(15,60){\framebox(15,15)}
\put(30,60){\framebox(15,15)}
\put(45,60){\framebox(15,15)}
\put(60,60){\framebox(15,15)}
\end{picture}
\end{center}
\vspace{-10pt}
\caption{Boxes which lie strictly below the diagonal.}
\label{pic:dimHess(N,h)}
\end{figure}
\end{example}

\begin{example} \label{example:Pet3}
Let $h=(2,3,3)$. Then $S_3^h=\{123, 213, 132, 321 \}$ 
where we use the standard one-line notation $w = w(1) \ w(2) \ \cdots \ w(n)$ for permutations in $S_n$ throughout this document.
In particular, we have geometrically that 
$$
\Hess(N,h) \cap X_w^{\circ}=\emptyset \iff w=312, 231.
$$
Since $\ell_h(123)=0$, $\ell_h(213)=\ell_h(132)=1$, $\ell_h(321)=2$ for each permutation in $S_3^h$, we have
\begin{align}\label{eq: Poinc for 233}
\Poin(\Hess(N,h),\q)=1+2\q+\q^2=(1+\q)^2.
\end{align}
\end{example}

\smallskip

\subsection{Regular semisimple Hessenberg varieties}
\label{subsec:Regular semisimple Hessenberg varieties}
Let $S$ be a regular semisimple matrix of size $n\times n$. In Jordan canonical form, it is given by
\begin{equation*}
S = 
\begin{pmatrix} 
c_1 &  & & \\
&  c_2 &  & \\ 
&  &  \ddots & \\
& &  & c_n 
\end{pmatrix}
\end{equation*}
where $c_1, c_2 \ldots, c_n$ are mutually distinct complex numbers.
For a Hessenberg function $h:[n]\rightarrow [n]$,  $\Hess(S,h)$ is called a \textbf{regular semisimple Hessenberg variety}. It is known that the topology of $\Hess(S,h)$ is independent of the choices of the (distinct) eigenvalues\footnote{For regular semisimple matrices $S$ and $S'$, the associated Hessenberg varieties $\Hess(S,h)$ and $\Hess(S',h)$ with a same Hessenberg function $h$ are diffeomorphic.}, and hence one may think that the topology of $\Hess(S,h)$ only depends on $h$. Based on this fact, 
we have the following two main examples for this class of Hessenberg varieties. 
For $h=(n,n,\ldots,n)$, $\Hess(S,h)$ is the flag variety $Fl(\C^n)$ itself.
For $h=(2,3,4,\ldots,n,n)$, $\Hess(S,h)$ is called the \textbf{permutohedral variety} which is the toric variety associated with the fan consisting of the collection of Weyl chambers of the root system of type $A_{n-1}$ (\cite[Theorem~11]{ma-pr-sh}). Similarly to the case for $\Hess(N,h)$ in Section \ref{subsec:Regular nilpotent Hessenberg varieties}, the Hessenberg function $h=(2,3,4,\ldots,n,n)$ gives the minimum for this class of Hessenberg varieties as well in the following sense. 

If we have $h(j)=j$ for some $j<n$, it is known that $\Hess(S,h)$ is not connected but equidimensional.
In fact, all the connected components are isomorphic, and each component is decomposed as the product of regular semisimple Hessenberg varieties of smaller sizes as in Section \ref{subsec:Regular nilpotent Hessenberg varieties}. See \cite{Teff11} for detail description of the connected components.
Hence the condition $h(j) \geq j+1$ for any $j<n$ is essential, and we may regard regular semisimple Hessenberg varieties $\Hess(S,h)$ as a (discrete) family of subvarieties of the flag variety connecting the permutohedral variety and the flag variety itself.
De Mari-Procesi-Shayman proved the following properties of $\Hess(S,h)$ (for arbitrary Lie type).

\begin{theorem}[\cite{ma-pr-sh}] \label{theorem:PropertyRegularSemisimple}
Let $\Hess(S,h)$ be a regular semisimple Hessenberg variety. 
Then the following hold. 
\begin{enumerate}
\item $\Hess(S,h)$ is smooth, and it is connected if and only if $h(j) \geq j+1$ for all $j<n$.
\item The $($complex$)$ dimension of $\Hess(S,h)$ is equal to $\sum _{j=1}^n (h(j)-j)$.
\item The Poincar\'{e} polynomial \eqref{eq:PoincarePolynomial} for $X=S$ has the following expression: 
$$
\Poin(\Hess(S,h),\q)=\displaystyle\sum _{w\in S_n} q^{\ell_h(w)}
$$
where $\ell_h(w)$ is defined in \eqref{eq:HessenbergLength}. 
\end{enumerate}
\end{theorem}

From Theorem \ref{theorem:PropertyRegularNilpotent} (2) and Theorem \ref{theorem:PropertyRegularSemisimple} (2), we see that
$$
\dim_\C \Hess(N,h)=\dim_\C \Hess(S,h)=\sum _{j=1}^n (h(j)-j).
$$

Unlike the situation for $\Hess(N,h)$, any regular semisimple Hessenberg variety $\Hess(S,h)$ intersects with all the Schubert cells $X_w^{\circ}$. 
It is known that the dimension of the intersection $X_w^{\circ} \cap \Hess(S,h)$ is equal to $\ell_h(w)$ given in \eqref{eq:HessenbergLength}, and that all of the intersections $X_w^{\circ} \cap \Hess(S,h)$ form a paving by affines of $\Hess(S,h)$.

\begin{example}\label{ex: Poin for perm}
Let $n=3$ and $h=(2,3,3)$. 
Since $\ell_h(123)=0$, $\ell_h(213)=\ell_h(132)=\ell_h(231)=\ell_h(312)=1$, $\ell_h(321)=2$ for each permutation in $S_3$, the Poincar\'e polynomial of $\Hess(S,h)$ is given by
\begin{align*}
\Poin(\Hess(S,h),\q)=1+4\q+\q^2.
\end{align*}
\end{example}

\bigskip

\section{Cohomology}\label{sect: cohomology}

In this section, we explain the structures of the cohomology rings of Hessenberg varieties, focusing on regular nilpotent Hessenberg varieties $\Hess(N,h)$ in Section \ref{subsec: coh of reg nilp}  and regular semisimple Hessenberg varieties $\Hess(S,h)$ in Section \ref{subsec: coh of reg ss}. We will also see that these two cohomology rings have an interesting relation in Section \ref{subsec: nil vs semi}.

\subsection{Cohomology rings of regular nilpotent Hessenberg varieties}\label{subsec: coh of reg nilp}

The cohomology ring of a regular nilpotent Hessenberg variety $\Hess(N,h)$ has been studied from various viewpoints (e.g. \cite{br-ca04}, \cite{ha-ty}, \cite{mb-ty13}, \cite{dr2015}, \cite{fu-ha-ma}, \cite{ha-ho-ma}, \cite{Insko}, \cite{Insko-Tymoczko}, \cite{AHHM}, \cite{AHMMS}).

In this section we explain an explicit presentation of the cohomology ring of a regular nilpotent Hessenberg variety given by \cite{AHHM} due to M. Harada, M. Masuda, and the authors.
We also discuss a relation between this presentation and Schubert polynomials along \cite{Hori}. 

We first recall an explicit presentation of the cohomology ring of the flag variety $Fl(\C^n)$. 
Let $E_i$ be the $i$-th tautological vector bundle over $Fl(\C^n)$; namely, $E_i$ is the subbundle of the trivial vector bundle $Fl(\C^n) \times \C^n$ over $Fl(\C^n)$ whose fiber over a point $V_\bullet =(V_1 \subset \cdots \subset V_n=\C^n) \in Fl(\C^n)$ is exactly $V_i$.
We denote the negative of the first Chern class of the quotient line bundle $E_i/E_{i-1}$ by $\bar x_i$, i.e.\ 
\begin{equation} \label{eq:ChernTautological}
\bar x_i:=-c_1(E_i/E_{i-1}) \in H^2(Fl(\C^n);\Q).
\end{equation}
We will also use the same symbol $\bar x_i$ for its restrictions to the cohomology of Hessenberg varieties by abuse of notation.
It is known that there is a ring isomorphism 
\begin{equation} \label{eq:CohomologyFl}
H^*(Fl(\C^n);\Q) \cong \Q[x_1,\ldots,x_n]/(e_1,\ldots,e_n)
\end{equation}
which sends $x_i$ to $\bar x_i$ where $e_i$ is the elementary symmetric polynomial of degree $i$ in the variables $x_1,\ldots,x_n$ (cf. \cite[p161, Proposition~3]{fult97}).

In order to describe the cohomology ring of a regular nilpotent Hessenberg variety $\Hess(N,h)$, we define polynomials $f_{i,j}$ for $1 \leq j \leq i \leq n$ as follows:
\begin{equation} \label{eq:f_{i,j}}
f_{i,j}:=\sum_{k=1}^{j} \big( \prod_{\ell=j+1}^i (x_k-x_{\ell})\big)x_k.
\end{equation}
Here, we take by convention $\prod_{\ell=j+1}^i (x_k-x_{\ell})=1$ whenever $i=j$. Note that this definition does not involve $n$.

\begin{theorem}[\cite{AHHM}] \label{theorem:AHHM}
The restriction map 
$$
H^*(Fl(\C^n);\Q) \to H^*(\Hess(N,h);\Q)
$$
is surjective, and there is a ring isomorphism
\begin{equation} \label{eq:AHHM} 
H^*(\Hess(N,h);\Q) \cong \Q[x_1,\ldots,x_n]/(f_{h(1),1},f_{h(2),2},\ldots,f_{h(n),n})
\end{equation}
which sends $x_i$ to $\bar x_i=-c_1(E_i/E_{i-1})|_{\Hess(N,h)}$.
\end{theorem}

\begin{remark}
The presentation \eqref{eq:AHHM} does not hold for the integral coefficients. See \cite[Remark~3]{AHHM2}.
\end{remark}

From the presentation $\eqref{eq:AHHM}$, we can see that the cohomology ring $H^*(\Hess(N,h);\Q)$ is a complete intersection since the number of generators of the polynomial ring $\Q[x_1,\ldots,x_n]$ is equal to the number of generators of the ideal $(f_{h(1),1},\ldots,f_{h(n),n})$.
This implies the following corollary.

\begin{corollary}[\cite{AHHM}] \label{corollary:Hess(N,h)PDA}
$H^*(\Hess(N,h);\Q)$ is a Poincar\'e duality algebra.
\end{corollary}

Note that a regular nilpotent Hessenberg variety is singular in general, but its cohomology is a Poincar\'e duality algebra. 
For arbitrary Lie type, the restriction map is surjective and Corollary~\ref{corollary:Hess(N,h)PDA} holds (\cite{AHMMS}).

\begin{example} \label{example:n=3}
Let $n=3$. We first assign polynomials $f_{i,j}$ to each box below the diagonal line.
\begin{figure}[h]
\begin{center}
\begin{picture}(350,70)
\put(0,10){\framebox(20,20){$f_{3,1}$}}
\put(20,10){\framebox(20,20){$f_{3,2}$}}
\put(40,10){\framebox(20,20){$f_{3,3}$}}
\put(0,30){\framebox(20,20){$f_{2,1}$}}
\put(20,30){\framebox(20,20){$f_{2,2}$}}
\put(40,30){\framebox(20,20)}
\put(0,50){\framebox(20,20){$f_{1,1}$}}
\put(20,50){\framebox(20,20)}
\put(40,50){\framebox(20,20)}

\put(90,55){$f_{1,1}=x_1$, \ \ \ $f_{2,2}=x_1+x_2$, \ \ \ $f_{3,3}=x_1+x_2+x_3$,}
\put(90,35){$f_{2,1}=(x_1-x_2)x_1$, \ \ \ $f_{3,2}=(x_1-x_3)x_1+(x_2-x_3)x_2$,}
\put(90,15){$f_{3,1}=(x_1-x_2)(x_1-x_3)x_1$.}

\put(0,10){\framebox(20,20){$f_{3,1}$}}
\put(20,10){\framebox(20,20){$f_{3,2}$}}
\put(40,10){\framebox(20,20){$f_{3,3}$}}
\put(0,30){\framebox(20,20){$f_{2,1}$}}
\put(20,30){\framebox(20,20){$f_{2,2}$}}
\put(40,30){\framebox(20,20)}
\put(0,50){\framebox(20,20){$f_{1,1}$}}
\put(20,50){\framebox(20,20)}
\put(40,50){\framebox(20,20)}

\put(0,10){\framebox(20,20){$f_{3,1}$}}
\put(20,10){\framebox(20,20){$f_{3,2}$}}
\put(40,10){\framebox(20,20){$f_{3,3}$}}
\put(0,30){\framebox(20,20){$f_{2,1}$}}
\put(20,30){\framebox(20,20){$f_{2,2}$}}
\put(40,30){\framebox(20,20)}
\put(0,50){\framebox(20,20){$f_{1,1}$}}
\put(20,50){\framebox(20,20)}
\put(40,50){\framebox(20,20)}
\end{picture}
\end{center}
\vspace{-10pt}
\caption{The polynomials $f_{i,j}$ for $1\leq i\leq j\leq 3$.}
\label{picture:the polynomials $f_{i,j}$ for $n=3$}
\end{figure}

\vspace{-5pt}
If we take $h=(3,3,3)$, then we obtain from Theorem~\ref{theorem:AHHM} an explicit presentation of the flag variety $\Hess(N,h)=Fl(\C^3)$:
\begin{align}\label{eq: cohomology 333}
H^*(\Hess(N,h);\Q) \cong \mathbb{Q}[x_1,x_2,x_3]/(f_{3,1},f_{3,2},f_{3,3}),
\end{align}
where the polynomials $f_{3,1},f_{3,2},f_{3,3}$ in this presentation are obtained by taking the bottom one in each column (See Figure \ref{picture:$h=(3,3,3)$ and $h'=(2,3,3)$}). It is straightforward to verify that the ideals in \eqref{eq:CohomologyFl} (with $n=3$) and \eqref{eq: cohomology 333} are the same.
If we take $h'=(2,3,3)$, then we also obtain an explicit presentation of the Peterson variety $\Hess(N,h')$ in $Fl(\C^3)$: 
\begin{align}\label{eq: cohomology 233}
H^*(\Hess(N,h');\Q) \cong \mathbb{Q}[x_1,x_2,x_3]/(f_{2,1},f_{3,2},f_{3,3}).
\end{align}

\vspace{-5pt}

\begin{figure}[h]
\begin{center}
\begin{picture}(150,70)
\put(0,10){\framebox(20,20){$f_{3,1}$}}
\put(20,10){\framebox(20,20){$f_{3,2}$}}
\put(40,10){\framebox(20,20){$f_{3,3}$}}
\put(0,30){\framebox(20,20){$f_{2,1}$}}
\put(20,30){\framebox(20,20){$f_{2,2}$}}
\put(40,30){\framebox(20,20)}
\put(0,50){\framebox(20,20){$f_{1,1}$}}
\put(20,50){\framebox(20,20)}
\put(40,50){\framebox(20,20)}

\put(0,53){\colorbox{gray}}
\put(0,58){\colorbox{gray}}
\put(0,62){\colorbox{gray}}
\put(0,67){\colorbox{gray}}
\put(4,53){\colorbox{gray}}
\put(4,58){\colorbox{gray}}
\put(4,62){\colorbox{gray}}
\put(4,67){\colorbox{gray}}
\put(9,53){\colorbox{gray}}
\put(9,58){\colorbox{gray}}
\put(9,62){\colorbox{gray}}
\put(9,67){\colorbox{gray}}
\put(14,53){\colorbox{gray}}
\put(14,58){\colorbox{gray}}
\put(14,62){\colorbox{gray}}
\put(14,67){\colorbox{gray}}

\put(20,53){\colorbox{gray}}
\put(20,58){\colorbox{gray}}
\put(20,62){\colorbox{gray}}
\put(20,67){\colorbox{gray}}
\put(24,53){\colorbox{gray}}
\put(24,58){\colorbox{gray}}
\put(24,62){\colorbox{gray}}
\put(24,67){\colorbox{gray}}
\put(29,53){\colorbox{gray}}
\put(29,58){\colorbox{gray}}
\put(29,62){\colorbox{gray}}
\put(29,67){\colorbox{gray}}
\put(34,53){\colorbox{gray}}
\put(34,58){\colorbox{gray}}
\put(34,62){\colorbox{gray}}
\put(34,67){\colorbox{gray}}

\put(40,53){\colorbox{gray}}
\put(40,58){\colorbox{gray}}
\put(40,62){\colorbox{gray}}
\put(40,67){\colorbox{gray}}
\put(44,53){\colorbox{gray}}
\put(44,58){\colorbox{gray}}
\put(44,62){\colorbox{gray}}
\put(44,67){\colorbox{gray}}
\put(49,53){\colorbox{gray}}
\put(49,58){\colorbox{gray}}
\put(49,62){\colorbox{gray}}
\put(49,67){\colorbox{gray}}
\put(54,53){\colorbox{gray}}
\put(54,58){\colorbox{gray}}
\put(54,62){\colorbox{gray}}
\put(54,67){\colorbox{gray}}

\put(0,33){\colorbox{gray}}
\put(0,38){\colorbox{gray}}
\put(0,42){\colorbox{gray}}
\put(0,47){\colorbox{gray}}
\put(4,33){\colorbox{gray}}
\put(4,38){\colorbox{gray}}
\put(4,42){\colorbox{gray}}
\put(4,47){\colorbox{gray}}
\put(9,33){\colorbox{gray}}
\put(9,38){\colorbox{gray}}
\put(9,42){\colorbox{gray}}
\put(9,47){\colorbox{gray}}
\put(14,33){\colorbox{gray}}
\put(14,38){\colorbox{gray}}
\put(14,42){\colorbox{gray}}
\put(14,47){\colorbox{gray}}

\put(20,33){\colorbox{gray}}
\put(20,38){\colorbox{gray}}
\put(20,42){\colorbox{gray}}
\put(20,47){\colorbox{gray}}
\put(24,33){\colorbox{gray}}
\put(24,38){\colorbox{gray}}
\put(24,42){\colorbox{gray}}
\put(24,47){\colorbox{gray}}
\put(29,33){\colorbox{gray}}
\put(29,38){\colorbox{gray}}
\put(29,42){\colorbox{gray}}
\put(29,47){\colorbox{gray}}
\put(34,33){\colorbox{gray}}
\put(34,38){\colorbox{gray}}
\put(34,42){\colorbox{gray}}
\put(34,47){\colorbox{gray}}

\put(40,33){\colorbox{gray}}
\put(40,38){\colorbox{gray}}
\put(40,42){\colorbox{gray}}
\put(40,47){\colorbox{gray}}
\put(44,33){\colorbox{gray}}
\put(44,38){\colorbox{gray}}
\put(44,42){\colorbox{gray}}
\put(44,47){\colorbox{gray}}
\put(49,33){\colorbox{gray}}
\put(49,38){\colorbox{gray}}
\put(49,42){\colorbox{gray}}
\put(49,47){\colorbox{gray}}
\put(54,33){\colorbox{gray}}
\put(54,38){\colorbox{gray}}
\put(54,42){\colorbox{gray}}
\put(54,47){\colorbox{gray}}

\put(0,13){\colorbox{gray}}
\put(0,18){\colorbox{gray}}
\put(0,22){\colorbox{gray}}
\put(0,27){\colorbox{gray}}
\put(4,13){\colorbox{gray}}
\put(4,18){\colorbox{gray}}
\put(4,22){\colorbox{gray}}
\put(4,27){\colorbox{gray}}
\put(9,13){\colorbox{gray}}
\put(9,18){\colorbox{gray}}
\put(9,22){\colorbox{gray}}
\put(9,27){\colorbox{gray}}
\put(14,13){\colorbox{gray}}
\put(14,18){\colorbox{gray}}
\put(14,22){\colorbox{gray}}
\put(14,27){\colorbox{gray}}

\put(20,13){\colorbox{gray}}
\put(20,18){\colorbox{gray}}
\put(20,22){\colorbox{gray}}
\put(20,27){\colorbox{gray}}
\put(24,13){\colorbox{gray}}
\put(24,18){\colorbox{gray}}
\put(24,22){\colorbox{gray}}
\put(24,27){\colorbox{gray}}
\put(29,13){\colorbox{gray}}
\put(29,18){\colorbox{gray}}
\put(29,22){\colorbox{gray}}
\put(29,27){\colorbox{gray}}
\put(34,13){\colorbox{gray}}
\put(34,18){\colorbox{gray}}
\put(34,22){\colorbox{gray}}
\put(34,27){\colorbox{gray}}

\put(40,13){\colorbox{gray}}
\put(40,18){\colorbox{gray}}
\put(40,22){\colorbox{gray}}
\put(40,27){\colorbox{gray}}
\put(44,13){\colorbox{gray}}
\put(44,18){\colorbox{gray}}
\put(44,22){\colorbox{gray}}
\put(44,27){\colorbox{gray}}
\put(49,13){\colorbox{gray}}
\put(49,18){\colorbox{gray}}
\put(49,22){\colorbox{gray}}
\put(49,27){\colorbox{gray}}
\put(54,13){\colorbox{gray}}
\put(54,18){\colorbox{gray}}
\put(54,22){\colorbox{gray}}
\put(54,27){\colorbox{gray}}

\put(0,10){\framebox(20,20){$f_{3,1}$}}
\put(20,10){\framebox(20,20){$f_{3,2}$}}
\put(40,10){\framebox(20,20){$f_{3,3}$}}
\put(0,30){\framebox(20,20){$f_{2,1}$}}
\put(20,30){\framebox(20,20){$f_{2,2}$}}
\put(40,30){\framebox(20,20)}
\put(0,50){\framebox(20,20){$f_{1,1}$}}
\put(20,50){\framebox(20,20)}
\put(40,50){\framebox(20,20)}
\put(30,20){\circle{17}}
\put(50,20){\circle{17}}
\put(10,20){\circle{17}}
\put(0,-5){$h=(3,3,3)$}

\put(90,10){\framebox(20,20){$f_{3,1}$}}
\put(110,10){\framebox(20,20){$f_{3,2}$}}
\put(130,10){\framebox(20,20){$f_{3,3}$}}
\put(90,30){\framebox(20,20){$f_{2,1}$}}
\put(110,30){\framebox(20,20){$f_{2,2}$}}
\put(130,30){\framebox(20,20)}
\put(90,50){\framebox(20,20){$f_{1,1}$}}
\put(110,50){\framebox(20,20)}
\put(130,50){\framebox(20,20)}

\put(90,53){\colorbox{gray}}
\put(90,58){\colorbox{gray}}
\put(90,62){\colorbox{gray}}
\put(90,67){\colorbox{gray}}
\put(94,53){\colorbox{gray}}
\put(94,58){\colorbox{gray}}
\put(94,62){\colorbox{gray}}
\put(94,67){\colorbox{gray}}
\put(99,53){\colorbox{gray}}
\put(99,58){\colorbox{gray}}
\put(99,62){\colorbox{gray}}
\put(99,67){\colorbox{gray}}
\put(104,53){\colorbox{gray}}
\put(104,58){\colorbox{gray}}
\put(104,62){\colorbox{gray}}
\put(104,67){\colorbox{gray}}

\put(110,53){\colorbox{gray}}
\put(110,58){\colorbox{gray}}
\put(110,62){\colorbox{gray}}
\put(110,67){\colorbox{gray}}
\put(114,53){\colorbox{gray}}
\put(114,58){\colorbox{gray}}
\put(114,62){\colorbox{gray}}
\put(114,67){\colorbox{gray}}
\put(119,53){\colorbox{gray}}
\put(119,58){\colorbox{gray}}
\put(119,62){\colorbox{gray}}
\put(119,67){\colorbox{gray}}
\put(124,53){\colorbox{gray}}
\put(124,58){\colorbox{gray}}
\put(124,62){\colorbox{gray}}
\put(124,67){\colorbox{gray}}

\put(130,53){\colorbox{gray}}
\put(130,58){\colorbox{gray}}
\put(130,62){\colorbox{gray}}
\put(130,67){\colorbox{gray}}
\put(134,53){\colorbox{gray}}
\put(134,58){\colorbox{gray}}
\put(134,62){\colorbox{gray}}
\put(134,67){\colorbox{gray}}
\put(139,53){\colorbox{gray}}
\put(139,58){\colorbox{gray}}
\put(139,62){\colorbox{gray}}
\put(139,67){\colorbox{gray}}
\put(144,53){\colorbox{gray}}
\put(144,58){\colorbox{gray}}
\put(144,62){\colorbox{gray}}
\put(144,67){\colorbox{gray}}

\put(90,33){\colorbox{gray}}
\put(90,38){\colorbox{gray}}
\put(90,42){\colorbox{gray}}
\put(90,47){\colorbox{gray}}
\put(94,33){\colorbox{gray}}
\put(94,38){\colorbox{gray}}
\put(94,42){\colorbox{gray}}
\put(94,47){\colorbox{gray}}
\put(99,33){\colorbox{gray}}
\put(99,38){\colorbox{gray}}
\put(99,42){\colorbox{gray}}
\put(99,47){\colorbox{gray}}
\put(104,33){\colorbox{gray}}
\put(104,38){\colorbox{gray}}
\put(104,42){\colorbox{gray}}
\put(104,47){\colorbox{gray}}

\put(110,33){\colorbox{gray}}
\put(110,38){\colorbox{gray}}
\put(110,42){\colorbox{gray}}
\put(110,47){\colorbox{gray}}
\put(114,33){\colorbox{gray}}
\put(114,38){\colorbox{gray}}
\put(114,42){\colorbox{gray}}
\put(114,47){\colorbox{gray}}
\put(119,33){\colorbox{gray}}
\put(119,38){\colorbox{gray}}
\put(119,42){\colorbox{gray}}
\put(119,47){\colorbox{gray}}
\put(124,33){\colorbox{gray}}
\put(124,38){\colorbox{gray}}
\put(124,42){\colorbox{gray}}
\put(124,47){\colorbox{gray}}

\put(130,33){\colorbox{gray}}
\put(130,38){\colorbox{gray}}
\put(130,42){\colorbox{gray}}
\put(130,47){\colorbox{gray}}
\put(134,33){\colorbox{gray}}
\put(134,38){\colorbox{gray}}
\put(134,42){\colorbox{gray}}
\put(134,47){\colorbox{gray}}
\put(139,33){\colorbox{gray}}
\put(139,38){\colorbox{gray}}
\put(139,42){\colorbox{gray}}
\put(139,47){\colorbox{gray}}
\put(144,33){\colorbox{gray}}
\put(144,38){\colorbox{gray}}
\put(144,42){\colorbox{gray}}
\put(144,47){\colorbox{gray}}

\put(110,13){\colorbox{gray}}
\put(110,18){\colorbox{gray}}
\put(110,22){\colorbox{gray}}
\put(110,27){\colorbox{gray}}
\put(114,13){\colorbox{gray}}
\put(114,18){\colorbox{gray}}
\put(114,22){\colorbox{gray}}
\put(114,27){\colorbox{gray}}
\put(119,13){\colorbox{gray}}
\put(119,18){\colorbox{gray}}
\put(119,22){\colorbox{gray}}
\put(119,27){\colorbox{gray}}
\put(124,13){\colorbox{gray}}
\put(124,18){\colorbox{gray}}
\put(124,22){\colorbox{gray}}
\put(124,27){\colorbox{gray}}

\put(130,13){\colorbox{gray}}
\put(130,18){\colorbox{gray}}
\put(130,22){\colorbox{gray}}
\put(130,27){\colorbox{gray}}
\put(134,13){\colorbox{gray}}
\put(134,18){\colorbox{gray}}
\put(134,22){\colorbox{gray}}
\put(134,27){\colorbox{gray}}
\put(139,13){\colorbox{gray}}
\put(139,18){\colorbox{gray}}
\put(139,22){\colorbox{gray}}
\put(139,27){\colorbox{gray}}
\put(144,13){\colorbox{gray}}
\put(144,18){\colorbox{gray}}
\put(144,22){\colorbox{gray}}
\put(144,27){\colorbox{gray}}

\put(90,10){\framebox(20,20){$f_{3,1}$}}
\put(110,10){\framebox(20,20){$f_{3,2}$}}
\put(130,10){\framebox(20,20){$f_{3,3}$}}
\put(90,30){\framebox(20,20){$f_{2,1}$}}
\put(110,30){\framebox(20,20){$f_{2,2}$}}
\put(130,30){\framebox(20,20)}
\put(90,50){\framebox(20,20){$f_{1,1}$}}
\put(110,50){\framebox(20,20)}
\put(130,50){\framebox(20,20)}
\put(120,20){\circle{17}}
\put(140,20){\circle{17}}
\put(100,40){\circle{17}}
\put(90,-5){$h'=(2,3,3)$}

\end{picture}
\end{center}
\caption{The bottom $f_{i,j}$'s for $h=(3,3,3)$ and $h'=(2,3,3)$.}
\label{picture:$h=(3,3,3)$ and $h'=(2,3,3)$}
\end{figure}

\newpage
One may notice that the polynomial $f_{3,1}$ in \eqref{eq: cohomology 333} is replaced by $f_{2,1}$ in \eqref{eq: cohomology 233} when we changed the Hessenberg function from $h$ to $h'$. 
More specifically, the polynomial $f_{2,1}$ does not vanish in $H^*(\Hess(N,h))$, but it does vanish in $H^*(\Hess(N,h'))$. 
This polynomial $f_{2,1}$ has the following expression as a linear combination of Schubert polynomials $\mathfrak{S}_w \ (w \in S_n)$:
\begin{equation} \label{eq:f2,1}
f_{2,1}=x_1^2-x_1x_2=\mathfrak{S}_{312}-\mathfrak{S}_{231}
\end{equation}
since $\mathfrak{S}_{312}=x_1^2$ and $\mathfrak{S}_{231}=x_1x_2$. As seen in Example~\ref{example:Pet3}, we also have
$$
\Hess(N,h') \cap X_w^{\circ} = \emptyset \iff w=312, \ 231.
$$
These permutations $w=312, 231$ are exactly the ones that
appeared in \eqref{eq:f2,1}.
\end{example}

In general, we have a similar interpretation of the presentation \eqref{eq:AHHM} in Theorem \ref{theorem:AHHM} by considering a smaller Hessenberg variety $\Hess(N,h)\supset \Hess(N,h')$ of codimension $1$, as suggested by the above example.\\

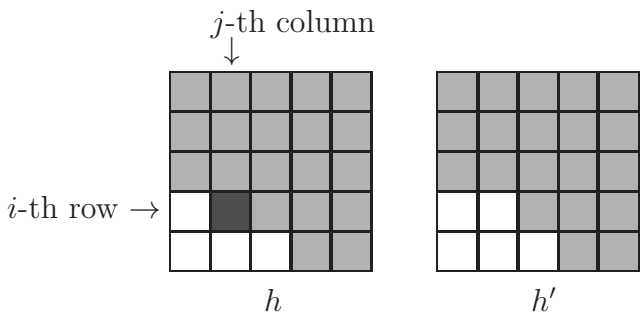
\begin{figure}[h]
\begin{center}
\begin{picture}(175,90)
\put(0,63){\colorbox{gray}}
\put(0,67){\colorbox{gray}}
\put(0,72){\colorbox{gray}}
\put(4,63){\colorbox{gray}}
\put(4,67){\colorbox{gray}}
\put(4,72){\colorbox{gray}}
\put(9,63){\colorbox{gray}}
\put(9,67){\colorbox{gray}}
\put(9,72){\colorbox{gray}}

\put(15,63){\colorbox{gray}}
\put(15,67){\colorbox{gray}}
\put(15,72){\colorbox{gray}}
\put(19,63){\colorbox{gray}}
\put(19,67){\colorbox{gray}}
\put(19,72){\colorbox{gray}}
\put(24,63){\colorbox{gray}}
\put(24,67){\colorbox{gray}}
\put(24,72){\colorbox{gray}}

\put(30,63){\colorbox{gray}}
\put(30,67){\colorbox{gray}}
\put(30,72){\colorbox{gray}}
\put(34,63){\colorbox{gray}}
\put(34,67){\colorbox{gray}}
\put(34,72){\colorbox{gray}}
\put(39,63){\colorbox{gray}}
\put(39,67){\colorbox{gray}}
\put(39,72){\colorbox{gray}}

\put(45,63){\colorbox{gray}}
\put(45,67){\colorbox{gray}}
\put(45,72){\colorbox{gray}}
\put(49,63){\colorbox{gray}}
\put(49,67){\colorbox{gray}}
\put(49,72){\colorbox{gray}}
\put(54,63){\colorbox{gray}}
\put(54,67){\colorbox{gray}}
\put(54,72){\colorbox{gray}}

\put(60,63){\colorbox{gray}}
\put(60,67){\colorbox{gray}}
\put(60,72){\colorbox{gray}}
\put(64,63){\colorbox{gray}}
\put(64,67){\colorbox{gray}}
\put(64,72){\colorbox{gray}}
\put(69,63){\colorbox{gray}}
\put(69,67){\colorbox{gray}}
\put(69,72){\colorbox{gray}}

\put(0,48){\colorbox{gray}}
\put(0,52){\colorbox{gray}}
\put(0,57){\colorbox{gray}}
\put(4,48){\colorbox{gray}}
\put(4,52){\colorbox{gray}}
\put(4,57){\colorbox{gray}}
\put(9,48){\colorbox{gray}}
\put(9,52){\colorbox{gray}}
\put(9,57){\colorbox{gray}}

\put(15,48){\colorbox{gray}}
\put(15,52){\colorbox{gray}}
\put(15,57){\colorbox{gray}}
\put(19,48){\colorbox{gray}}
\put(19,52){\colorbox{gray}}
\put(19,57){\colorbox{gray}}
\put(24,48){\colorbox{gray}}
\put(24,52){\colorbox{gray}}
\put(24,57){\colorbox{gray}}

\put(30,48){\colorbox{gray}}
\put(30,52){\colorbox{gray}}
\put(30,57){\colorbox{gray}}
\put(34,48){\colorbox{gray}}
\put(34,52){\colorbox{gray}}
\put(34,57){\colorbox{gray}}
\put(39,48){\colorbox{gray}}
\put(39,52){\colorbox{gray}}
\put(39,57){\colorbox{gray}}

\put(45,48){\colorbox{gray}}
\put(45,52){\colorbox{gray}}
\put(45,57){\colorbox{gray}}
\put(49,48){\colorbox{gray}}
\put(49,52){\colorbox{gray}}
\put(49,57){\colorbox{gray}}
\put(54,48){\colorbox{gray}}
\put(54,52){\colorbox{gray}}
\put(54,57){\colorbox{gray}}

\put(60,48){\colorbox{gray}}
\put(60,52){\colorbox{gray}}
\put(60,57){\colorbox{gray}}
\put(64,48){\colorbox{gray}}
\put(64,52){\colorbox{gray}}
\put(64,57){\colorbox{gray}}
\put(69,48){\colorbox{gray}}
\put(69,52){\colorbox{gray}}
\put(69,57){\colorbox{gray}}

\put(0,33){\colorbox{gray}}
\put(0,37){\colorbox{gray}}
\put(0,42){\colorbox{gray}}
\put(4,33){\colorbox{gray}}
\put(4,37){\colorbox{gray}}
\put(4,42){\colorbox{gray}}
\put(9,33){\colorbox{gray}}
\put(9,37){\colorbox{gray}}
\put(9,42){\colorbox{gray}}

\put(15,33){\colorbox{gray}}
\put(15,37){\colorbox{gray}}
\put(15,42){\colorbox{gray}}
\put(19,33){\colorbox{gray}}
\put(19,37){\colorbox{gray}}
\put(19,42){\colorbox{gray}}
\put(24,33){\colorbox{gray}}
\put(24,37){\colorbox{gray}}
\put(24,42){\colorbox{gray}}

\put(30,33){\colorbox{gray}}
\put(30,37){\colorbox{gray}}
\put(30,42){\colorbox{gray}}
\put(34,33){\colorbox{gray}}
\put(34,37){\colorbox{gray}}
\put(34,42){\colorbox{gray}}
\put(39,33){\colorbox{gray}}
\put(39,37){\colorbox{gray}}
\put(39,42){\colorbox{gray}}

\put(45,33){\colorbox{gray}}
\put(45,37){\colorbox{gray}}
\put(45,42){\colorbox{gray}}
\put(49,33){\colorbox{gray}}
\put(49,37){\colorbox{gray}}
\put(49,42){\colorbox{gray}}
\put(54,33){\colorbox{gray}}
\put(54,37){\colorbox{gray}}
\put(54,42){\colorbox{gray}}

\put(60,33){\colorbox{gray}}
\put(60,37){\colorbox{gray}}
\put(60,42){\colorbox{gray}}
\put(64,33){\colorbox{gray}}
\put(64,37){\colorbox{gray}}
\put(64,42){\colorbox{gray}}
\put(69,33){\colorbox{gray}}
\put(69,37){\colorbox{gray}}
\put(69,42){\colorbox{gray}}

\put(15,18){\colorbox{Gray}}
\put(15,22){\colorbox{Gray}}
\put(15,27){\colorbox{Gray}}
\put(19,18){\colorbox{Gray}}
\put(19,22){\colorbox{Gray}}
\put(19,27){\colorbox{Gray}}
\put(24,18){\colorbox{Gray}}
\put(24,22){\colorbox{Gray}}
\put(24,27){\colorbox{Gray}}

\put(30,18){\colorbox{gray}}
\put(30,22){\colorbox{gray}}
\put(30,27){\colorbox{gray}}
\put(34,18){\colorbox{gray}}
\put(34,22){\colorbox{gray}}
\put(34,27){\colorbox{gray}}
\put(39,18){\colorbox{gray}}
\put(39,22){\colorbox{gray}}
\put(39,27){\colorbox{gray}}

\put(45,18){\colorbox{gray}}
\put(45,22){\colorbox{gray}}
\put(45,27){\colorbox{gray}}
\put(49,18){\colorbox{gray}}
\put(49,22){\colorbox{gray}}
\put(49,27){\colorbox{gray}}
\put(54,18){\colorbox{gray}}
\put(54,22){\colorbox{gray}}
\put(54,27){\colorbox{gray}}

\put(60,18){\colorbox{gray}}
\put(60,22){\colorbox{gray}}
\put(60,27){\colorbox{gray}}
\put(64,18){\colorbox{gray}}
\put(64,22){\colorbox{gray}}
\put(64,27){\colorbox{gray}}
\put(69,18){\colorbox{gray}}
\put(69,22){\colorbox{gray}}
\put(69,27){\colorbox{gray}}

\put(45,3){\colorbox{gray}}
\put(45,7){\colorbox{gray}}
\put(45,12){\colorbox{gray}}
\put(49,3){\colorbox{gray}}
\put(49,7){\colorbox{gray}}
\put(49,12){\colorbox{gray}}
\put(54,3){\colorbox{gray}}
\put(54,7){\colorbox{gray}}
\put(54,12){\colorbox{gray}}

\put(60,3){\colorbox{gray}}
\put(60,7){\colorbox{gray}}
\put(60,12){\colorbox{gray}}
\put(64,3){\colorbox{gray}}
\put(64,7){\colorbox{gray}}
\put(64,12){\colorbox{gray}}
\put(69,3){\colorbox{gray}}
\put(69,7){\colorbox{gray}}
\put(69,12){\colorbox{gray}}

\put(0,0){\framebox(15,15)}
\put(15,0){\framebox(15,15)}
\put(30,0){\framebox(15,15)}
\put(45,0){\framebox(15,15)}
\put(60,0){\framebox(15,15)}
\put(0,15){\framebox(15,15)}
\put(15,15){\framebox(15,15)}
\put(30,15){\framebox(15,15)}
\put(45,15){\framebox(15,15)}
\put(60,15){\framebox(15,15)}
\put(0,30){\framebox(15,15)}
\put(15,30){\framebox(15,15)}
\put(30,30){\framebox(15,15)}
\put(45,30){\framebox(15,15)}
\put(60,30){\framebox(15,15)}
\put(0,45){\framebox(15,15)}
\put(15,45){\framebox(15,15)}
\put(30,45){\framebox(15,15)}
\put(45,45){\framebox(15,15)}
\put(60,45){\framebox(15,15)}
\put(0,60){\framebox(15,15)}
\put(15,60){\framebox(15,15)}
\put(30,60){\framebox(15,15)}
\put(45,60){\framebox(15,15)}
\put(60,60){\framebox(15,15)}

\put(100,63){\colorbox{gray}}
\put(100,67){\colorbox{gray}}
\put(100,72){\colorbox{gray}}
\put(104,63){\colorbox{gray}}
\put(104,67){\colorbox{gray}}
\put(104,72){\colorbox{gray}}
\put(109,63){\colorbox{gray}}
\put(109,67){\colorbox{gray}}
\put(109,72){\colorbox{gray}}

\put(115,63){\colorbox{gray}}
\put(115,67){\colorbox{gray}}
\put(115,72){\colorbox{gray}}
\put(119,63){\colorbox{gray}}
\put(119,67){\colorbox{gray}}
\put(119,72){\colorbox{gray}}
\put(124,63){\colorbox{gray}}
\put(124,67){\colorbox{gray}}
\put(124,72){\colorbox{gray}}

\put(130,63){\colorbox{gray}}
\put(130,67){\colorbox{gray}}
\put(130,72){\colorbox{gray}}
\put(134,63){\colorbox{gray}}
\put(134,67){\colorbox{gray}}
\put(134,72){\colorbox{gray}}
\put(139,63){\colorbox{gray}}
\put(139,67){\colorbox{gray}}
\put(139,72){\colorbox{gray}}

\put(145,63){\colorbox{gray}}
\put(145,67){\colorbox{gray}}
\put(145,72){\colorbox{gray}}
\put(149,63){\colorbox{gray}}
\put(149,67){\colorbox{gray}}
\put(149,72){\colorbox{gray}}
\put(154,63){\colorbox{gray}}
\put(154,67){\colorbox{gray}}
\put(154,72){\colorbox{gray}}

\put(160,63){\colorbox{gray}}
\put(160,67){\colorbox{gray}}
\put(160,72){\colorbox{gray}}
\put(164,63){\colorbox{gray}}
\put(164,67){\colorbox{gray}}
\put(164,72){\colorbox{gray}}
\put(169,63){\colorbox{gray}}
\put(169,67){\colorbox{gray}}
\put(169,72){\colorbox{gray}}

\put(100,48){\colorbox{gray}}
\put(100,52){\colorbox{gray}}
\put(100,57){\colorbox{gray}}
\put(104,48){\colorbox{gray}}
\put(104,52){\colorbox{gray}}
\put(104,57){\colorbox{gray}}
\put(109,48){\colorbox{gray}}
\put(109,52){\colorbox{gray}}
\put(109,57){\colorbox{gray}}

\put(115,48){\colorbox{gray}}
\put(115,52){\colorbox{gray}}
\put(115,57){\colorbox{gray}}
\put(119,48){\colorbox{gray}}
\put(119,52){\colorbox{gray}}
\put(119,57){\colorbox{gray}}
\put(124,48){\colorbox{gray}}
\put(124,52){\colorbox{gray}}
\put(124,57){\colorbox{gray}}

\put(130,48){\colorbox{gray}}
\put(130,52){\colorbox{gray}}
\put(130,57){\colorbox{gray}}
\put(134,48){\colorbox{gray}}
\put(134,52){\colorbox{gray}}
\put(134,57){\colorbox{gray}}
\put(139,48){\colorbox{gray}}
\put(139,52){\colorbox{gray}}
\put(139,57){\colorbox{gray}}

\put(145,48){\colorbox{gray}}
\put(145,52){\colorbox{gray}}
\put(145,57){\colorbox{gray}}
\put(149,48){\colorbox{gray}}
\put(149,52){\colorbox{gray}}
\put(149,57){\colorbox{gray}}
\put(154,48){\colorbox{gray}}
\put(154,52){\colorbox{gray}}
\put(154,57){\colorbox{gray}}

\put(160,48){\colorbox{gray}}
\put(160,52){\colorbox{gray}}
\put(160,57){\colorbox{gray}}
\put(164,48){\colorbox{gray}}
\put(164,52){\colorbox{gray}}
\put(164,57){\colorbox{gray}}
\put(169,48){\colorbox{gray}}
\put(169,52){\colorbox{gray}}
\put(169,57){\colorbox{gray}}

\put(100,33){\colorbox{gray}}
\put(100,37){\colorbox{gray}}
\put(100,42){\colorbox{gray}}
\put(104,33){\colorbox{gray}}
\put(104,37){\colorbox{gray}}
\put(104,42){\colorbox{gray}}
\put(109,33){\colorbox{gray}}
\put(109,37){\colorbox{gray}}
\put(109,42){\colorbox{gray}}

\put(115,33){\colorbox{gray}}
\put(115,37){\colorbox{gray}}
\put(115,42){\colorbox{gray}}
\put(119,33){\colorbox{gray}}
\put(119,37){\colorbox{gray}}
\put(119,42){\colorbox{gray}}
\put(124,33){\colorbox{gray}}
\put(124,37){\colorbox{gray}}
\put(124,42){\colorbox{gray}}

\put(130,33){\colorbox{gray}}
\put(130,37){\colorbox{gray}}
\put(130,42){\colorbox{gray}}
\put(134,33){\colorbox{gray}}
\put(134,37){\colorbox{gray}}
\put(134,42){\colorbox{gray}}
\put(139,33){\colorbox{gray}}
\put(139,37){\colorbox{gray}}
\put(139,42){\colorbox{gray}}

\put(145,33){\colorbox{gray}}
\put(145,37){\colorbox{gray}}
\put(145,42){\colorbox{gray}}
\put(149,33){\colorbox{gray}}
\put(149,37){\colorbox{gray}}
\put(149,42){\colorbox{gray}}
\put(154,33){\colorbox{gray}}
\put(154,37){\colorbox{gray}}
\put(154,42){\colorbox{gray}}

\put(160,33){\colorbox{gray}}
\put(160,37){\colorbox{gray}}
\put(160,42){\colorbox{gray}}
\put(164,33){\colorbox{gray}}
\put(164,37){\colorbox{gray}}
\put(164,42){\colorbox{gray}}
\put(169,33){\colorbox{gray}}
\put(169,37){\colorbox{gray}}
\put(169,42){\colorbox{gray}}

\put(130,18){\colorbox{gray}}
\put(130,22){\colorbox{gray}}
\put(130,27){\colorbox{gray}}
\put(134,18){\colorbox{gray}}
\put(134,22){\colorbox{gray}}
\put(134,27){\colorbox{gray}}
\put(139,18){\colorbox{gray}}
\put(139,22){\colorbox{gray}}
\put(139,27){\colorbox{gray}}

\put(145,18){\colorbox{gray}}
\put(145,22){\colorbox{gray}}
\put(145,27){\colorbox{gray}}
\put(149,18){\colorbox{gray}}
\put(149,22){\colorbox{gray}}
\put(149,27){\colorbox{gray}}
\put(154,18){\colorbox{gray}}
\put(154,22){\colorbox{gray}}
\put(154,27){\colorbox{gray}}

\put(160,18){\colorbox{gray}}
\put(160,22){\colorbox{gray}}
\put(160,27){\colorbox{gray}}
\put(164,18){\colorbox{gray}}
\put(164,22){\colorbox{gray}}
\put(164,27){\colorbox{gray}}
\put(169,18){\colorbox{gray}}
\put(169,22){\colorbox{gray}}
\put(169,27){\colorbox{gray}}

\put(145,3){\colorbox{gray}}
\put(145,7){\colorbox{gray}}
\put(145,12){\colorbox{gray}}
\put(149,3){\colorbox{gray}}
\put(149,7){\colorbox{gray}}
\put(149,12){\colorbox{gray}}
\put(154,3){\colorbox{gray}}
\put(154,7){\colorbox{gray}}
\put(154,12){\colorbox{gray}}

\put(160,3){\colorbox{gray}}
\put(160,7){\colorbox{gray}}
\put(160,12){\colorbox{gray}}
\put(164,3){\colorbox{gray}}
\put(164,7){\colorbox{gray}}
\put(164,12){\colorbox{gray}}
\put(169,3){\colorbox{gray}}
\put(169,7){\colorbox{gray}}
\put(169,12){\colorbox{gray}}

\put(100,0){\framebox(15,15)}
\put(115,0){\framebox(15,15)}
\put(130,0){\framebox(15,15)}
\put(145,0){\framebox(15,15)}
\put(160,0){\framebox(15,15)}
\put(100,15){\framebox(15,15)}
\put(115,15){\framebox(15,15)}
\put(130,15){\framebox(15,15)}
\put(145,15){\framebox(15,15)}
\put(160,15){\framebox(15,15)}
\put(100,30){\framebox(15,15)}
\put(115,30){\framebox(15,15)}
\put(130,30){\framebox(15,15)}
\put(145,30){\framebox(15,15)}
\put(160,30){\framebox(15,15)}
\put(100,45){\framebox(15,15)}
\put(115,45){\framebox(15,15)}
\put(130,45){\framebox(15,15)}
\put(145,45){\framebox(15,15)}
\put(160,45){\framebox(15,15)}
\put(100,60){\framebox(15,15)}
\put(115,60){\framebox(15,15)}
\put(130,60){\framebox(15,15)}
\put(145,60){\framebox(15,15)}
\put(160,60){\framebox(15,15)}

\put(-65,20){ $i$-th row $\rightarrow$}
\put(20,80){$\downarrow$}
\put(15,90){$j$-th column}

\put(35,-15){$h$}
\put(135,-15){$h'$}
\end{picture}
\end{center}
\vspace{5pt}
\caption{The pictures of $h$ and $h'$.}
\label{picture:h and h'}
\end{figure}
\noindent
The unique difference in the generators of the ideal appearing in the presentation \eqref{eq:AHHM} for $H^*(\Hess(N,h);\Q)$ and $H^*(\Hess(N,h');\Q)$ is the polynomial $f_{i-1,j}$.
The second author showed that this polynomial can be written as an alternating sum of Schubert polynomials $\mathfrak{S}_{w}$ where the set of permutations $w$ appearing in this sum coincides with the set of minimal length permutations $w$ in $S_n$ satisfying
\begin{align*}
\Hess(N,h) \cap X_w^{\circ} \neq \emptyset \ \ \  {\rm and} \ \  \Hess(N,h') \cap X_w^{\circ} = \emptyset.
\end{align*}
See \cite{Hori} for details.

\smallskip

\subsection{The cohomology rings of regular semisimple Hessenberg varieties}\label{sect: regular semisimple}\label{subsec: coh of reg ss}
Let $\Hess(S,h)$ be a regular semisimple Hessenberg variety.
One of the most interesting feature of $\Hess(S,h)$ is the $\Sn$-representation on its cohomology $H^*(\Hess(S,h);\C)$ constructed by J. Tymoczko (\cite{tymo08}). She first constructed an $\Sn$-representation on a torus equivariant cohomology $H^*_T(\Hess(S,h);\C)$ via a combinatorial description called the GKM-presentation, and she showed that this representation descends to the ordinary cohomology $H^*(\Hess(S,h);\C)$. An alternative geometric construction via a monodromy action of the fundamental group of the space of regular semisimple matrices is explained in Brosnan-Chow (\cite{br-ch}). 
Following the construction \cite{tymo08}, N. Teff started to analyze this $\Sn$-representation in \cite{Teff11, Teff13}, and Shareshian-Wachs (\cite{sh-wa11}, \cite{sh-wa14}) announced a beautiful conjecture on this representation using chromatic quasisymmetric functions. In this manuscript, we explain this representation along the construction due to Tymoczko.

Let $\Tn\subset {\rm GL}(n, \C)$ be the maximal torus consisting of diagonal elements of ${\rm GL}(n, \C)$. The flag variety $\Flags(\C^n)$ has a natural action of ${\rm GL}(n, \C)$, and hence the torus $\Tn$ acts on $\Flags(\C^n)$ via its restriction. This $\Tn$-action preserves the regular semisimple Hessenberg variety $\Hess(S,h)$ since all the elements of $\Tn$ commute with the diagonal matrix $S$. It is known that $\Hess(S,h)$ contains all the $\Tn$-fixed points of the flag variety $\Flags(\C^n)$ (\cite[Proposition~3]{ma-pr-sh}) so that
$
\Hess(S,h)^{\Tn}=\Flags(\C^n)^{\Tn}\cong \Sn.
$
Here, the last bijection corresponds $w\in \Sn$ and the permutation flag associated with $w$. 

Let $H^*_{\Tn}(\Hess(S,h);\C)$ be the $\Tn$-equivariant cohomology of $\Hess(S,h)$. Recalling that $\Hess(S,h)$ has no odd-degree cohomology from Theorem \ref{theorem:Tymoczko}, we can apply localization techniques for $T$-equivariant cohomology which we refer \cite{Go-Ko-Ma,tymo08} for details. As a conclusion, we obtain the so-called GKM presentation of $H^*_{\Tn}(\Hess(S,h);\C)$ as a subring of a direct sum of polynomial rings $\bigoplus _{w \in \Sn} \C[t_1,\dots,t_n]$. 

\begin{proposition} $($\cite[Proposition~4.7]{tymo08}$)$\label{proposition:GKM}
The equivariant cohomology $H^*_{\Tn}(\Hess(S,h);\C)$ is isomorphic $($as rings$)$ to
\begin{equation}\label{eq:GKM}
\left\{ \alpha \in \bigoplus_{w\in\Sn} \C[t_1,\dots,t_n] \left|
\begin{matrix} \text{ $\alpha(w)-\alpha(w')$ is divisible by $t_{w(i)}-t_{w(j)}$} \\
\text{ if $w'=w(j \ i)$ for some $j<i$ with $i \leq h(j)$} 
 \end{matrix} \right.\right\}
\end{equation}
where $\alpha(w)$ is the $w$-component of $\alpha$ and $(j \ i)\in S_n$ is the transposition of $j$ and $i$. 
\end{proposition}

We can visualize elements of the subring appearing in \eqref{eq:GKM} in terms of the so-called GKM graph whose vertex set is $\Sn$ and there is an edge between vertices $w, w'\in\Sn$ if there exists $1\leq j<i\leq n$ with $i\leq h(j)$ satisfying $w'=w(j\ i)$. Additionally, we equip such an edge with the data of the polynomial $\pm(t_{w(i)}-t_{w(j)})$ (up to sign) arising in \eqref{eq:GKM}.
This labeled graph is called \textbf{the GKM graph} of $\Hess(S,h)$, and we denote it by $\Gamma(h)$. 
In this language, the condition in \eqref{eq:GKM} says that the collection of polynomials $(\alpha(w))_{w\in \Sn}$ satisfies the following: if $w$ and $w'$ are connected in $\Gamma(h)$ by an edge labeled by $t_{w(i)}-t_{w(j)}$, then the difference of the polynomials assigned for $w$ and $w'$ must be divisible by the label.

\begin{example} \label{example:GKM graph}
Let $n=3$. 
For $h=(3,3,3)$ and $h'=(2,3,3)$, the corresponding GKM graphs are depicted in Figure $\ref{pic:GKM graphs}$, where we use the one-line notation for each vertex.
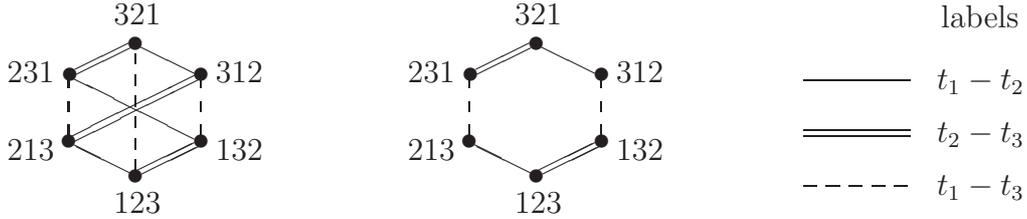
\begin{figure}[h]
\begin{center}
\begin{picture}(380,70)
        \put(50,14){\circle*{5}}
        \put(50,64){\circle*{5}}
        \put(74.5,27){\circle*{5}}
        \put(74.5,52){\circle*{5}}
        \put(25,27){\circle*{5}}
        \put(25.5,52){\circle*{5}}

        \put(72,28){\line(-2,1){48}}
        \put(48,15){\line(-2,1){21}}
        \put(72,52.5){\line(-2,1){20}}
        \put(23.5,27.5){\line(2,1){49}}
        \put(24.5,25.5){\line(2,1){49}}
        \put(51,16){\line(2,1){21}}
        \put(52.5,14.5){\line(2,1){21}}
        \put(27.5,54){\line(2,1){20}}
        \put(28,52){\line(2,1){20}}
        \put(50,16){\line(0,1){4}}
        \put(50,24){\line(0,1){4}}
        \put(50,32){\line(0,1){4}}
        \put(50,40){\line(0,1){4}}
        \put(50,48){\line(0,1){4}}
        \put(50,56){\line(0,1){4}}
        \put(25,29){\line(0,1){4}}
        \put(25,38){\line(0,1){4}}
        \put(25,47){\line(0,1){4}}
        \put(74.3,29){\line(0,1){4}}
        \put(74.3,38){\line(0,1){4}}
        \put(74.3,47){\line(0,1){4}}

        \put(42,0){$123$}
        \put(42,71){$321$}
        \put(80,20){$132$}
        \put(80,49){$312$}
        \put(2,20){$213$}
        \put(2,49){$231$}


        \put(200,14){\circle*{5}}
        \put(200,64){\circle*{5}}
        \put(224.5,27){\circle*{5}}
        \put(224.5,52){\circle*{5}}
        \put(175,27){\circle*{5}}
        \put(175.5,52){\circle*{5}}

        \put(198,15){\line(-2,1){21}}
        \put(222,52.5){\line(-2,1){20}}
        \put(201,16){\line(2,1){21}}
        \put(202.5,14.5){\line(2,1){21}}
        \put(177.5,54){\line(2,1){20}}
        \put(178,52){\line(2,1){20}}
        \put(175,29){\line(0,1){4}}
        \put(175,38){\line(0,1){4}}
        \put(175,47){\line(0,1){4}}
        \put(224.3,29){\line(0,1){4}}
        \put(224.3,38){\line(0,1){4}}
        \put(224.3,47){\line(0,1){4}}

        \put(192,0){$123$}
        \put(192,71){$321$}
        \put(230,20){$132$}
        \put(230,49){$312$}
        \put(152,20){$213$}
        \put(152,49){$231$}


\put(351.5,70){{\rm labels}}        
\put(300,50){\line(1,0){40}}
\put(350,46.5){$t_1-t_2$}
\put(300,31){\line(1,0){40}}
\put(300,29){\line(1,0){40}}
\put(350,26.5){$t_2-t_3$}
\put(300,10){\line(1,0){5}}
\put(308.5,10){\line(1,0){5}}
\put(317,10){\line(1,0){5}}
\put(325.5,10){\line(1,0){5}}
\put(334,10){\line(1,0){5}}
\put(350,6.5){$t_1-t_3$}
\end{picture}
\end{center}
\vspace{-5pt}
\caption{The GKM graphs $\Gamma(h)$ and $\Gamma(h')$.}
\label{pic:GKM graphs}
\end{figure}
\vspace{-5pt}

\noindent
For example, one can verify that the tuples of polynomials in Figure $\ref{pic:element xi}$ are elements of $H^*_{\Tn}(\Hess(S,h);\C)=H^*_{\Tn}(Fl(\C^3);\C)$ (and hence of $H^*_{\Tn}(\Hess(S,h');\C)$). 
\begin{figure}[h]
\begin{center}
\begin{picture}(400,70)
        \put(65,14){\circle*{5}}
        \put(65,64){\circle*{5}}
        \put(89.5,27){\circle*{5}}
        \put(89.5,52){\circle*{5}}
        \put(40,27){\circle*{5}}
        \put(40.5,52){\circle*{5}}

        \put(87,28){\line(-2,1){48}}
        \put(63,15){\line(-2,1){21}}
        \put(87,52.5){\line(-2,1){20}}
        \put(38.5,27.5){\line(2,1){49}}
        \put(39.5,25.5){\line(2,1){49}}
        \put(66,16){\line(2,1){21}}
        \put(67.5,14.5){\line(2,1){21}}
        \put(42.5,54){\line(2,1){20}}
        \put(43,52){\line(2,1){20}}
        \put(65,16){\line(0,1){4}}
        \put(65,24){\line(0,1){4}}
        \put(65,32){\line(0,1){4}}
        \put(65,40){\line(0,1){4}}
        \put(65,48){\line(0,1){4}}
        \put(65,56){\line(0,1){4}}
        \put(40,29){\line(0,1){4}}
        \put(40,38){\line(0,1){4}}
        \put(40,47){\line(0,1){4}}
        \put(89.3,29){\line(0,1){4}}
        \put(89.3,38){\line(0,1){4}}
        \put(89.3,47){\line(0,1){4}}

        \put(62,0){$t_1$}
        \put(62,71){$t_3$}
        \put(95,22){$t_1$}
        \put(95,50){$t_3$}
        \put(27,22){$t_2$}
        \put(27,50){$t_2$}
        \put(-2,35){$\bar x_1=$}   
        \put(115,35){,}


        \put(205,14){\circle*{5}}
        \put(205,64){\circle*{5}}
        \put(229.5,27){\circle*{5}}
        \put(229.5,52){\circle*{5}}
        \put(180,27){\circle*{5}}
        \put(180.5,52){\circle*{5}}

        \put(227,28){\line(-2,1){48}}
        \put(203,15){\line(-2,1){21}}
        \put(227,52.5){\line(-2,1){20}}
        \put(178.5,27.5){\line(2,1){49}}
        \put(179.5,25.5){\line(2,1){49}}
        \put(206,16){\line(2,1){21}}
        \put(207.5,14.5){\line(2,1){21}}
        \put(182.5,54){\line(2,1){20}}
        \put(183,52){\line(2,1){20}}
        \put(205,16){\line(0,1){4}}
        \put(205,24){\line(0,1){4}}
        \put(205,32){\line(0,1){4}}
        \put(205,40){\line(0,1){4}}
        \put(205,48){\line(0,1){4}}
        \put(205,56){\line(0,1){4}}
        \put(180,29){\line(0,1){4}}
        \put(180,38){\line(0,1){4}}
        \put(180,47){\line(0,1){4}}
        \put(229.3,29){\line(0,1){4}}
        \put(229.3,38){\line(0,1){4}}
        \put(229.3,47){\line(0,1){4}}

        \put(202,0){$t_2$}
        \put(202,71){$t_2$}
        \put(235,22){$t_3$}
        \put(235,50){$t_1$}
        \put(166,22){$t_1$}
        \put(166,50){$t_3$}
        \put(137,35){$\bar x_2=$}   
        \put(255,35){,}   


        \put(343,14){\circle*{5}}
        \put(343,64){\circle*{5}}
        \put(367.5,27){\circle*{5}}
        \put(367.5,52){\circle*{5}}
        \put(318,27){\circle*{5}}
        \put(318.5,52){\circle*{5}}

        \put(365,28){\line(-2,1){48}}
        \put(341,15){\line(-2,1){21}}
        \put(365,52.5){\line(-2,1){20}}
        \put(316.5,27.5){\line(2,1){49}}
        \put(317.5,25.5){\line(2,1){49}}
        \put(344,16){\line(2,1){21}}
        \put(345.5,14.5){\line(2,1){21}}
        \put(320.5,54){\line(2,1){20}}
        \put(321,52){\line(2,1){20}}
        \put(343,16){\line(0,1){4}}
        \put(343,24){\line(0,1){4}}
        \put(343,32){\line(0,1){4}}
        \put(343,40){\line(0,1){4}}
        \put(343,48){\line(0,1){4}}
        \put(343,56){\line(0,1){4}}
        \put(318,29){\line(0,1){4}}
        \put(318,38){\line(0,1){4}}
        \put(318,47){\line(0,1){4}}
        \put(367.3,29){\line(0,1){4}}
        \put(367.3,38){\line(0,1){4}}
        \put(367.3,47){\line(0,1){4}}

        \put(340,0){$t_3$}
        \put(340,71){$t_1$}
        \put(373,22){$t_2$}
        \put(373,50){$t_2$}
        \put(305,22){$t_3$}
        \put(305,50){$t_1$}
        \put(277,35){$\bar x_3=$}   

\end{picture}
\end{center}
\caption{Some elements of $H^*_{\Tn}(\Hess(S,h))$.}
\label{pic:element xi}
\end{figure}
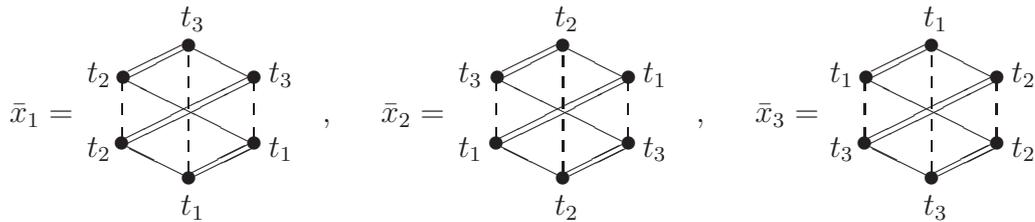

\noindent
We use the symbols $\bar{x}_1, \bar{x}_2, \bar{x}_3$ by abuse of notation (cf.\ \eqref{eq:ChernTautological}) because of the reason which we will explain soon later (See Example \ref{example:Sn-rep for xi} below). Also, the elements $\bar{y}_1, \bar{y}_2, \bar{y}_3$ given by the tuples of polynomials in Figure $\ref{pic:element}$ are elements of $H^*_{\Tn}(\Hess(S,h');\C)$ but not of $H^*_{\Tn}(\Hess(S,h);\C)=H^*_{\Tn}(Fl(\C^3);\C)$. 
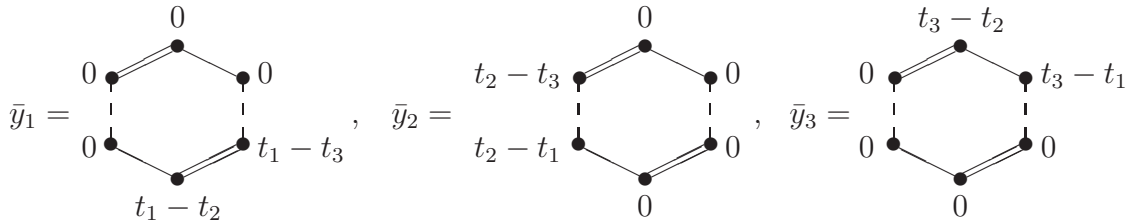
\begin{figure}[h]
\begin{center}
\begin{picture}(400,70)
        \put(50,14){\circle*{5}}
        \put(50,64){\circle*{5}}
        \put(74.5,27){\circle*{5}}
        \put(74.5,52){\circle*{5}}
        \put(25,27){\circle*{5}}
        \put(25.5,52){\circle*{5}}

        \put(48,15){\line(-2,1){21}}
        \put(72,52.5){\line(-2,1){20}}
        \put(51,16){\line(2,1){21}}
        \put(52.5,14.5){\line(2,1){21}}
        \put(27.5,54){\line(2,1){20}}
        \put(28,52){\line(2,1){20}}
        \put(25,29){\line(0,1){4}}
        \put(25,38){\line(0,1){4}}
        \put(25,47){\line(0,1){4}}
        \put(74.3,29){\line(0,1){4}}
        \put(74.3,38){\line(0,1){4}}
        \put(74.3,47){\line(0,1){4}}

        \put(34,0){$t_1-t_2$}
        \put(47,71){$0$}
        \put(80,22){$t_1-t_3$}
        \put(80,50){$0$}
        \put(14,22){$0$}
        \put(14,50){$0$}
        \put(-13,35){$\bar y_1=$}   
        \put(115,35){,}


        \put(225,14){\circle*{5}}
        \put(225,64){\circle*{5}}
        \put(249.5,27){\circle*{5}}
        \put(249.5,52){\circle*{5}}
        \put(200,27){\circle*{5}}
        \put(200.5,52){\circle*{5}}

        \put(223,15){\line(-2,1){21}}
        \put(247,52.5){\line(-2,1){20}}
        \put(226,16){\line(2,1){21}}
        \put(227.5,14.5){\line(2,1){21}}
        \put(202.5,54){\line(2,1){20}}
        \put(203,52){\line(2,1){20}}
        \put(200,29){\line(0,1){4}}
        \put(200,38){\line(0,1){4}}
        \put(200,47){\line(0,1){4}}
        \put(249.3,29){\line(0,1){4}}
        \put(249.3,38){\line(0,1){4}}
        \put(249.3,47){\line(0,1){4}}

        \put(222,0){$0$}
        \put(222,71){$0$}
        \put(255,22){$0$}
        \put(255,50){$0$}
        \put(161,22){$t_2-t_1$}
        \put(161,50){$t_2-t_3$}
        \put(130,35){$\bar y_2=$}   
        \put(265.5,35){,}   


        \put(343,14){\circle*{5}}
        \put(343,64){\circle*{5}}
        \put(367.5,27){\circle*{5}}
        \put(367.5,52){\circle*{5}}
        \put(318,27){\circle*{5}}
        \put(318.5,52){\circle*{5}}

        \put(341,15){\line(-2,1){21}}
        \put(365,52.5){\line(-2,1){20}}
        \put(344,16){\line(2,1){21}}
        \put(345.5,14.5){\line(2,1){21}}
        \put(320.5,54){\line(2,1){20}}
        \put(321,52){\line(2,1){20}}
        \put(318,29){\line(0,1){4}}
        \put(318,38){\line(0,1){4}}
        \put(318,47){\line(0,1){4}}
        \put(367.3,29){\line(0,1){4}}
        \put(367.3,38){\line(0,1){4}}
        \put(367.3,47){\line(0,1){4}}

        \put(340,0){$0$}
        \put(327,71){$t_3-t_2$}
        \put(373,22){$0$}
        \put(373,50){$t_3-t_1$}
        \put(305,22){$0$}
        \put(305,50){$0$}
        \put(279,35){$\bar y_3=$}   

\end{picture}
\end{center}
\vspace{-5pt}
\caption{Some elements of $H^*_{\Tn}(\Hess(S,h'))$.}
\label{pic:element}
\end{figure}
\end{example}

In what follows, we identify $H^*_{\Tn}(\Hess(S,h);\C)$ and the presentation \eqref{eq:GKM}, and we do not distinguish them. For each $i=1,\ldots,n$, it is clear that a collection $(t_i)_{w\in\Sn}$ lies in \eqref{eq:GKM}. For simplicity, we also write this element as $t_i$ by abuse of notation. Then the theory of $T$-equivariant cohomology also shows that there is a ring isomorphism
\begin{align}\label{eq: recovering ordinary cohomology}
H^*(\Hess(S,h);\C) \cong H^*_{\Tn}(\Hess(S,h);\C)/(t_1,\ldots,t_n).
\end{align}
This means that one can study the ordinary cohomology ring $H^*(\Hess(S,h);\C)$ from the equivariant cohomology ring $H^*_{\Tn}(\Hess(S,h);\C)$. 

We now describe the  $\Sn$-action on $H^*_{\Tn}(\Hess(S,h);\C)$ constructed by Tymoczko \cite{tymo08}. 
For $v \in \Sn$ and $\alpha = (\alpha(w))\in \bigoplus _{w \in \Sn} \C[t_1,\dots,t_n]$, we define an element $v \cdot \alpha$ by the formula 
\begin{equation}\label{eq:Tymoczko} 
(v\cdot \alpha)(w) := v \cdot \alpha(v^{-1}w) \ \textup{ for all } w \in \Sn
\end{equation}
where $v \cdot f(t_{1},\ldots,t_{n})=f(t_{v(1)},\ldots,t_{v(n)})$ for $f(t_{1},\ldots,t_{n}) \in \C[t_1,\dots,t_n]$ in the right-hand side. This $\Sn$-action preserves the subset \eqref{eq:GKM}, and hence it defines an $\Sn$-action on the equivariant cohomology $H^*_{\Tn}(\Hess(S,h);\C)$ by Proposition~\ref{proposition:GKM}. 
Since we have $v\cdot t_i=t_{v(i)}$ for the classes
$t_i=(t_i)_{w \in \Sn}$ defined above,
the $\Sn$-action on $H^*_{\Tn}(\Hess(S,h);\C)$ induces an $\Sn$-action on the ordinary cohomology $H^*(\Hess(S,h);\C)$ via \eqref{eq: recovering ordinary cohomology}. By construction this $\Sn$-representation preserves the cup product of $H^*(\Hess(S,h);\C)$. 

\begin{example} \label{example:Sn-rep for xi}
Let $n=3$ and $h=(3,3,3)$.
For $\bar{x}_1,\bar{x}_2,\bar{x}_3\in H^2_{\Tn}(\Hess(S,h);\C)=H^2_{\Tn}(Fl(\C^3);\C)$ given in Figure \ref{pic:element xi}, 
one can easily verify from the definition \eqref{eq:Tymoczko}  that they are invariant under the $S_3$-action; $$w\cdot \bar{x}_i=\bar{x}_{i}\quad (1\leq i\leq 3)$$ for any $w\in S_3$. Under the isomorphism \eqref{eq: recovering ordinary cohomology}, it follows that these equivariant cohomology classes corresponds to $\bar{x}_i=-c_1(E_i/E_{i-1}) \in H^2(Fl(\C^3);\Q)$ introduced in \eqref{eq:ChernTautological} which gives a justification for our notation.
From \eqref{eq:CohomologyFl} (or \eqref{eq: cohomology 333}), this means that the $S_3$-representation on $H^*(Fl(\C^3);\C)$ is trivial, and the same claim holds for the case $H^*(Fl(\C^n);\C)$ in general (\cite[Proposition 4.4]{tymo08}).
\end{example}

\begin{example} \label{example:Sn-rep for yi}
Let $n=3$ and $h=(2,3,3)$.
For $\bar{y}_1,\bar{y}_2,\bar{y}_3\in H^2_{\Tn}(\Hess(S,h);\C)$ given in Figure \ref{pic:element}, 
one can also verify that these classes are naturally permuted by the $S_3$-action;
$$w\cdot \bar{y}_i=\bar{y}_{w(i)}\quad (1\leq i\leq 3)$$ for any $w\in S_3$ from the definition \eqref{eq:Tymoczko}.
\end{example}

Compared to the situation for $\Hess(N,h)$ (e.g.\ Theorem \ref{theorem:AHHM}), the restriction map $H^*(\Flags(\C^n);\C)\rightarrow H^*(\Hess(S,h);\C)$ is not surjective in general. Hence, to describe the ring $H^*(\Hess(S,h);\C)$ in terms of ring generators and relations among them, we need to find some cohomology classes of $\Hess(S,h)$ which do not come from $H^*(\Flags(\C^n);\C)$ by restriction. However, since we have a surjection $H^*_{\Tn}(\Hess(S,h);\C)\rightarrow H^*(\Hess(S,h);\C)$ by \eqref{eq: recovering ordinary cohomology}, the graphical presentation of the equivariant cohomology ring $H^*_{\Tn}(\Hess(S,h);\C)$ can be used to seek those classes as we saw in Example \ref{example:GKM graph}. 
For the case $h=(h(1),n,\ldots,n)$, M. Masuda and the authors showed that we can explicitly describe the integral cohomology ring in terms of ring generators and their relations by this approach. See \cite{ab-ho-ma} for details. 
For the case $h=(2,3,4,\ldots,n,n)$, the cohomology ring $H^*(\Hess(S,h);\C)$ is well-understood since $\Hess(S,h)$ is a non-singular projective toric variety in this case (\cite{Klyachko}, \cite{proc90}, \cite{Abe}).
However, the ring structure of $H^*(\Hess(S,h);\C)$ for general $h$ is not well-understood at this moment.

\begin{example}\label{ex: cohomology of regular semisimple}
Let $n=3$ and $h=(2,3,3)$. 
Recall from Figure \ref{pic:element} that we have three classes $\bar{y}_1,\bar{y}_2,\bar{y}_3$ in the equivariant cohomology $H^2_{\Tn}(\Hess(S,h);\C)$.
We also denote by the same symbol $\bar{y}_i\in H^2(\Hess(S,h);\C)$ the image of $\bar{y}_i$ under the isomorphism \eqref{eq: recovering ordinary cohomology} by abuse of notation.
Then the presentation of $H^*(\Hess(S,h);\C)$ due to \cite{ab-ho-ma} is given by
\begin{align}\label{eq: presentation for toric}
H^*(\Hess(S,h);\C) 
\cong \C[x_1,x_2,x_3,y_1,y_2,y_3]/J, 
\end{align}
where $x_i$ and $y_i$ correspond to $\bar{x}_i=-c_1(E_i/E_{i-1})|_{\Hess(S,h)}$ and $\bar{y}_i$ respectively, and the ideal $J$ is generated by 
\begin{align*}
&y_{k} y_{k'} \quad \text{for } 1\leq k\neq k'\leq 3, \\
&x_1y_k \quad \text{for } 1\leq k\leq 3, \\
&x_3y_k+x_2x_3 \quad \text{for } 1\leq k\leq 3, \\
&y_1+y_2+y_3 - (x_1-x_2), \\
&x_1+x_2+x_3, \ x_1x_2+x_1x_3+x_2x_3, \ x_1x_2x_3.
\end{align*}
Also, the $S_3$-action on $H^*(\Hess(S,h);\C)$ is given by $w\cdot x_i = x_{i}$ and $w\cdot y_i = y_{w(i)}$ for $i=1,2,3$ and $w\in S_3$, and the irreducible decomposition as $S_3$-representation is given by
\begin{align}\label{eq: Sn-rep for toric}
H^*(\Hess(S,h);\C) 
\cong 
S^{(3)} \oplus \big((S^{(3)})^{\oplus 2} \oplus S^{(2,1)} \big)\, \q \oplus S^{(3)}\, \q^2,
\end{align}
where $S^{\lambda}$ for a partition $\lambda$ of $3$ is the irreducible representation of $S_3$ corresponding to $\lambda$ (cf.\ \cite{fult97}), and $\q$ is a formal symbol standing for the cohomology grading with $\deg (\q)=2$. Note that $S^{(3)}$ is the trivial representation. In particular, this recovers the Poincar\'e polynomial of $\Hess(S,h)$ in Example~\ref{ex: Poin for perm}.
\end{example}

\smallskip

\subsection{Regular nilpotent vs. regular semisimple}\label{subsec: nil vs semi}
In the last two sections, we have described the cohomology of regular nilpotent Hessenberg varieties $\Hess(N,h)$ and regular semisimple Hessenberg varieties $\Hess(S,h)$, and we saw that the latter cohomology $H^*(\Hess(S,h);\C)$ admits the $S_n$-representation constructed by Tymoczko.
Since this representation preserves the cup product, the invariant subgroup $H^*(\Hess(S,h);\C)^{\Sn}$ in fact forms a subring of $H^*(\Hess(S,h);\C)$. M. Harada, M. Masuda, and the authors (\cite{AHHM}) showed that the $\Sn$-representation provides a connection between the topology of $\Hess(N,h)$ and $\Hess(S,h)$ as follows.

\begin{theorem} $($\cite{AHHM}$)$\label{theorem:Masuda isomorphism}
Let $\Hess(N,h)$ and $\Hess(S,h)$ be a regular nilpotent Hessenberg variety and a regular semisimple Hessenberg variety, respectively.
Then, there is a ring isomorphism
\begin{align*}
H^*(\Hess(N,h);\C) \cong H^*(\Hess(S,h);\C)^{\Sn}
\end{align*}
which sends $\bar{x}_i=-c_1(E_i/E_{i-1})|_{\Hess(N,h)}$ to $\bar{x}_i=-c_1(E_i/E_{i-1})|_{\Hess(S,h)}$.
\end{theorem}

For the case $h=(n,n,\ldots,n)$, we have $\Hess(N,h)=\Hess(S,h)=Fl(\C^n)$, and the isomorphism in Theorem \ref{theorem:Masuda isomorphism} is obvious since the $S_n$-representation on the cohomology is trivial in this case (See Example \ref{example:Sn-rep for xi}).
For the case $h=(2,3,4,\ldots,n,n)$, 
explicit presentations for the rings $H^*(\Hess(N,h);\C)$ and $H^*(\Hess(S,h);\C)^{\Sn}$ were given in \cite{fu-ha-ma,ha-ho-ma} and \cite{Klyachko} respectively, and those presentations are in fact identical although it was not mentioned. In this case, it means that the cohomology ring of the Peterson variety is isomorphic to the $\Sn$-invariant subring of the cohomology of the permutohedral variety. See the work of P. Brosnan and T. Chow (\cite{br-ch}) for the geometry behind this phenomenon.
We also note that Theorem \ref{theorem:Masuda isomorphism} holds for arbitrary Lie type (\cite{AHMMS}).

\begin{example}
Let $n=3$ and $h=(2,3,3)$. Then from \eqref{eq: Poinc for 233} and \eqref{eq: Sn-rep for toric}, it is clear that we have the following equalities for dimension;
$$
\dim H^{2k}(\Hess(N,h);\C) = \dim H^{2k}(\Hess(S,h);\C)^{S_3}
$$
for all $k$. However, Theorem \ref{theorem:Masuda isomorphism} states more; these are isomorphic as (graded) rings. For this case $h=(2,3,3)$, one can directly construct this isomorphism by using the ring presentations \eqref{eq: cohomology 233} and \eqref{eq: presentation for toric} (see \cite[Remark 4.8.]{ab-ho-ma}).
\end{example}

\vspace{10pt}
\bigskip

\section{Combinatorics}\label{sect: combinatorics}
In this section, we explain combinatorial objects which are related to Hessenberg varieties.
More specifically, we will see how hyperplane arrangements arise to describe the structure of the cohomology rings $H^*(\Hess(N,h);\R)$ of regular nilpotent Hessenberg varieties and how Stanley's chromatic symmetric function of graphs determines the $S_n$-representation on the cohomology rings $H^*(\Hess(S,h);\C)$ of regular semisimple Hessenberg varieties.

\smallskip

\subsection{Hyperplane arrangements} \label{subsect:Hyperplane arrangements}

In this section we explain a connection between Hessenberg varieties and hyperplane arrangements established in \cite{AHMMS}. 
Originally, E. Sommers and J. Tymoczko pointed out that Hessenberg varieties are related to hyperplane arrangements, and they conjectured that the Ponicar\'e polynomial of a regular nilpotent Hessenberg variety $\Hess(N,h)$ can be described in terms of certain hyperplane arrangement (\cite{so-ty}).
This conjecture was verified for some Lie types by Sommers-Tymoczko, G. R\"ohrle, A. Schauenburg (\cite{so-ty}, \cite{Roe}). After this, an explicit presentation of the cohomology ring of $\Hess(N,h)$ was provided by M. Harada, M. Masuda, and the authors (\cite{AHHM}).
Motivated by all these works, 
T. Abe, M. Masuda, S. Murai, T. Sato, and the second author proved that the cohomology ring of a regular nilpotent Hessenberg variety is isomorphic to a ring coming from those hyperplane arrangements (\cite{AHMMS}).
By taking Hilbert series of both sides of this isomorphism, one sees that the conjecture of Sommers and Tymoczko is true. For general reference about hyperplane arrangements, see \cite{OT}.

\newpage
Let $V$ be a real vector space of finite dimension.
A \textbf{(central) hyperplane arrangement} $\A$ in $V$ is a finite set of linear hyperplanes in $V$. 
As we see in the following example, Hessenberg functions naturally determine hyperplane arrangements.
\begin{example} \label{ex:arrangement} 
Let $V=\{(x_1,\ldots,x_n) \in \R^n \mid x_1+\cdots+x_n =0 \},$ and consider hyperplanes in $V$ given by $H_{i,j}=\{x_j -x_i=0 \}$ for $1 \leq j < i \leq n$. For each Hessenberg function $h:[n]\rightarrow [n]$, the set
$$\A_h:=\{H_{i,j} \mid 1 \leq j < i \leq h(j) \}$$ 
is called an \textbf{ideal arrangement} associated with $h$. 
In particular, if $h=(n,n,\ldots,n)$, then $\A_h$ is called the \textbf{Weyl arrangement} (of type $A_{n-1}$). 
\vspace{-5pt}
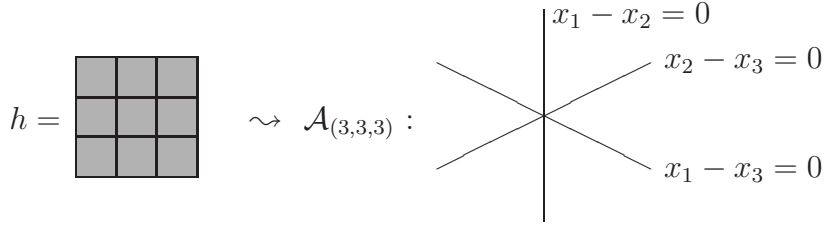
\begin{figure}[h]
\begin{center}
\begin{picture}(250,80)

\put(0,45){\colorbox{gray}}
\put(0,49){\colorbox{gray}}
\put(0,54){\colorbox{gray}}
\put(4,45){\colorbox{gray}}
\put(4,49){\colorbox{gray}}
\put(4,54){\colorbox{gray}}
\put(9,45){\colorbox{gray}}
\put(9,49){\colorbox{gray}}
\put(9,54){\colorbox{gray}}

\put(15,45){\colorbox{gray}}
\put(15,49){\colorbox{gray}}
\put(15,54){\colorbox{gray}}
\put(19,45){\colorbox{gray}}
\put(19,49){\colorbox{gray}}
\put(19,54){\colorbox{gray}}
\put(24,45){\colorbox{gray}}
\put(24,49){\colorbox{gray}}
\put(24,54){\colorbox{gray}}

\put(30,45){\colorbox{gray}}
\put(30,49){\colorbox{gray}}
\put(30,54){\colorbox{gray}}
\put(34,45){\colorbox{gray}}
\put(34,49){\colorbox{gray}}
\put(34,54){\colorbox{gray}}
\put(39,45){\colorbox{gray}}
\put(39,49){\colorbox{gray}}
\put(39,54){\colorbox{gray}}

\put(0,30){\colorbox{gray}}
\put(0,34){\colorbox{gray}}
\put(0,39){\colorbox{gray}}
\put(4,30){\colorbox{gray}}
\put(4,34){\colorbox{gray}}
\put(4,39){\colorbox{gray}}
\put(9,30){\colorbox{gray}}
\put(9,34){\colorbox{gray}}
\put(9,39){\colorbox{gray}}

\put(15,30){\colorbox{gray}}
\put(15,34){\colorbox{gray}}
\put(15,39){\colorbox{gray}}
\put(19,30){\colorbox{gray}}
\put(19,34){\colorbox{gray}}
\put(19,39){\colorbox{gray}}
\put(24,30){\colorbox{gray}}
\put(24,34){\colorbox{gray}}
\put(24,39){\colorbox{gray}}

\put(30,30){\colorbox{gray}}
\put(30,34){\colorbox{gray}}
\put(30,39){\colorbox{gray}}
\put(34,30){\colorbox{gray}}
\put(34,34){\colorbox{gray}}
\put(34,39){\colorbox{gray}}
\put(39,30){\colorbox{gray}}
\put(39,34){\colorbox{gray}}
\put(39,39){\colorbox{gray}}

\put(0,15){\colorbox{gray}}
\put(0,19){\colorbox{gray}}
\put(0,24){\colorbox{gray}}
\put(4,15){\colorbox{gray}}
\put(4,19){\colorbox{gray}}
\put(4,24){\colorbox{gray}}
\put(9,15){\colorbox{gray}}
\put(9,19){\colorbox{gray}}
\put(9,24){\colorbox{gray}}

\put(15,15){\colorbox{gray}}
\put(15,19){\colorbox{gray}}
\put(15,24){\colorbox{gray}}
\put(19,15){\colorbox{gray}}
\put(19,19){\colorbox{gray}}
\put(19,24){\colorbox{gray}}
\put(24,15){\colorbox{gray}}
\put(24,19){\colorbox{gray}}
\put(24,24){\colorbox{gray}}

\put(30,15){\colorbox{gray}}
\put(30,19){\colorbox{gray}}
\put(30,24){\colorbox{gray}}
\put(34,15){\colorbox{gray}}
\put(34,19){\colorbox{gray}}
\put(34,24){\colorbox{gray}}
\put(39,15){\colorbox{gray}}
\put(39,19){\colorbox{gray}}
\put(39,24){\colorbox{gray}}

\put(0,12){\framebox(15,15)}
\put(15,12){\framebox(15,15)}
\put(30,12){\framebox(15,15)}
\put(0,27){\framebox(15,15)}
\put(15,27){\framebox(15,15)}
\put(30,27){\framebox(15,15)}
\put(0,42){\framebox(15,15)}
\put(15,42){\framebox(15,15)}
\put(30,42){\framebox(15,15)}

\put(175,35){\line(0,1){40}}
\put(175,35){\line(0,-1){40}}
\put(175,35){\line(2,1){40}}
\put(175,35){\line(-2,-1){40}}
\put(175,35){\line(2,-1){40}}
\put(175,35){\line(-2,1){40}}

\put(178,70){$x_1-x_2=0$}
\put(220,53){$x_2-x_3=0$}
\put(220,12){$x_1-x_3=0$}

\put(65,30){$\leadsto$}
\put(85,30){$\mathcal{A}_{(3,3,3)}$ :}

\put(-25,30){$h=$}
\end{picture}
\end{center}
\caption{The Weyl arrangement for $n=3$.}
\label{picture:the Weyl arrangement for $n=3$}
\end{figure}
\vspace{-10pt}
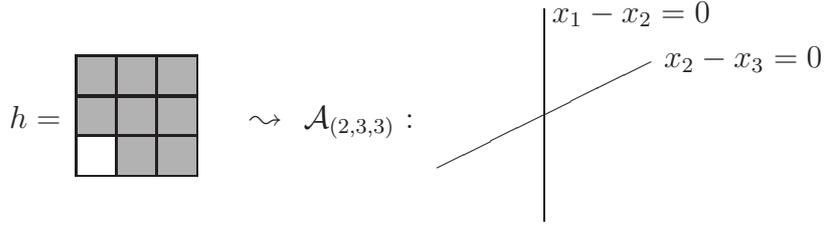
\begin{figure}[h]
\begin{center}
\begin{picture}(250,80)
\put(0,45){\colorbox{gray}}
\put(0,49){\colorbox{gray}}
\put(0,54){\colorbox{gray}}
\put(4,45){\colorbox{gray}}
\put(4,49){\colorbox{gray}}
\put(4,54){\colorbox{gray}}
\put(9,45){\colorbox{gray}}
\put(9,49){\colorbox{gray}}
\put(9,54){\colorbox{gray}}

\put(15,45){\colorbox{gray}}
\put(15,49){\colorbox{gray}}
\put(15,54){\colorbox{gray}}
\put(19,45){\colorbox{gray}}
\put(19,49){\colorbox{gray}}
\put(19,54){\colorbox{gray}}
\put(24,45){\colorbox{gray}}
\put(24,49){\colorbox{gray}}
\put(24,54){\colorbox{gray}}

\put(30,45){\colorbox{gray}}
\put(30,49){\colorbox{gray}}
\put(30,54){\colorbox{gray}}
\put(34,45){\colorbox{gray}}
\put(34,49){\colorbox{gray}}
\put(34,54){\colorbox{gray}}
\put(39,45){\colorbox{gray}}
\put(39,49){\colorbox{gray}}
\put(39,54){\colorbox{gray}}

\put(0,30){\colorbox{gray}}
\put(0,34){\colorbox{gray}}
\put(0,39){\colorbox{gray}}
\put(4,30){\colorbox{gray}}
\put(4,34){\colorbox{gray}}
\put(4,39){\colorbox{gray}}
\put(9,30){\colorbox{gray}}
\put(9,34){\colorbox{gray}}
\put(9,39){\colorbox{gray}}

\put(15,30){\colorbox{gray}}
\put(15,34){\colorbox{gray}}
\put(15,39){\colorbox{gray}}
\put(19,30){\colorbox{gray}}
\put(19,34){\colorbox{gray}}
\put(19,39){\colorbox{gray}}
\put(24,30){\colorbox{gray}}
\put(24,34){\colorbox{gray}}
\put(24,39){\colorbox{gray}}

\put(30,30){\colorbox{gray}}
\put(30,34){\colorbox{gray}}
\put(30,39){\colorbox{gray}}
\put(34,30){\colorbox{gray}}
\put(34,34){\colorbox{gray}}
\put(34,39){\colorbox{gray}}
\put(39,30){\colorbox{gray}}
\put(39,34){\colorbox{gray}}
\put(39,39){\colorbox{gray}}

\put(15,15){\colorbox{gray}}
\put(15,19){\colorbox{gray}}
\put(15,24){\colorbox{gray}}
\put(19,15){\colorbox{gray}}
\put(19,19){\colorbox{gray}}
\put(19,24){\colorbox{gray}}
\put(24,15){\colorbox{gray}}
\put(24,19){\colorbox{gray}}
\put(24,24){\colorbox{gray}}

\put(30,15){\colorbox{gray}}
\put(30,19){\colorbox{gray}}
\put(30,24){\colorbox{gray}}
\put(34,15){\colorbox{gray}}
\put(34,19){\colorbox{gray}}
\put(34,24){\colorbox{gray}}
\put(39,15){\colorbox{gray}}
\put(39,19){\colorbox{gray}}
\put(39,24){\colorbox{gray}}

\put(0,12){\framebox(15,15)}
\put(15,12){\framebox(15,15)}
\put(30,12){\framebox(15,15)}
\put(0,27){\framebox(15,15)}
\put(15,27){\framebox(15,15)}
\put(30,27){\framebox(15,15)}
\put(0,42){\framebox(15,15)}
\put(15,42){\framebox(15,15)}
\put(30,42){\framebox(15,15)}

\put(175,35){\line(0,1){40}}
\put(175,35){\line(0,-1){40}}
\put(175,35){\line(2,1){40}}
\put(175,35){\line(-2,-1){40}}

\put(178,70){$x_1-x_2=0$}
\put(220,53){$x_2-x_3=0$}

\put(65,30){$\leadsto$}
\put(85,30){$\mathcal{A}_{(2,3,3)}$ :}

\put(-25,30){$h=$}

\put(0,45){\colorbox{gray}}
\put(0,49){\colorbox{gray}}
\put(0,54){\colorbox{gray}}
\put(4,45){\colorbox{gray}}
\put(4,49){\colorbox{gray}}
\put(4,54){\colorbox{gray}}
\put(9,45){\colorbox{gray}}
\put(9,49){\colorbox{gray}}
\put(9,54){\colorbox{gray}}

\put(15,45){\colorbox{gray}}
\put(15,49){\colorbox{gray}}
\put(15,54){\colorbox{gray}}
\put(19,45){\colorbox{gray}}
\put(19,49){\colorbox{gray}}
\put(19,54){\colorbox{gray}}
\put(24,45){\colorbox{gray}}
\put(24,49){\colorbox{gray}}
\put(24,54){\colorbox{gray}}

\put(30,45){\colorbox{gray}}
\put(30,49){\colorbox{gray}}
\put(30,54){\colorbox{gray}}
\put(34,45){\colorbox{gray}}
\put(34,49){\colorbox{gray}}
\put(34,54){\colorbox{gray}}
\put(39,45){\colorbox{gray}}
\put(39,49){\colorbox{gray}}
\put(39,54){\colorbox{gray}}

\put(15,30){\colorbox{gray}}
\put(15,34){\colorbox{gray}}
\put(15,39){\colorbox{gray}}
\put(19,30){\colorbox{gray}}
\put(19,34){\colorbox{gray}}
\put(19,39){\colorbox{gray}}
\put(24,30){\colorbox{gray}}
\put(24,34){\colorbox{gray}}
\put(24,39){\colorbox{gray}}

\put(30,30){\colorbox{gray}}
\put(30,34){\colorbox{gray}}
\put(30,39){\colorbox{gray}}
\put(34,30){\colorbox{gray}}
\put(34,34){\colorbox{gray}}
\put(34,39){\colorbox{gray}}
\put(39,30){\colorbox{gray}}
\put(39,34){\colorbox{gray}}
\put(39,39){\colorbox{gray}}

\put(30,15){\colorbox{gray}}
\put(30,19){\colorbox{gray}}
\put(30,24){\colorbox{gray}}
\put(34,15){\colorbox{gray}}
\put(34,19){\colorbox{gray}}
\put(34,24){\colorbox{gray}}
\put(39,15){\colorbox{gray}}
\put(39,19){\colorbox{gray}}
\put(39,24){\colorbox{gray}}

\put(0,12){\framebox(15,15)}
\put(15,12){\framebox(15,15)}
\put(30,12){\framebox(15,15)}
\put(0,27){\framebox(15,15)}
\put(15,27){\framebox(15,15)}
\put(30,27){\framebox(15,15)}
\put(0,42){\framebox(15,15)}
\put(15,42){\framebox(15,15)}
\put(30,42){\framebox(15,15)}

\end{picture}
\end{center}
\caption{The ideal arrangement associated with $h=(2,3,3)$.}
\label{picture:the ideal arrangement associated with $h=(2,3,3)$}
\end{figure} 
\end{example}
\vspace{-10pt}

Let $\CR=\mbox{Sym}(V^*)$ be the symmetric algebra of $V^*$, where $V^*$ is the dual space of $V$. We regard $\CR$ as an algebra of polynomial functions on $V$.
A map $\theta: \CR \to \CR$ is an \textbf{$\R$-derivation} if it satisfies 
\begin{enumerate}  
\item $\theta$ is $\R$-linear, 
\item $\theta(f \cdot g)=\theta(f) \cdot g+f \cdot \theta(g)$ for all $f,g \in \CR$.
\end{enumerate}
We denote the module of $\R$-derivations $\theta: \CR \to \CR$ by $\Der\CR$. If one chooses a linear coordinate system $x_1,\dots, x_m$ on $V$, i.e., $x_1,\ldots,x_m$ is a basis for $V^*$, then the $\CR$-module $\Der \CR$ can be expressed as $\bigoplus_{i=1}^m \CR \, \frac{\partial}{\partial{x_i}}$ where $\frac{\partial}{\partial{x_i}}$ denotes the partial derivative with respect to $x_i$.

Let $\A$ be a hyperplane arrangement in $V$.
For each $H \in \A$, let $\alpha_H \in V^*$ be the defining linear form of $H$ so that $H=\ker(\alpha_H)$. The \textbf{logarithmic derivation module} 
$D(\A)$ of $\A$ is an $\CR$-module defined by
\[
D(\A):=\{\theta \in \Der \CR \mid 
\theta(\alpha_H) \in \CR \alpha_H \ \text{ for all }H \in \A\}.
\]
Geometrically, this consists of polynomial vector fields on $V$ tangent to $\A$.

\begin{example} \label{ex:logarithmic}
Let $V=\R^2$ and $H=\{x=0 \}$ a hyperplane. We consider a hyperplane arrangement $\A=\{H \}$.
An element $\theta \in \Der\CR$ can be written as 
$$
\theta=f \, \frac{\partial}{\partial x} + g \, \frac{\partial}{\partial y}
$$ 
for some $f, g \in \CR=\R[x,y]$. 
Since we can take $\alpha_H=x$, we have $\theta(\alpha_H)=\theta(x)=f$.
Hence, $\theta$ belongs to the logarithmic derivation module $D(\A)$ if and only if $f$ is divisible by $x$.
From this, 
it is straightforward to see that 
$D(\A)$ is free module over $\R[x,y]$ with a basis $x \, \frac{\partial}{\partial x}$ and $\frac{\partial}{\partial y}$. Hence, $D(\A)$ consists of polynomial vector fields on $\R^2$ tangent to $\A$, as desired.\vspace{-7pt}
\begin{figure}[h]
\begin{center}
\begin{picture}(200,110)
  \put(0,70){\line(0,1){40}}
  \put(0,70){\line(0,-1){40}}
  \put(0,40){\circle{1}}
  \put(0,55){\circle{1}} 
  \put(0,70){\circle{1}}
  \put(0,85){\circle{1}}
  \put(0,100){\circle{1}}   
  \put(10,40){\vector(1,0){10}}
  \put(10,55){\vector(1,0){10}} 
  \put(10,70){\vector(1,0){10}}
  \put(10,85){\vector(1,0){10}}
  \put(10,100){\vector(1,0){10}}   
  \put(25,40){\vector(1,0){25}}
  \put(25,55){\vector(1,0){25}} 
  \put(25,70){\vector(1,0){25}}
  \put(25,85){\vector(1,0){25}}
  \put(25,100){\vector(1,0){25}}          
  \put(-10,40){\vector(-1,0){10}}
  \put(-10,55){\vector(-1,0){10}} 
  \put(-10,70){\vector(-1,0){10}}
  \put(-10,85){\vector(-1,0){10}}
  \put(-10,100){\vector(-1,0){10}}   
  \put(-25,40){\vector(-1,0){25}}
  \put(-25,55){\vector(-1,0){25}} 
  \put(-25,70){\vector(-1,0){25}}
  \put(-25,85){\vector(-1,0){25}}
  \put(-25,100){\vector(-1,0){25}}
  \put(-3,15){$x=0$}
  \put(-50,0){the vector field $x \, \frac{\partial}{\partial x}$}

  \put(150,70){\line(0,1){40}}
  \put(150,70){\line(0,-1){40}}
  \put(150,35){\vector(0,1){10}} 
  \put(150,55){\vector(0,1){10}}
  \put(150,75){\vector(0,1){10}}
  \put(150,95){\vector(0,1){10}}   
  \put(165,35){\vector(0,1){10}} 
  \put(165,55){\vector(0,1){10}}
  \put(165,75){\vector(0,1){10}}
  \put(165,95){\vector(0,1){10}}   
  \put(180,35){\vector(0,1){10}}
  \put(180,55){\vector(0,1){10}}
  \put(180,75){\vector(0,1){10}}
  \put(180,95){\vector(0,1){10}}          
  \put(135,35){\vector(0,1){10}}
  \put(135,55){\vector(0,1){10}}
  \put(135,75){\vector(0,1){10}}
  \put(135,95){\vector(0,1){10}}   
  \put(120,35){\vector(0,1){10}}
  \put(120,55){\vector(0,1){10}}
  \put(120,75){\vector(0,1){10}}
  \put(120,95){\vector(0,1){10}}
  \put(147,15){$x=0$}
  \put(105,0){the vector field $\frac{\partial}{\partial y}$}
\end{picture}
\end{center}
\caption{A basis of the logarithmic derivation module of the hyperplane arrangement $\{x=0 \}$ in $\R^2$.}
\label{picture:a basis of the logarithmic derivation module of the hyperplane arrangement}
\end{figure}
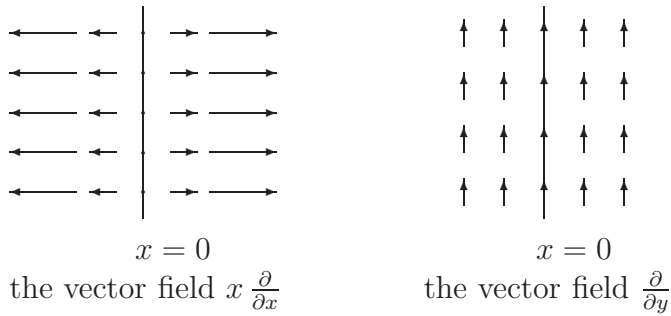
\end{example}

\vspace{-15pt}
In Example~\ref{ex:logarithmic}, we saw that the logarithmic derivation module $D(\A)$ is a free $\CR$-module. 
However, the logarithmic derivation module $D(\A)$ of a hyperplane arrangement $\A$ is not free in general (cf. \cite[Example~$4.34$]{OT}). A hyperplane arrangement $\A$ is a \textbf{free arrangement} if its logarithmic derivation module $D(\A)$ is a free module over $\CR$. Note that if $\A$ is a free arrangement in $V$, then $D(\A)$ has a basis consisting of $m$ homogeneous elements where $m:=\dim V$ (cf. \cite[Proposition~4.18]{OT}). Here, a nonzero element $\theta\in \Der \CR=\CR\otimes V$ is homogeneous if $\theta=\sum_{k=1}^{\ell}f_k\otimes v_k$ $(f_k\in \CR,\ v_k\in V)$ and all non-zero $f_k$'s are homogeneous polynomials of the same degree. 

We now explain the connection with the cohomology rings of regular nilpotent Hessenberg varieties. 
Let $h$ be a Hessenberg function and $\A_h$ the ideal arrangement in $V=\{(x_1,\ldots,x_n) \in \R^n \mid x_1+\cdots+x_n =0 \}$ given in Example \ref{ex:arrangement}.
First, it is known to be free (for arbitrary Lie type) by 
T. Abe, M. Barakat, M. Cuntz, T. Hoge, and H. Terao.

\begin{theorem}[\cite{ABCHT}] \label{theorem:idealfree}
Any ideal arrangement $\A_h$ is free.
\end{theorem}

We next define an ideal of $\CR$ from the logarithmic derivation module $D(\A_h)$ as follows. Let $Q$ be an $S_n$-invariant non-degenerate quadratic form on $V$, which is unique up to a non-zero scalar multiple.
We may take $Q=x_1^2+\cdots+x_n^2 \in \mbox{Sym}^2(V^*)^{S_n}$. We define an ideal $\mathfrak{a}(h)$ of $\CR$ by
$$
\mathfrak{a}(h):=\{\theta(Q) \in \CR \mid \theta \in D(\A_h) \}.
$$
T. Abe, M. Masuda, S. Murai, T. Sato, and the second author proved the following.

\begin{theorem}[\cite{AHMMS}] \label{theorem:AHMMS}
There  is a ring isomorphism
$$
H^*(\Hess(N,h);\R) \cong \CR/\mathfrak{a}(h)
$$
which sends $x_i$ to $\bar x_{i}=-c_1(E_i/E_{i-1})|_{\Hess(N,h)}$ of \eqref {eq:ChernTautological}.
\end{theorem}

In arbitrary Lie type, 
Theorem~\ref{theorem:AHMMS} is proved in \cite{AHMMS}. It is known that the Poincar\'e polynomial of $\Hess(N,h)$ has a summation formula such as \eqref{eq:PoinHess(N,h)Sum} (\cite{precup13a}).
On the other hand, the Hilbert series of the quotient ring $\CR/\mathfrak{a}(h)$ has a product formula such as \eqref{eq:PoinHess(N,h)Prod}.
Theorem~\ref{theorem:AHMMS} gives an affirmative answer to a conjecture of Sommers and Tymoczko in \cite{so-ty} which states that these formulas are equal. 

\begin{remark}
The quotient ring in the right-hand side of Theorem~\ref{theorem:AHMMS} is an example of  Solomon-Terao algebras studied in \cite{ES} and \cite{AMMN}, where \cite{ES} considered more general hypersurface singularities.
The Solomon-Terao algebra $ST(\A,\eta)$ is defined by a hyperplane arrangement $\A$ and a homogeneous polynomial $\eta$. 
Motivated by the work of L. Solomon and H. Terao in \cite{SoloT}, it is proved in \cite{AMMN} that $ST(\A,\eta)$ for a generic $\eta$ is a complete intersection if and only if the hyperplane arrangement $\A$ is free. 
In \cite{ES} the same equivalence is proved more generally for hypersurface singularities which are holonomic in the sense of K. Saito \cite[(3.8)]{S2}. 
This generalization was obtained independently and published slightly earlier. 
\end{remark}

Theorem~\ref{theorem:AHMMS} also tells us that if we can find an explicit $\CR$-basis of $D(\A)$, then we obtain an explicit presentation of the cohomology ring of $\Hess(N,h)$. 
In fact, we can recover the presentation \eqref{eq:AHHM} as follows. First, we recall a well-known criterion for bases of logarithmic derivation modules.

\begin{theorem}[Saito's criterion, \cite{S2}, see also \cite{OT}] \label{theorem:Saito's criterion}
Let $\A$ be a hyperplane arrangement in $V$ and let $\theta_1,\ldots,\theta_m \in D(\A)$ be homogeneous derivations. Then $\theta_1,\ldots,
\theta_m$ form an $\CR$-basis of $D(\A)$ if and only if 
$\theta_1,\ldots,\theta_m$ are $\CR$-independent and $\sum_{i=1}^m \deg \theta_i=|\A|$.
\label{Saito}
\end{theorem}
 
To recover the presentation \eqref{eq:AHHM} by Saito's criterion, 
we define derivations $\psi_{i,j}$ for $1 \leq j \leq i \leq n$ on the ambient vector space $V=\{(x_1,\ldots,x_n) \in \R^n \mid x_1+\cdots+x_n =0 \}$ of ideal arrangements as follows: 
\vspace{-3pt}
$$
\psi_{i,j} := \sum_{k=1}^j \left( \prod_{\ell=j+1}^{i} (x_k-x_\ell)\right)\left(\frac{\partial}{\partial x_k}-\frac 1 n \overline \partial \right) \in \Der \CR\vspace{-3pt}
$$
where $\overline \partial:=\frac{\partial}{\partial x_1}+\cdots+\frac{\partial}{\partial x_n}$ and we take by convention $\prod_{\ell=j+1}^i (x_k-x_{\ell})=1$ whenever $i=j$.
Note that $\frac{\partial}{\partial x_k}$ is not an element of $\Der \CR$ but $(\frac{\partial}{\partial x_k}-\frac 1 n \overline \partial )$ is, since $\CR=\R[x_1,\ldots,x_n]/(x_1+\cdots+x_n)$. Using Theorem~\ref{theorem:Saito's criterion}, one can verify that $\{\psi_{h(j),j}\mid 1 \leq j \leq n-1 \}$ form an $\CR$-basis of $D(\A_h)$ (\cite[Proposition~10.3]{AHMMS}).
Since $\psi_{i,j}(Q)=\psi_{i,j}(x_1^2+\cdots+x_n^2)=2f_{i,j}$ in $\CR$, we recover the presentation \eqref{eq:AHHM} from Theorem~\ref{theorem:AHMMS}.

Abe, Barakat, Cuntz, Hoge, and Terao in \cite{ABCHT} gave a theoretical method to construct a basis of $D(\A_h)$ over $\CR$ for arbitrary Lie type by classification-free proof. For a construction of an explicit basis  for each Lie type, Barakat, Cuntz, and Hoge provided ones for types $E$ and $F$ by computer when the work of \cite{ABCHT} was in progress. 
Also, Terao and Abe worked for types $A$ and $B$, respectively. 
In \cite{AHMMS}, an explicit basis was constructed for types $A, B, C, G$. 
Motivated by this, Enokizono, Nagaoka, Tsuchiya, and the second author in \cite{EHNT1} introduced and studied uniform bases for the logarithmic derivation modules of the ideal arrangements.
In particular, from Theorem~\ref{theorem:AHMMS} which is valid for all Lie types, we obtained explicit presentations of the cohomology rings of regular nilpotent Hessenberg
varieties in all Lie types.

\smallskip
\vspace{10pt}

\subsection{Chromatic symmetric functions in graph theory}
In Section \ref{sect: regular semisimple}, we explained the representation of $S_n$ on the cohomology ring $H^*(\Hess(S,h);\C)$ of regular semisimple Hessenberg varieties constructed by J. Tymoczko. 
Recall that there is a natural correspondence between representations of $S_n$ and symmetric functions of degree $n$ (cf.\ \cite[Section 7]{fult97}). In this section, we describe the symmetric function corresponding to the $\Sn$-representation on $H^*(\Hess(S,h);\C)$ which was conjectured by J. Shareshian and M. Wachs (\cite[Conjecture 1.2]{sh-wa11}, \cite[Conjecture 1.4]{sh-wa14}) in terms of chromatic symmetric functions of a graph determined by $h$. This conjecture is proved by P. Brosnan and T. Chow in \cite{br-ch} and soon after by M. Guay-Paquet (\cite{gu}). 
We refer \cite{fult97} for notations and elementary knowledge of symmetric functions.

Let $h:[n]\rightarrow[n]$ be a Hessenberg function and $G_h$  a graph on the vertex set $[n]$ defined as follows; there is an edge between the vertices $i, j\in[n]$ if $j<i\leq h(j)$. 
Namely, if we visualize $h$ as a configuration of boxes as in Example \ref{example:HessenbergFunction}, we have an edge between $i$ and $j$ if and only if we have a box in $(i,j)$-th position which is strictly below the diagonal (cf.\ ~Figure \ref{pic:dimHess(N,h)}).

\vspace{5pt}

\begin{example}\label{ex: SW graph}
If we take $h=(n,n,\ldots,n)$, then the graph $G_h$ is the complete graph\footnote{A complete graph is a graph in which every pair of distinct vertices is connected by an edge.} on the vertex set $[n]$.
If we take $h=(2,3,4,\ldots,n,n)$, then the graph $G_h$ is the path graph on the vertex set $[n]$.
\vspace{10pt}
\begin{figure}[h]
\begin{center}
\begin{picture}(250,60)

\put(0,45){\colorbox{gray}}
\put(0,49){\colorbox{gray}}
\put(0,54){\colorbox{gray}}
\put(4,45){\colorbox{gray}}
\put(4,49){\colorbox{gray}}
\put(4,54){\colorbox{gray}}
\put(9,45){\colorbox{gray}}
\put(9,49){\colorbox{gray}}
\put(9,54){\colorbox{gray}}

\put(15,45){\colorbox{gray}}
\put(15,49){\colorbox{gray}}
\put(15,54){\colorbox{gray}}
\put(19,45){\colorbox{gray}}
\put(19,49){\colorbox{gray}}
\put(19,54){\colorbox{gray}}
\put(24,45){\colorbox{gray}}
\put(24,49){\colorbox{gray}}
\put(24,54){\colorbox{gray}}

\put(30,45){\colorbox{gray}}
\put(30,49){\colorbox{gray}}
\put(30,54){\colorbox{gray}}
\put(34,45){\colorbox{gray}}
\put(34,49){\colorbox{gray}}
\put(34,54){\colorbox{gray}}
\put(39,45){\colorbox{gray}}
\put(39,49){\colorbox{gray}}
\put(39,54){\colorbox{gray}}

\put(0,30){\colorbox{gray}}
\put(0,34){\colorbox{gray}}
\put(0,39){\colorbox{gray}}
\put(4,30){\colorbox{gray}}
\put(4,34){\colorbox{gray}}
\put(4,39){\colorbox{gray}}
\put(9,30){\colorbox{gray}}
\put(9,34){\colorbox{gray}}
\put(9,39){\colorbox{gray}}

\put(15,30){\colorbox{gray}}
\put(15,34){\colorbox{gray}}
\put(15,39){\colorbox{gray}}
\put(19,30){\colorbox{gray}}
\put(19,34){\colorbox{gray}}
\put(19,39){\colorbox{gray}}
\put(24,30){\colorbox{gray}}
\put(24,34){\colorbox{gray}}
\put(24,39){\colorbox{gray}}

\put(30,30){\colorbox{gray}}
\put(30,34){\colorbox{gray}}
\put(30,39){\colorbox{gray}}
\put(34,30){\colorbox{gray}}
\put(34,34){\colorbox{gray}}
\put(34,39){\colorbox{gray}}
\put(39,30){\colorbox{gray}}
\put(39,34){\colorbox{gray}}
\put(39,39){\colorbox{gray}}

\put(0,15){\colorbox{gray}}
\put(0,19){\colorbox{gray}}
\put(0,24){\colorbox{gray}}
\put(4,15){\colorbox{gray}}
\put(4,19){\colorbox{gray}}
\put(4,24){\colorbox{gray}}
\put(9,15){\colorbox{gray}}
\put(9,19){\colorbox{gray}}
\put(9,24){\colorbox{gray}}

\put(15,15){\colorbox{gray}}
\put(15,19){\colorbox{gray}}
\put(15,24){\colorbox{gray}}
\put(19,15){\colorbox{gray}}
\put(19,19){\colorbox{gray}}
\put(19,24){\colorbox{gray}}
\put(24,15){\colorbox{gray}}
\put(24,19){\colorbox{gray}}
\put(24,24){\colorbox{gray}}

\put(30,15){\colorbox{gray}}
\put(30,19){\colorbox{gray}}
\put(30,24){\colorbox{gray}}
\put(34,15){\colorbox{gray}}
\put(34,19){\colorbox{gray}}
\put(34,24){\colorbox{gray}}
\put(39,15){\colorbox{gray}}
\put(39,19){\colorbox{gray}}
\put(39,24){\colorbox{gray}}

\put(0,12){\framebox(15,15)}
\put(15,12){\framebox(15,15)}
\put(30,12){\framebox(15,15)}
\put(0,27){\framebox(15,15)}
\put(15,27){\framebox(15,15)}
\put(30,27){\framebox(15,15)}
\put(0,42){\framebox(15,15)}
\put(15,42){\framebox(15,15)}
\put(30,42){\framebox(15,15)}

\put(145,33){\circle*{5}}
\put(175,33){\circle*{5}}
\put(205,33){\circle*{5}}
\put(145,33){\line(1,0){30}}
\put(175,33){\line(1,0){30}}
\put(175,33){\oval(60,40)[t]}

\put(142,18){$1$}
\put(172,18){$2$}
\put(202,18){$3$}

\put(65,30){$\leadsto$}
\put(85,30){$G_{(3,3,3)}$ :}

\put(-25,30){$h=$}
\end{picture}
\end{center}
\vspace{-10pt}
\caption{The complete graph on the vertex set $[3]$.}
\label{pic:complete graph}
\end{figure}
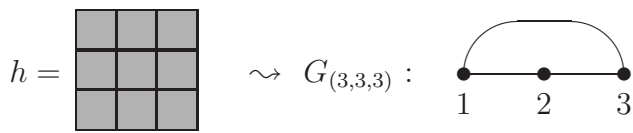
\begin{figure}[h]
\begin{center}
\begin{picture}(250,60)
\put(0,45){\colorbox{gray}}
\put(0,49){\colorbox{gray}}
\put(0,54){\colorbox{gray}}
\put(4,45){\colorbox{gray}}
\put(4,49){\colorbox{gray}}
\put(4,54){\colorbox{gray}}
\put(9,45){\colorbox{gray}}
\put(9,49){\colorbox{gray}}
\put(9,54){\colorbox{gray}}

\put(15,45){\colorbox{gray}}
\put(15,49){\colorbox{gray}}
\put(15,54){\colorbox{gray}}
\put(19,45){\colorbox{gray}}
\put(19,49){\colorbox{gray}}
\put(19,54){\colorbox{gray}}
\put(24,45){\colorbox{gray}}
\put(24,49){\colorbox{gray}}
\put(24,54){\colorbox{gray}}

\put(30,45){\colorbox{gray}}
\put(30,49){\colorbox{gray}}
\put(30,54){\colorbox{gray}}
\put(34,45){\colorbox{gray}}
\put(34,49){\colorbox{gray}}
\put(34,54){\colorbox{gray}}
\put(39,45){\colorbox{gray}}
\put(39,49){\colorbox{gray}}
\put(39,54){\colorbox{gray}}

\put(0,30){\colorbox{gray}}
\put(0,34){\colorbox{gray}}
\put(0,39){\colorbox{gray}}
\put(4,30){\colorbox{gray}}
\put(4,34){\colorbox{gray}}
\put(4,39){\colorbox{gray}}
\put(9,30){\colorbox{gray}}
\put(9,34){\colorbox{gray}}
\put(9,39){\colorbox{gray}}

\put(15,30){\colorbox{gray}}
\put(15,34){\colorbox{gray}}
\put(15,39){\colorbox{gray}}
\put(19,30){\colorbox{gray}}
\put(19,34){\colorbox{gray}}
\put(19,39){\colorbox{gray}}
\put(24,30){\colorbox{gray}}
\put(24,34){\colorbox{gray}}
\put(24,39){\colorbox{gray}}

\put(30,30){\colorbox{gray}}
\put(30,34){\colorbox{gray}}
\put(30,39){\colorbox{gray}}
\put(34,30){\colorbox{gray}}
\put(34,34){\colorbox{gray}}
\put(34,39){\colorbox{gray}}
\put(39,30){\colorbox{gray}}
\put(39,34){\colorbox{gray}}
\put(39,39){\colorbox{gray}}

\put(15,15){\colorbox{gray}}
\put(15,19){\colorbox{gray}}
\put(15,24){\colorbox{gray}}
\put(19,15){\colorbox{gray}}
\put(19,19){\colorbox{gray}}
\put(19,24){\colorbox{gray}}
\put(24,15){\colorbox{gray}}
\put(24,19){\colorbox{gray}}
\put(24,24){\colorbox{gray}}

\put(30,15){\colorbox{gray}}
\put(30,19){\colorbox{gray}}
\put(30,24){\colorbox{gray}}
\put(34,15){\colorbox{gray}}
\put(34,19){\colorbox{gray}}
\put(34,24){\colorbox{gray}}
\put(39,15){\colorbox{gray}}
\put(39,19){\colorbox{gray}}
\put(39,24){\colorbox{gray}}

\put(0,12){\framebox(15,15)}
\put(15,12){\framebox(15,15)}
\put(30,12){\framebox(15,15)}
\put(0,27){\framebox(15,15)}
\put(15,27){\framebox(15,15)}
\put(30,27){\framebox(15,15)}
\put(0,42){\framebox(15,15)}
\put(15,42){\framebox(15,15)}
\put(30,42){\framebox(15,15)}

\put(145,33){\circle*{5}}
\put(175,33){\circle*{5}}
\put(205,33){\circle*{5}}
\put(145,33){\line(1,0){30}}
\put(175,33){\line(1,0){30}}

\put(142,18){$1$}
\put(172,18){$2$}
\put(202,18){$3$}

\put(65,30){$\leadsto$}
\put(85,30){$G_{(2,3,3)}$ :}

\put(-25,30){$h=$}

\put(0,45){\colorbox{gray}}
\put(0,49){\colorbox{gray}}
\put(0,54){\colorbox{gray}}
\put(4,45){\colorbox{gray}}
\put(4,49){\colorbox{gray}}
\put(4,54){\colorbox{gray}}
\put(9,45){\colorbox{gray}}
\put(9,49){\colorbox{gray}}
\put(9,54){\colorbox{gray}}

\put(15,45){\colorbox{gray}}
\put(15,49){\colorbox{gray}}
\put(15,54){\colorbox{gray}}
\put(19,45){\colorbox{gray}}
\put(19,49){\colorbox{gray}}
\put(19,54){\colorbox{gray}}
\put(24,45){\colorbox{gray}}
\put(24,49){\colorbox{gray}}
\put(24,54){\colorbox{gray}}

\put(30,45){\colorbox{gray}}
\put(30,49){\colorbox{gray}}
\put(30,54){\colorbox{gray}}
\put(34,45){\colorbox{gray}}
\put(34,49){\colorbox{gray}}
\put(34,54){\colorbox{gray}}
\put(39,45){\colorbox{gray}}
\put(39,49){\colorbox{gray}}
\put(39,54){\colorbox{gray}}

\put(15,30){\colorbox{gray}}
\put(15,34){\colorbox{gray}}
\put(15,39){\colorbox{gray}}
\put(19,30){\colorbox{gray}}
\put(19,34){\colorbox{gray}}
\put(19,39){\colorbox{gray}}
\put(24,30){\colorbox{gray}}
\put(24,34){\colorbox{gray}}
\put(24,39){\colorbox{gray}}

\put(30,30){\colorbox{gray}}
\put(30,34){\colorbox{gray}}
\put(30,39){\colorbox{gray}}
\put(34,30){\colorbox{gray}}
\put(34,34){\colorbox{gray}}
\put(34,39){\colorbox{gray}}
\put(39,30){\colorbox{gray}}
\put(39,34){\colorbox{gray}}
\put(39,39){\colorbox{gray}}

\put(30,15){\colorbox{gray}}
\put(30,19){\colorbox{gray}}
\put(30,24){\colorbox{gray}}
\put(34,15){\colorbox{gray}}
\put(34,19){\colorbox{gray}}
\put(34,24){\colorbox{gray}}
\put(39,15){\colorbox{gray}}
\put(39,19){\colorbox{gray}}
\put(39,24){\colorbox{gray}}

\put(0,12){\framebox(15,15)}
\put(15,12){\framebox(15,15)}
\put(30,12){\framebox(15,15)}
\put(0,27){\framebox(15,15)}
\put(15,27){\framebox(15,15)}
\put(30,27){\framebox(15,15)}
\put(0,42){\framebox(15,15)}
\put(15,42){\framebox(15,15)}
\put(30,42){\framebox(15,15)}

\end{picture}
\end{center}
\vspace{-10pt}
\caption{The path graph on the vertex set $[3]$.}
\label{pic:path graph}
\end{figure}
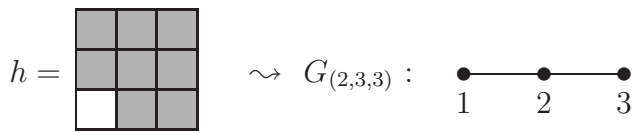 
\end{example}

For a graph $G=(V,E)$ on the vertex set $V=[n]$, Shareshian-Wachs (\cite{sh-wa11, sh-wa14}) introduced 
the chromatic quasisymmetric function $X_{G}(\text{\rm \textbf{x}},\q)$ of $G$ as
\begin{align}\label{eq: chromatic quasisymmetric function}
X_{G}(\textbf{x},\q) = \sum_{\kappa} \big(\prod_{i=1}^n x_{\kappa(i)}\big)\, \q^{\text{asc}(\kappa)}, 
\end{align}
where the summation runs over all proper colorings\footnote{A map $\kappa:[n]\rightarrow \mathbb{N}$ is called a proper coloring of $G$ if $\kappa(i)\neq\kappa(j)$ for all pair of vertices $i$ and $j$ which are connected by an edge.} $\kappa:[n]\rightarrow \mathbb{N}=\{1,2,3,\ldots\}$ of $G$ and $\text{asc}(\kappa)\coloneqq |\{(j,i) \in E\mid j<i, \ \kappa(j)<\kappa(i)\}|$ is the number of ascents of $\kappa$. Here, the variable $\q$ stands for the grading, and this is a graded version of Stanley's chromatic symmetric function $X_{G}(\textbf{x})$ of $G$ (\cite{Stanley1995}). In general, $X_{G}(\textbf{x},\q)$ is quasisymmetric in $\textbf{x}$-variables but may not be symmetric. However, for our graph $G_h$, it is known that $X_{G_h}(\textbf{x},\q)$ is in fact symmetric (\cite[Theorem 4.5]{sh-wa14}). In \cite[Theroem 6.3]{sh-wa14}, the Schur basis expansion of $X_{G_h}(\textbf{x},\q)$ is determined in terms of combinatorics, where the non-graded version was originally obtained by V. Gasharov (\cite{gash}).

\begin{example}
Let $n=3$ and $h=(2,3,3)$. Then, one can verify that 
\begin{align}\label{eq: ch qu-sy for toric}
X_{G_h}(\text{\rm \textbf{x}},\q) = s_{(1,1,1)}(\text{\rm \textbf{x}})+\big(2s_{(1,1,1)}(\text{\rm \textbf{x}}) + s_{(2,1)}(\text{\rm \textbf{x}})\big)\, \q+s_{(1,1,1)}(\text{\rm \textbf{x}})\, \q^2.
\end{align}
Here, $s_{\lambda}(\text{\rm \textbf{x}})$ is the Schur function corresponding to a partition $\lambda$ of $3$. 
\end{example}

The following theorem determines the $\Sn$-representation on $H^*(\Hess(S,h);\C)$ in terms of the combinatorics of the graph $G_h$.
This beautiful theorem was conjectured by Shareshian-Wachs (\cite[Conjecture 1.2]{sh-wa11}, \cite[Conjecture 1.4]{sh-wa14}), and it was proved by Brosnan-Chow (\cite{br-ch}) and soon after by Guay-Paquet (\cite{gu}). We denote by $\text{\rm ch}$ the Frobenius character under which symmetric functions of degree $n$ corresponds to representations of $S_n$.

\begin{theorem} $($\cite{br-ch}, \cite{gu}$)$\label{thm: Shareshian-Wachs conjecture}
Let $\Hess(S,h)$ be a regular semisimple Hessenberg variety.~Then, 
\begin{align*}
\omega X_{G_h}(\text{\rm \textbf{x}},\q) = \sum_{k=0}^{m} \text{\rm ch} H^{2k}(\Hess(S,h);\C) \, \q^k .
\end{align*}
where $m=\dim_{\C} \Hess(S,h)$ and $\omega$ is the involution on the ring of symmetric functions in $\text{\rm \textbf{x}}$-variables sending each Schur function $s_{\lambda}(\text{\rm \textbf{x}})$ to $s_{\tilde{\lambda}}(\text{\rm \textbf{x}})$ associated with the transpose $\tilde{\lambda}$.
\end{theorem}

For the case $h=(n,n,\ldots,n)$, recall from Section \ref{subsec: coh of reg ss} that the $S_n$-representation on $H^*(\Hess(S,h);\C)=H^*(Fl(\C^n);\C)$ is trivial, and it is straightforward to see that $X_{G_h}(\text{\rm \textbf{x}},\q)=e_{n}(\textbf{x})\sum_{\kappa}q^{\text{asc}(\kappa)}$ 
in this case where the latter sum is equal to the Poincar$\acute{\text{e}}$ polynomial of $Fl(\C^n)$. Since the complete symmetric function $h_{n}(\textbf{x})=\omega (e_{n}(\textbf{x}))$ corresponds to the $1$-dimensional trivial representation, the equality in Theorem \ref{thm: Shareshian-Wachs conjecture} holds in this case. For the case $h=(2,3,4,\ldots,n,n)$, both of $X_{G_h}(\text{\rm \textbf{x}},\q)$ and the $\Sn$-representation on $H^*(\Hess(S,h);\C)$ were well-studied in \cite{Stanley1989,proc90,ste}, and the equality was known. See \cite[Section 1]{sh-wa14} for details.
The proof of Theorem \ref{thm: Shareshian-Wachs conjecture} given by Brosnan-Chow is geometric in the sense that they used the theory of monodromy actions, whereas Guay-Paquet provided a combinatorial proof using the theory of Hopf algebras. 

\begin{example}
Let $n=3$ and $h=(2,3,3)$. Then, from \eqref{eq: Sn-rep for toric} and \eqref{eq: ch qu-sy for toric}, we see that
\begin{align*}
\omega X_{G_h}(\text{\rm \textbf{x}},\q) 
&= s_{(3)}(\text{\rm \textbf{x}})+\big(2s_{(3)}(\text{\rm \textbf{x}}) + s_{(2,1)}(\text{\rm \textbf{x}})\big)\, \q+s_{(3)}(\text{\rm \textbf{x}})\, \q^2 \\
&=\text{\rm ch} H^{0}(\Hess(S,h);\C) + \text{\rm ch} H^{2}(\Hess(S,h);\C)\, \q + \text{\rm ch} H^{4}(\Hess(S,h);\C)\, \q^2.
\end{align*}
\end{example}

Theorem \ref{thm: Shareshian-Wachs conjecture} is related with the Stanley-Stembridge conjecture in graph theory as we explain in what follows.
Recall that the chromatic quasisymmetric function $X_{G}(\text{\rm \textbf{x}},\q)$ is a graded version of Stanley's chromatic symmetric function $X_{G}(\text{\rm \textbf{x}})$ (\cite{Stanley1995}) which we obtain by forgetting the grading parameter $\q$ in \eqref{eq: chromatic quasisymmetric function}. 
Given a poset $P$, we can construct its incomparability graph which has as its vertices the elements of $P$ and we have an edge between two vertices if they are not comparable in $P$. The Stanley-Stembridge conjecture is stated as follows. 

\begin{conjecture}\emph{(}\emph{Stanley-Stembridge conjecture} \cite[Conjecture 5.5]{st-ste}, \cite[Conjecture 5.1]{Stanley1995}\emph{)}\label{conjecture:h-posi}
Let $G$ be the incomparability graph of a $(3+1)$-free poset. Then $X_{G}(\text{\rm \textbf{x}})$ is $e$-positive.
\end{conjecture}

Here, a poset is called \textbf{$(r+s)$-free} if the poset does not contain an induced subposet isomorphic to the direct sum of an $r$ element chain and an $s$ element chain, and a symmetric function of degree $n$ is called \textbf{$e$-positive} if it is a non-negative sum of elementary symmetric functions $e_{\lambda}(\text{\rm \textbf{x}})$ for partition $\lambda$ of $n$.  
See \cite{sh-wa14}, for more information on Stanley-Stembridge conjecture.

\begin{example}\label{ex: e-positive}
Take a poset on $[3]=\{1,2,3\}$ for which $1<3$ is the unique comparable pair. It is obviously (3+1)-free.
In this case, the incomparability graph $G$ is the path graph given in Figure \ref{pic:path graph}. Hence, from \eqref{eq: ch qu-sy for toric}, we see that 
\begin{align*}
X_{G}(\text{\rm \textbf{x}}) 
&= s_{(1,1,1)}(\text{\rm \textbf{x}})+\big(2s_{(1,1,1)}(\text{\rm \textbf{x}}) + s_{(2,1)}(\text{\rm \textbf{x}})\big)+s_{(1,1,1)}(\text{\rm \textbf{x}}) \\
&= 3e_{(3)}(\text{\rm \textbf{x}})+e_{(2,1)}(\text{\rm \textbf{x}}),
\end{align*}
as desired since we have $e_{(3)}(\text{\rm \textbf{x}})=s_{(1,1,1)}(\text{\rm \textbf{x}})$ and $e_{(2,1)}(\text{\rm \textbf{x}})=s_{(2,1)}(\text{\rm \textbf{x}})+s_{(1,1,1)}(\text{\rm \textbf{x}})$.
\end{example}
 
Guay-Paquet (\cite{gu13}) showed that, to solve the Stanley-Stembridge conjecture, it is enough to solve it for posets which are $(3+1)$-free and $(2+2)$-free, and this is precisely the case that the corresponding incomparability graph can be identified with $G_h$ for some Hessenberg function $h:[n]\rightarrow[n]$, where $n$ is the number of the elements of the poset, as we encountered in Example \ref{ex: e-positive} (See \cite[Section 4]{sh-wa11} for details). Thus, the Stanley-Stembridge conjecture is now reduced to the following conjecture by Theorem~\ref{thm: Shareshian-Wachs conjecture}.

\begin{conjecture}\emph{(}\cite[Conjecture 5.4]{sh-wa11}, \cite[Conjecture 10.4]{sh-wa14}\emph{)}\label{conjecture:permutation rep}
The $\Sn$-representation on $H^*(\Hess(S,h);\C)$ constructed by Tymoczko is a permutation representation, i.e. a direct sum of  induced representations of the trivial representation from $S_{\lambda}$ to $\Sn$ where $S_{\lambda}=S_{\lambda_1}\times \dots \times S_{\lambda_{\ell}}$ for $\lambda=(\lambda_1,\ldots,\lambda_{\ell})$.
\end{conjecture}

Motivated by the connection to Stanely-Stembridge conjecture, M.~Harada and M.~Precup verified Conjecture \ref{conjecture:h-posi} for so-called \textit{abelian} regular semisimple Hessenberg varieties, and they also derived a set of linear relations satisfied by the multiplicities of certain permutation representations. See \cite{ha-pr17,ha-pr18} for details.

\bigskip

\section{More topics}
\subsection{Semisimple Hessenberg varieties.}
E. Insko and M. Precup studied Hessenberg varieties associated with semisimple matrices which may not be regular semisimple (\cite{InskoPrecup}).
They determined the irreducible components of semisimple Hessenberg varieties 
for $h=(2,3,4,\ldots,n,n)$ in arbitrary Lie type. They also proved that irreducible components are smooth and gave an explicit description of their intersections.

\subsection{Hessenberg varieties for the minimal nilpotent orbit.}
P. Crooks and the first author studied Hessenberg varieties associated with minimal nilpotent matrices (i.e. nilpotent matrices with Jordan blocks of size $1$ and a single Jordan block of size $2$). Explicit descriptions of their Poincar$\acute{\text{e}}$ polynomials and irreducible components are described, and a certain presentation of their cohomology rings are also provided in terms of Schubert classes. See \cite{Abe-Crooks} for details.

\subsection{Regular Hessenberg varieties}
An $n\times n$ matrix $R$ is called \textbf{regular} if the Jordan blocks of $R$ have distinct eigenvalues. For a regular matrix $R$, $\Hess(R,h)$ is called a \textbf{regular Hessenberg variety}. This class of Hessenberg varieties contains regular nilpotent Hessenberg varieties and regular semisimple Hessenberg varieties, and $\Hess(R,h)$ plays an important role in the work of P. Brosnan and T. Chow (\cite{br-ch}) proving Shareshian-Wachs conjecture (Theorem \ref{thm: Shareshian-Wachs conjecture}). 
They are singular varieties in general, however 
M. Precup proved that the Betti numbers of $\Hess(R,h)$ are palindromic (\cite{precup18}).
N. Fujita, H. Zeng, and the first author proved that higher cohomology groups of their structure sheaves vanish and that they degenerate to the regular nilpotent Hessenberg variety $\Hess(N,h)$ if $h(i)\geq i+1$ for all $1\leq i\leq n$ (\cite{ab-fu-ze}), which was motivated by the works of D. Anderson and J. Tymoczo (\cite{an-ty}) and L. DeDieu, F. Galetto, M. Harada, and the first author (\cite{ab-de-ga-ha}).

\subsection{Poincar$\acute{\text{e}}$ dual of Hessenberg varieties}
In \cite{ab-fu-ze}, the Poincar$\acute{\text{e}}$ dual of a regular Hessenberg variety $\Hess(R,h)$ in $H^*(Fl(\C^n))$ was computed in terms of positive roots associated the Hessenberg function $h$, and E. Insko, J. Tymoczko, and A. Woo gave a combinatorial formula for this class  using Schubert polynomials (\cite{InskoTymoczkoWoo}). Also, the cohomology class $[\Hess(R,h)]\in H^*(Fl(\C^n))$ does not depend on a choice of a regular matrix $R$ if $h(i)\geq i+1$ for all $1\leq i\leq n$ (See for details \cite{ab-fu-ze,InskoTymoczkoWoo}).

\subsection{Additive bases of the cohomology rings of regular nilpotent Hessenberg varieties}
M. Enokizono, T. Nagaoka, A. Tsuchiya, and the second author constructed in \cite{EHNT2} an additive basis of the cohomology ring of a regular nilpotent Hessenberg variety $\Hess(N,h)$. 
This basis is obtained by extending the Poincar\'e duals $[\Hess(N,h')]\in H^*(\Hess(N,h))$ of smaller regular nilpotent Hessenberg varieties $\Hess(N,h')$ with $h' \subset h$.
In particular, all of the classes $[\Hess(N,h')]\in H^*(\Hess(N,h))$ with $h' \subset h$, are linearly independent. 
On the other hand, M. Harada, S. Murai, M. Precup, J. Tymoczko, and the second author derive in \cite{HHMPT} a filtration on the cohomology ring $H^*(\Hess(N,h))$ of regular nilpotent Hessenberg varieties, from which they obtain a monomial basis for $H^*(\Hess(N,h))$. 
This basis is different from the one obtained in \cite{EHNT2}. 
From the filtration they additionally obtain an inductive formula for the Poincar\'e polynomials of $\Hess(N,h)$; moreover, the monomial basis has an interpretation in terms of Schubert calculus.

\subsection{The volume polynomials of Hessenberg varieties.}
Recall from Corollary \ref{corollary:Hess(N,h)PDA} that the cohomology ring $H^*(\Hess(N,h);\Q)$ of a regular nilpotent Hessenberg variety is a Poincar$\acute{\text{e}}$ duality algebra. 
In \cite{AHMMS}, T. Abe, M. Masuda, S. Murai, T. Sato, and the second author computed the volume polynomial of this ring $H^*(\Hess(N,h);\Q)$, and it precisely gives  the volume of a certain embedding of any regular Hessenberg variety associated with $h$ into a projective space (\cite{ab-de-ga-ha, ab-fu-ze}). 
M. Harada, M. Masuda, S. Park, and the second author provided a combinatorial formula for this polynomial in terms of the volumes of certain faces of the Gelfand-Zetlin polytope (\cite{ha-ho-ma-pa}).

\subsection{Hessenberg varieties of parabolic type.}
J. Tymoczko and M. Precup showed that the Betti numbers of parabolic Hessenberg varieties decompose into a combination of those of Springer fibers and Schubert varieties associated to the parabolic (\cite{pr-ty17}). As a corollary, they deduced that the Betti numbers of some parabolic Hessenberg varieties in Lie type A are equal to those of a specific union of Schubert varieties.

\subsection{Twins for regular semisimple Hessenberg varieties.}
Given a Hessenberg function $h:[n]\rightarrow[n]$, A. Ayzenberg and V. Buchstaber introduced a smooth $T$-manifold $X_h$, where $T$ is the maximal torus of $\text{GL}(n,\C)$ given in Section \ref{sect: regular semisimple}. This manifold is similar to the regular semisimple Hessenberg variety $\Hess(S,h)$ in some sense. For example, they have the same Betti numbers and their $T$-equivariant cohomology rings are isomorphic \textit{as rings}. See \cite{ay-buHess} for details.

\subsection{Integrable systems and Hessenberg varieties.}
For a Hessenberg function $h:[n]\rightarrow [n]$, we denote by $\mathcal{X}(h)$ the family of Hessenberg varieties associated with $h$. Then, $\mathcal{X}(h)$ in fact have a structure of a vector bundle over the flag variety $Fl(\C^n)$. For the case $h=(2,3,4,\ldots,n,n)$, the family $\mathcal{X}(h)$ contains the Peterson variety and the permutohedral variety. In this special case, P. Crooks and the first author showed that $\mathcal{X}(h)$ admits a Poisson structure with an open dense symplectic leaf, and that there is a completely integrable system on $\mathcal{X}(h)$ which contains the Toda lattice as a sub-system (\cite{Abe-Crooks 2018}). 

\subsection{The poset of Hessenberg varieties.}
E. Drellich studied the poset of the Hessenberg varieties $\Hess(X,h)$ in $Fl(\C^n)$ for a given $n\times n$ matrix $X$. She proved that if $X$ is not a scalar multiple of the identity matrix, then the Hessenberg functions determine distinct Hessenberg varieties.
See \cite{dr2017} for details. 

\subsection{Springer correspondence for symmetric spaces}
As mentioned in Remark~\ref{rem: def of Hess}, Hessenberg varieties can be defined in a more general setting (\cite{G-Ko-Ma06}), and those Hessenberg varieties also appear in the works of T. H. Chen, K. Vilonen, and T. Xue (\cite{ch-vi-xu15-1, ch-vi-xu15-2, ch-vi-xu18}) on Springer correspondence for symmetric spaces.


\begin{thebibliography}{9}

\bibitem{Abe}
   H. Abe, 
   \emph{Young diagrams and intersection numbers for toric manifolds associated with Weyl chambers},
   Electron. J. Combin. \textbf{22(2)} (2015), no. 2, Paper 2.4, 24 pp.

\bibitem{ab-de-ga-ha}
H. Abe, L. DeDieu, F. Galetto, and M. Harada,
\emph{Geometry of Hessenberg varieties with applications to Newton-Okounkov bodies}, 
Selecta Math. (N.S.) \textbf{24} (2018), no. 3, 2129--2163.

\bibitem{Abe-Crooks}
H.~Abe and P.~Crooks, 
\newblock \emph{Hessenberg varieties for the minimal nilpotent orbit}, 
Pure Appl. Math. Q. \textbf{12} (2016), no. 2, 183--223.

\bibitem{Abe-Crooks 2018}
H. Abe and P. Crooks, 
\emph{Hessenberg varieties, Slodowy slices, and integrable systems}, 
Math. Z., \textbf{291}(3) (2019), 1093--1132.

\bibitem{ab-fu-ze}
H. Abe, N. Fujita, and H. Zeng,
\emph{Geometry of regular Hessenberg varieties}, to appear in Transform. Groups.

\bibitem{AHHM}
H. Abe, M. Harada, T. Horiguchi, and M. Masuda,
\emph{The cohomology rings of regular nilpotent Hessenberg varieties in Lie type A}, 
Int. Math. Res. Not., \textbf{2019}(17) (2019), 5316--5388.

\bibitem{AHHM2}
H. Abe, M. Harada, T. Horiguchi, and M. Masuda, 
   \emph{The equivariant cohomology rings of regular nilpotent Hessenberg varieties in Lie type A: Research Announcement}, 
   Morfismos \textbf{18} (2014), No. 2, pp. 51--65.

\bibitem{ab-ho-ma}
H. Abe, T. Horiguchi, and M. Masuda,
\emph{The cohomology rings of regular semisimple Hessenberg varieties for $h=(h(1),n,\ldots,n)$},
J. Comb. \textbf{10} (1) 27--59, 2019.

\bibitem{ABCHT} 
T. Abe, M. Barakat, M. Cuntz, T. Hoge, and H. Terao, 
\emph{The freeness of ideal subarrangements of Weyl arrangements},
J. European Math. Soc. \textbf{18} (2016), 1339--1348.

\bibitem{AHMMS}
T. Abe, T. Horiguchi, M. Masuda, S. Murai, and T. Sato,
\emph{Hessenberg varieties and hyperplane arrangements},
to appear in J. Reine Angew. Math., DOI: 10.1515/crelle-2018-0039.

\bibitem{AMMN}
T. Abe, T. Maeno, S. Murai, and Y. Numata,
\emph{Solomon-Terao algebra of hyperplane arrangements},
J. Math. Soc. Japan 71 (2019), no. 4, 1027--1047.

\bibitem{an-ty}
D. Anderson and J. Tymoczko, 
\emph{Schubert polynomials and classes of Hessenberg varieties}, 
J. Algebra \textbf{323} (2010), no. 10, 2605--2623. 

\bibitem{ay-buHess}
A. Ayzenberg and V. Buchstaber,
\emph{Manifolds of isospectral matrices and Hessenberg varieties},
arXiv:1803.01132.

\bibitem{Balibanu}
	A. Balibanu, 
	\emph{The Peterson Variety and the Wonderful Compactification}, 
	Represent. Theory \textbf{21} (2017), 132--150.

\bibitem{br-ca04} 
M. Brion and J. Carrell, 
\emph{The equivariant cohomology ring of regular varieties}, 
Michigan Math. J. \textbf{52} (2004), no. 1, 189--203. 

\bibitem{br-ch} 
P. Brosnan and T. Chow, \emph{Unit interval orders and the dot action on the cohomology of regular semisimple Hessenberg varieties}, 
Adv. Math. \textbf{329} (2018), 955--1001. 

\bibitem{ch-vi-xu15-1}
T. H. Chen, K. Vilonen, and T. Xue,
\emph{Springer Correspondence for Symmetric Spaces},
arXiv:1510.05986.
	
\bibitem{ch-vi-xu15-2}
T. H. Chen, K. Vilonen, and T. Xue,
\emph{Hessenberg varieties, intersections of quadrics, and the Springer correspondence},
arXiv:1511.00617.

\bibitem{ch-vi-xu18}
T. H. Chen, K. Vilonen, and T. Xue,
\emph{Springer correspondence for the split symmetric pair in type A},
Compos. Math. \textbf{154} (2018), no. 11, 2403--2425. 

\bibitem{De Concini-Lusztig-Procesi}
C. De Concini, G. Lusztig, C. Procesi, 
\emph{Homology of the zero-set of a nilpotent vector field on a flag manifold}, 
J. Amer. Math. Soc. \textbf{1} (1988), no. 1, 15--34. 

\bibitem{ma-pr-sh}
F. De Mari, C. Procesi, and M. A. Shayman,
\emph{Hessenberg varieties}, 
Trans. Amer. Math. Soc. {\bf 332} (1992),  no. 2, 529--534. 

\bibitem{ma-sh}
F. De Mari and M. A. Shayman, 
\emph{Generalized Eulerian numbers and the topology of the Hessenberg variety of a matrix}, 
Acta Appl. Math. \textbf{12} (1988), no. 3, 213--235.

\bibitem{dr2015}
E. Drellich, 
\emph{Monk's rule and Giambelli's formula for Peterson varieties of all Lie types}, 
J. Algebraic Combin. \textbf{41} (2015), no. 2, 539--575. 

\bibitem{Dre1}
E. Drellich,
\emph{Combinatorics of equivariant cohomology: Flags and regular nilpotent Hessenberg varieties}, 
PhD thesis, University of Massachusetts, 2015.

\bibitem{dr2017}
E. Drellich,
\emph{The Containment Poset of Type A Hessenberg Varieties},
arXiv:1710.05412.

\bibitem{EHNT1}
M. Enokizono, T. Horiguchi, T. Nagaoka, and A. Tsuchiya,
\emph{Uniform bases for ideal arrangements},
arXiv:1912.02448.

\bibitem{EHNT2}
M. Enokizono, T. Horiguchi, T. Nagaoka, and A. Tsuchiya,
\emph{An additive basis for the cohomology rings of regular nilpotent Hessenberg varieties},
in preparation.

\bibitem{ES}
R. Epure and M. Schulze,
\emph{A Saito criterion for holonomic divisors}, 
Manuscripta Math. \textbf{160} (2019), no. 1-2, 1--8.

\bibitem{fu-ha-ma}
Y. Fukukawa, M. Harada, and M. Masuda, 
\emph{The equivariant cohomology rings of Peterson varieties}, 
J. Math. Soc. Japan \textbf{67} (2015), no. 3, 1147--1159. 

\bibitem{fult97}
W. Fulton, \emph{Young Tableaux}, 
London Mathematical Society Student Texts, 35. Cambridge University Press, Cambridge.

\bibitem{gash}
V. Gasharov, 
\emph{Incomparability graphs of (3+1)-free posets are s-positive},
Proceedings of the 6th Conference on Formal Power Series and Algebraic Combinatorics (New Brunswick, NJ, 1994). Discrete Math. \textbf{157} (1996), no. 1-3, 193--197.

\bibitem{Go-Ko-Ma} 
M. Goresky, R. Kottwitz, and R. MacPherson,
\emph{Equivariant cohomology, Koszul duality, and the localization theorem}, 
Invent. Math. \textbf{131} (1998), no. 1, 25--83. 
   
\bibitem{G-Ko-Ma06}
M. Goresky, R. Kottwitz, and R. MacPherson,
\emph{Purity of equivalued affine Springer fibers},
Represent. Theory \textbf{10} (2006), 130--146.

\bibitem{gu13}
M. Guay-Paquet,
\emph{A modular relation for the chromatic symmetric functions of (3+1)-free posets}, arXiv:1306.2400.

\bibitem{gu}
M. Guay-Paquet,
\emph{A second proof of the Shareshian-Wachs conjecture, by way of a new Hopf algebra}, arXiv:1601.05498.

\bibitem{ha-ho-ma}
M. Harada, T. Horiguchi, and M. Masuda, \emph{The equivariant cohomology rings of Peterson varieties in all Lie types}, 
Canad. Math. Bull. \textbf{58} (2015), no. 1, 80--90. 

\bibitem{ha-ho-ma-pa}   
M. Harada, T. Horiguchi, M. Masuda, and S. Park, 
\emph{The volume polynomial of regular semisimple Hessenberg varieties and the Gelfand-Zetlin polytope},
Proceedings of the Steklov Institute of Mathematics, \textbf{305} (2019), 318-–344. 

\bibitem{HHMPT}
M. Harada, T. Horiguchi, S. Murai, M. Precup, and J. Tymoczko, 
\emph{A filtration on the cohomology rings of regular nilpotent Hessenberg varieties},
in preparation.
   
\bibitem{ha-pr17}
M. Harada and M. Precup,
\emph{The cohomology of abelian Hessenberg varieties and the Stanley-Stembridge conjecture},
arXiv:1709.06736,
to be published in J. of Alg. Comb.

\bibitem{ha-pr18}
M. Harada and M. Precup,
\emph{Upper-triangular linear relations on multiplicities and the Stanley-Stembridge conjecture},
arXiv:1812.09503.

\bibitem{ha-ty}
M. Harada and J. Tymoczko, \emph{A positive Monk formula in the $S^1$-equivariant cohomology of type A Peterson varieties}, Proc. Lond. Math. Soc. (3) \textbf{103} (2011), no. 1, 40--72. 

\bibitem{Hori}
T. Horiguchi, 
\emph{The cohomology rings of regular nilpotent Hessenberg varieties and Schubert polynomials},
Proc. Japan Acad. Ser. A Math. Sci \textbf{94} (2018) 87--92.

\bibitem{Insko}
	E. Insko,
	\emph{Schubert calculus and the homology of the Peterson variety}, 
	Electron. J. Combin. \textbf{22} (2015), no.2, Paper 2.26, 12 pp.

\bibitem{InskoPrecup}
E. Insko and M. Precup,
\emph{The singular locus of semisimple Hessenberg varieties},
J. Algebra \textbf{521} (2019), 65--96. 

\bibitem{Insko-Tymoczko}
	E. Insko and J. Tymoczko,
	\emph{Intersection theory of the Peterson variety and certain singularities of Schubert varieties},
	Geom. Dedicata \textbf{180} (2016), 95--116. 

\bibitem{InskoTymoczkoWoo}
E. Insko, J. Tymoczko, and A. Woo
\emph{A formula for the cohomology and K-class of a regular Hessenberg variety}, 
J. Pure Appl. Algebra., DOI: https://doi.org/10.1016/j.jpaa.2019.106230.

\bibitem{InskoYong} 
E. Insko and A. Yong,  
\emph{Patch ideals and Peterson varieties}, 
Transform. Groups \textbf{17} (2012), no. 4, 1011--1036. 
   
\bibitem{Klyachko}
A. Klyachko,
\emph{Orbits of a maximal torus on a flag space}, 
Functional Anal. Appl. {\bf 19} (1985),  no. 2, 65--66. 

\bibitem{Ko}
B. Kostant, 
\emph{Flag manifold quantum cohomology, the Toda lattice, and the representation with highest weight $\rho$}, 
Selecta Math. (N.S.) \textbf{2} (1996), 43--91. 

\bibitem{mb-ty13}
A. Mbirika and J. Tymoczko, 
\emph{Generalizing Tanisaki's ideal via ideals of truncated symmetric functions}, 
J. Algebraic Combin. \textbf{37} (2013), no. 1, 167--199. 
  
\bibitem{OT} 
P. Orlik and H. Terao, 
\textit{Arrangements of hyperplanes},
Grundlehren der Mathematischen Wissenschaften, 
\textbf{300}, Springer-Verlag, Berlin, 1992.

\bibitem{precup13a}
M. Precup, \emph{Affine pavings of Hessenberg varieties for semisimple groups}, 
Sel. Math. New Series \textbf{19} (2013), 903--922. 
 
\bibitem{precup18}
M. Precup,
\emph{The Betti numbers of regular Hessenberg varieties are palindromic}, 
Transform. Groups 23 (2018), no. 2, 491--499. 

\bibitem{pr-ty17}
M. Precup and J. Tymoczko,
\emph{Hessenberg varieties of parabolic type},
arXiv:1701.04140.

\bibitem{proc90}
C. Procesi, \emph{The toric variety associated to Weyl chambers}, Mots, 153-161, Lang. Raison. Calc., Herm\'es, Paris, 1990. 

\bibitem{R}
K. Rietsch, 
\emph{Totally positive Toeplitz matrices and quantum cohomology of partial flag varieties}, 
J. Amer. Math. Soc. \textbf{16} (2003), 363--392.

\bibitem{Roe}
G. R\"ohrle, 
\emph{Arrangements of ideal type}, 
J. Algebra \textbf{484} (2017), 126--167. 

\bibitem{S2} 
K. Saito, 
\emph{Theory of logarithmic differential forms and logarithmic vector fields}, 
J. Fac. Sci. Univ. Tokyo Sect.IA Math. \textbf{27} (1980), 265--291.

\bibitem{sh-wa11}
J. Shareshian and M. L. Wachs, 
\emph{Chromatic quasisymmetric functions and Hessenberg varieties}, 
Configuration spaces, 433--460, CRM Series, \textbf{14}, Ed. Norm., Pisa, 2012. 

\bibitem{sh-wa14}
J. Shareshian and M. L. Wachs, 
\emph{Chromatic quasisymmetric functions}, 
Adv. Math. \textbf{295} (2016), 497--551.

\bibitem{SoloT}
L. Solomon and H. Terao, 
\emph{A formula for the characteristic polynomial of an arrangement}, 
Adv. in Math. 64 (1987), no.3, 305--325.
  
\bibitem{so-ty}  
E. Sommers and J. Tymoczko,
\emph{Exponents for B-stable ideals}, 
Trans. Amer. Math. Soc. \textbf{358} (2006), no. 8, 3493--3509. 

\bibitem{Spa}
N. Spaltenstein, N.
\emph{The fixed point set of a unipotent transformation on the flag manifold}, 
Nederl. Akad. Wetensch. Proc. Ser. A \textbf{79}, Indag. Math. \textbf{38} (1976), no. 5, 452--456. 

\bibitem{Spr1}
T. A. Springer, 
\emph{Trigonometric sums, Green functions of finite groups and representations of Weyl
groups}, 
Invent. Math. \textbf{36} (1976), 173--207.

\bibitem{Spr2}
T. A. Springer, 
\emph{A construction of representations of Weyl groups}, 
Invent. Math. \textbf{44} (1978), 279--293.
  
\bibitem{Stanley1989}
R. P. Stanley, 
\emph{Log-concave and unimodal sequences in algebra, combinatorics, and geometry. Graph theory and its applications: East and West (Jinan, 1986)}, 500--535, 
Ann. New York Acad. Sci., \textbf{576}, New York Acad. Sci., New York, 1989. 

\bibitem{Stanley1995}
R. P. Stanley, 
\emph{A symmetric function generalization of the chromatic polynomial of a graph},
Adv. Math. \textbf{111} (1995), no. 1, 166--194. 

\bibitem{st-ste}
R. P. Stanley and J. R. Stembridge, 
\emph{On immanants of Jacobi-Trudi matrices and permutations with restricted position}, 
J. Combin. Theory Ser. A \textbf{62} (1993), no. 2, 261--279. 

\bibitem{ste}
J. Stembridge, 
\emph{Eulerian numbers, tableaux, and the Betti numbers of a toric variety}, 
Discrete Math. \textbf{99} (1992), no. 1-3, 307--320. 

\bibitem{Teff11}
N. Teff, 
\emph{Representations on Hessenberg varieties and Young's rule}, 
23rd International Conference on Formal Power Series and Algebraic Combinatorics (FPSAC 2011), 903--914, Discrete Math. Theor. Comput. Sci. Proc., AO, Assoc. Discrete Math. Theor. Comput. Sci., Nancy, 2011.

\bibitem{Teff13}
N. Teff, 
\emph{A divided difference operator for the highest root Hessenberg variety}, 
25th International Conference on Formal Power Series and Algebraic Combinatorics (FPSAC 2013), 993--1004, 
Discrete Math. Theor. Comput. Sci. Proc., AS, Assoc. Discrete Math. Theor. Comput. Sci., Nancy, 2013. 

\bibitem{ty}
J. Tymoczko, 
\emph{Linear conditions imposed on flag varieties}, 
Amer. J. Math. \textbf{128} (2006), no. 6, 1587--1604. 

\bibitem{Ty2}
J. Tymoczko, 
\emph{Paving Hessenberg varieties by affines}, 
Selecta Math. (N.S.) \textbf{13} (2007), 353--367.  

\bibitem{tymo08}
J. Tymoczko, 
\emph{Permutation actions on equivariant cohomology of flag varieties}, 
Toric topology, 365--384, Contemp. Math., \textbf{460}, Amer. Math. Soc., Providence, RI, 2008. 

\end{thebibliography}
\end{document}